\documentclass[12pt]{article}
\usepackage[T1]{fontenc}
\usepackage{amsfonts}
\usepackage{amssymb}
\usepackage{makeidx}
\usepackage{amsmath}
\usepackage{amsthm}
\usepackage{graphics}
\usepackage{epsfig}
\usepackage{multirow}
\usepackage{lscape}
\usepackage{algorithmic}
\usepackage{algorithm}
\usepackage{slashbox}
\usepackage{color}
\usepackage{vmargin}
\setmarginsrb{25mm}{25mm}{-1mm}{-1cm}{0cm}{0cm}{0cm}{1.5cm}
\newcommand{\n}{\mathbb{N}}

\newcommand{\Sp}{\textrm{Sp}}

\newtheorem{propriete}{Property}

\newtheorem{proposition}{Proposition}

\newtheorem{lemme}{Lemma}

\newtheorem{theorem}{Theorem}

\newtheorem*{theorem_annexe*}{Theorem \ref{vpP-2}}

\newtheorem{cor}{Corollary}

\newenvironment{demo}{{\textbf{Proof}. }}{\begin{flushright}$\Box$\end{flushright}}

\newtheorem*{theorem*}{Theorem}

\title{The lollipop graph is determined by its spectrum}

\author{R. Boulet, B. Jouve\\
\small Institut de Mathématiques \\[-0.8ex]
\small Université de Toulouse et CNRS (UMR 5219)\\[-0.8ex]
\small \texttt{\{boulet,jouve\}@univ-tlse2.fr}}

\date{ Preprint submitted to \textit{The Electronic Journal of Combinatorics}, Feb 2008 }

\begin{document}

\maketitle

\begin{abstract}
An even (resp. odd) lollipop is  the coalescence of a cycle of even (resp. odd) length and a path with pendant vertex as distinguished vertex. It is known that the odd lollipop is determined by its spectrum and the question is asked by W. Haemers, X. Liu and Y. Zhang for the even lollipop. We revisit the proof for odd lollipop, generalize it for even lollipop and therefore answer the question. Our proof is essentially based on a method of counting  closed walks.
\end{abstract}

\definecolor{gris}{gray}{0.5}
\marginpar{\makebox[0cm][l]{\rotatebox{90}{\Huge \textcolor{gris}{Preprint submitted to \emph{The Electronic Journal of Combinatorics}} \normalsize}}}

\section{Introduction}

Let $G$ be a simple graph with $n$ vertices and $A$ its adjacency matrix, $Q_G(X)$ denotes its characteristic polynomial and $\lambda_1(G)\geq \lambda_2(G) \geq \cdots \geq \lambda_n(G)$ the associated eigenvalues; $\lambda_1(G)$ is the spectral radius of $G$. It is known that some informations about the graph structure can be deduced from these eigenvalues such as the number of edges or the length of the shortest odd cycle; but the reverse question \emph{Which graphs are determined by their spectrum ?} (asked, among others, in \cite{RePEc:dgr:kubcen:200266}) is far from being solved; some partial results exist  \cite{The_complement_of_the_path, Z_n, T-shape} which contribute to answer this question. 

Let us remind that the coalescence of two graphs $G_1$ with distinguished vertex $v_1$ and $G_2$ with distinguished vertex $v_2$, is formed by identifying vertices $v_1$ and $v_2$ that is, the vertices $v_1$ and $v_2$ are replaced by a single vertex $v$ adjacent to the same vertices in $G_1$ as $v_1$ and the same vertices in $G_2$ as $v_2$. If it is not necessary $v_1$ or $v_2$ may not be specified.

\begin{figure}[htbp]
\begin{center}
\includegraphics[scale=0.5]{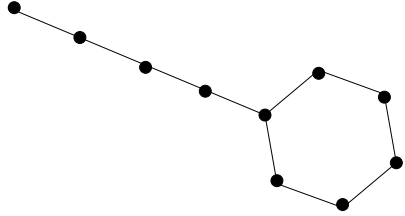}
\caption{Lollipop L(6,4)}
\label{lasso}
\end{center}
\end{figure}

A lollipop $L(p,k)$ is the coalescence of a cycle $C_p$ with $p\geq3$ vertices and a path $P_{k+1}$ with $k+1\geq 2$ vertices with one of its vertex of degree one as distinguished vertex, figure \ref{lasso} shows an example of a lollipop. The lollipop $L(p,0)$ is $C_p$. An even (resp. odd) lollipop has a cycle of even (resp. odd) length.
In this paper we shall show that the lollipop graph is determined by its spectrum, answering to an open question asked in \cite{lollipop, developments} for even lollipop. It is known \cite{lollipop} that the odd lollipop is determined by its spectrum, but the proof given in \cite{lollipop} cannot be generalized for even lollipops. We revisit here this proof in order to generalize it to even lollipops.

We describe in section \ref{section_basic} some basic results of spectral graph theory we shall use in the following of the paper. We also explain the method we  use to count closed walks in a graph and revisit two proofs of results about lollipops. The main section of the paper (section \ref{sec_pair}) shows that the even lollipop is determined by its spectrum; the proof is based on two points: connectivity of a graph cospectral with an even lollipop and existence of a $4$-cycle in a graph cospectral with a $L(4,k)$.


To fix notations, the disjoint union of two graphs $G$ and $H$ is noted $G\cup H$.

As defined in \cite{T-shape} a T-shape tree $S_{a,b,c}$ ($a,b,c>0$) is a tree with one and only one vertex $v$ of degree $3$ such that $S_{a,b,c}\backslash \{v\}=P_a\cup P_b\cup P_c$. We extend this notation for all $b,c\in\n$ by $S_{0,b,c}=P_{b+c+1}$.

\begin{figure}[htbp]
\begin{center}
\includegraphics[scale=0.5]{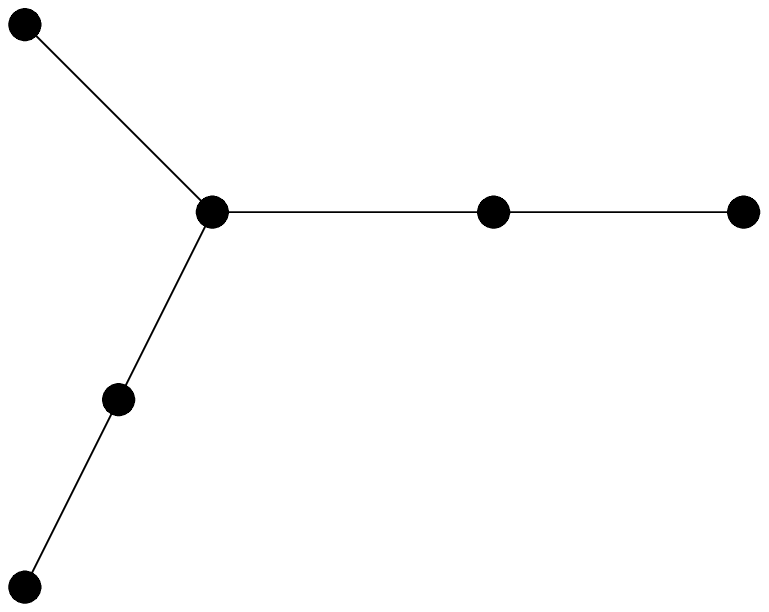}
\caption{$S_{1,2,2}$}
\label{Sabc}
\end{center}
\end{figure}

\begin{figure}[htbp]
\begin{center}
\includegraphics[scale=0.5]{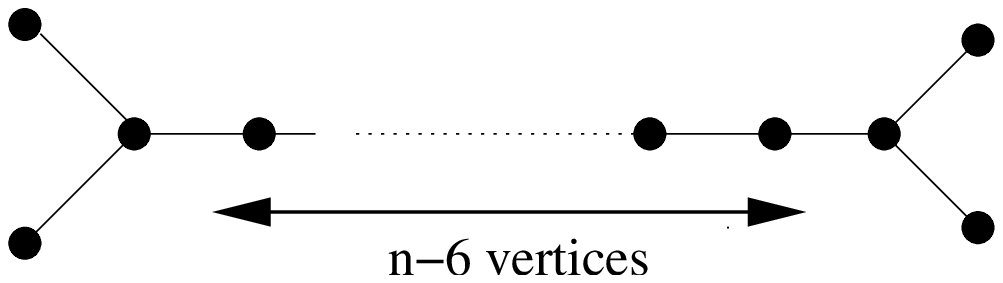}
\caption{$T_n$}
\label{Tn}
\end{center}
\end{figure}

By $S_{n-1}$ we denote the star with n vertices and by $T_n$ the tree with $n$ vertices drawn on figure \ref{Tn}.

Finally let $d(u,v)$ be the distance (the length of a shortest path) between two vertices $u$ and $v$ and $\delta(v)$ the degree of a vertex $v$.

\section{Basic results and revisited proofs}\label{section_basic}

\subsection{Counting the closed walks}\label{sec_marches}

It is a classical result that the number of closed walks of length $k\geq 2$ is $\sum_i\lambda_i^k$

We describe  here a method to count the number of closed walks of a given length within a graph.


Let $M$ be a graph, a \emph{$k$-covering closed walk} in $M$ is a closed walk of length $k$ in $M$ running through \emph{all} the edges at least once. Let $G$ be a graph, $M(G)$ denotes the set of all distinct subgraphs (not necessarily induced) of $G$ isomorphic to $M$ and $|M(G)|$ is the number of elements of $M(G)$. According to that point of view, $M$ may be called a motif (or a pattern).
The number of $k$-covering closed walks in a motif $M$ is denoted by $w_k(M)$ and we define the set $\mathcal{M}_k(G)=\{M,\ w_k(M)>0\}$ which is finite if $G$ is a finite graph. 

As a consequence, the number of closed walks of length $k$ in $G$ is: 
\begin{equation}
\sum_i \lambda_i^k=\sum_{M\in\mathcal{M}_k(G)} w_k(M)|M(G)|
\label{eq_closed_walk}
\end{equation}

In practice, there are at least two methods to determine $w_k(M)$: on one hand a combinatorial way which counts the number of covering closed walks of length $k$ in $M$, on the other hand an algebraic method which uses the following straightforward formula:
$$w_k(M)=\sum_{\lambda_i\in\Sp(M)}\lambda_i^k-\sum_{M'\in\mathcal{M}_k(M), \\ M'\neq M } w_k(M')|M'(M)|$$
where $\Sp(M)$ denotes the spectrum of the adjacency matrix of $M$.

%
%

Using equation (\ref{eq_closed_walk}) and table \ref{tab_marches} in appendix, we have the following proposition:

\begin{proposition}\label{sum6}
i) If $G$ is a graph without triangles and $C_5$ then:
\begin{eqnarray*}
\sum_i\lambda_i^6&=&12 |C_6(G)|+2|P_2(G)|+12|P_3(G)|+6|P_4(G)|+12|S_{1,1,1}(G)|\\
&&+48|C_4(G)|+12|L(4,1)(G)|
\end{eqnarray*}

ii)
If $G$ is a graph without $C_p$, $p\in\{3,5,6,7\}$ and of maximal degree $3$ then:
\begin{eqnarray*}
\sum_i\lambda_i^8&=&2|P_2(G)|+28|P_3(G)|+32|P_4(G)|+8|P_5(G)|+72|S_{1,1,1}(G)|+16|S_{1,1,2}(G)|\\
&&+264|C_4(G)|+112|L(4,1)(G)|+16|L(4,2)(G)|+16|C_8(G)|
\end{eqnarray*}
iii)
If $G$ is a graph without $C_p$, $p\in\{3,5,6,7,8,9\}$, of maximal degree $3$ and such that $\delta(u)=\delta(v)=3,\ u\neq v \Rightarrow d(u,v)>1$, then:
\begin{eqnarray*}
\sum_i\lambda_i^{10}&=&2|P_2(G)|+60|P_3(G)|+120|P_4(G)|+60|P_5(G)|+10|P_6(G)|+300|S_{1,1,1}(G)|\\
&&+140|S_{1,1,2}(G)|+20|S_{1,2,2}(G)|+20|S_{1,1,3}(G)|+1320|C_4(G)|\\
&&+840|L(4,1)(G)|+180|L(4,2)(G)|+20|L(4,3)(G)|+20|C_{10}(G)|
\end{eqnarray*}
\end{proposition}

In this paper we shall have to count all the $|M(G)|$, $M\in\mathcal{M}_i(G)$ of a given unicyclic graph $G$. For that aim we describe here the steps of the process we follow to count the $P_k(G)$ which are the only motifs hard to denombrate. Let $p$ be the length of the cycle of $G$.
\vspace{0.2cm}
\noindent\rule{10cm}{0.5mm}
\\\textbf{ALGORITHM to count $P_k(G)$:}
\begin{algorithmic}
\STATE  set $H=G$
\STATE set $|P_k(G)|=0$.
\WHILE{there exists a pendant vertex $u$ in $H$}
\STATE count the number $q$ of paths $P_k$ of $H$ containing $u$
\STATE let $|P_k(G)|=|P_k(G)|+q$
\STATE let $H=H\backslash \{u\}$
\ENDWHILE

\IF{$p\geq k$}
\STATE $|P_k(G)|=|P_k(G)|+p$
\ENDIF
\\\textbf{return} $|P_k(G)|$
\end{algorithmic}

\vspace{0.2cm}
\rule{10cm}{0.5mm}

\subsection{Known results}\label{sec_connus}


\begin{proposition}\label{sum4}\cite{spectra_of_unicyclic_graphs}
Let $G$ be a graph with $n$ vertices and $m$ edges and let $\lambda_i$ its associated eigenvalues. We have:  $\sum_i\lambda_i^4=8|C_4(G)|+2m+4|P_3(G)|$. 
Let $n_k$ be the number of vertices of degree $k$ in $G$, we have: $$\sum_i\lambda_i^4=8c_4+\sum_k kn_k+4\sum_{k\geq2} \frac{k(k-1)}{2}n_k$$
\end{proposition}
%


The following result relates the coefficients of the characteristic polynomial of a graph with structural properties of this graph:
\begin{theorem}\label{coef_poly}\cite{Doob}
Let $Q_G(X)=X^n+a_1X^{n-1}+a_2X^{n-2}+...+a_n$ be the characteristic polynomial of a graph $G$. We call an "elementary figure" the graph $P_2$ or the graphs $C_q,q>0$. We call a "basic figure" $U$ every graph all of whose components are elementary figures. Let $p(U)$ be the number of connected components of $U$ and $c(U)$ the number of cycles in $U$. We note $\mathcal{U}_i$ the set of basic figures with $i$ vertices. Then
$$a_i=\sum_{U\in\mathcal{U}_i}(-1)^{p(U)}2^{c(U)}\ ,\ i=1,2,...,n$$
\end{theorem}

It follows this theorem:
\begin{theorem}\label{coef_poly_impair}\cite{Doob}
Let $Q_G(X)=X^n+a_1X^{n-1}+a_2X^{n-2}+...+a_n$  be the characteristic polynomial of a graph $G$. The length of the shortest odd cycle in $G$ is given by the smallest odd index $p$ such that $a_p\neq0$ and the value of $a_p$ gives the number of $p$-cycles in $G$.
\end{theorem}

It ensues that a bipartite graph (\textit{ie} a graph with no odd cycles) cannot be cospectral with a non-bipartite graph.

%
The following result is useful at many time in the paper, for instance to find bounds on eigenvalues:
 
\begin{theorem}[Interlacing theorem]\cite{Godsil}\label{interlacement}
Let $G$ be a graph with $n$ vertices and associated eigenvalues $\lambda_1\geq \lambda_2\geq ... \geq \lambda_{n}$ and let $H$ be an induced subgraph of $G$ with $m$ vertices and associates eigenvalues $\mu_1\geq \mu_2\geq...\geq \mu_{m}$. Then for $i=1,...,m$,  $\lambda_{n-m+i}\leq \mu_{i}\leq\lambda_{i}$.
\end{theorem}

%
%
%

The next theorems give a way to compute the characteristic polynomial of a graph by deleting a vertex or an edge:
\begin{theorem}\label{isthme}\cite{Doob}
Let $G$ be a graph obtained by joining by an edge a vertex $x$ of a graph $G_1$ and a vertex $y$ of a graph $G_2$. Then
$$Q_G(X)=Q_{G_1}(X)Q_{G_2}(X)-Q_{G_1\backslash x}(X)Q_{G_2\backslash y}(X)$$
\end{theorem}


\begin{theorem}\label{del_vertex}\cite{Doob}
Let $G$ be a graph and $x$ a vertex of $G$, then:
$$Q_G(X)=XQ_{G\backslash x}(X)-\sum_{y\sim x}Q_{G\backslash\{x,y\}}(X)-2\sum_{C,\ x\in C}Q_{G\backslash C}(X)$$
where $y\sim x$ means that $yx$ is an edge of $G$ and the second sum is on the set of the cycles $C$ containing $x$.
\end{theorem}

\begin{theorem}\label{del_pendant}\cite{Doob}
Let $G$ be a graph and $x$ a pendant vertex of $G$. Then:
$$Q_G(X)=XQ_{G\backslash x}(X)-Q_{G\backslash {x,y}}(X)$$
where $y$ is the neighbor of $x$.
\end{theorem}

\begin{propriete}\label{rec_} We have the following equalities:\\
$Q_{C_p}(X)=XQ_{P_{p-1}}(X)-2Q_{P_{p-2}}(X)-2$\\
$Q_{P_p}(X)=XQ_{P_{p-1}}(X)-Q_{P_{p-2}}(X)$
\end{propriete}

\begin{demo}
A direct consequence of theorems \ref{isthme} and \ref{del_pendant}. 
\end{demo}

The following theorem relates the behavior of the spectral radius of a graph by subdividing an edge. An internal path of a graph $G$ is an elementary path $x_0x_1\cdots x_k$ (\textit{ie} $x_i\neq x_j$ for all $i\neq j$ but eventually $x_0=x_k$) of $G$ with $\delta(x_0)>2,\delta(x_k)>2,\delta(x_i)=2$ for all other $i$'s.

\begin{theorem}\label{homeo}\cite{On_the_largest_eigenvalue_of_some_homeomorphic_graphs,radii}
Let $xy$ be an edge of a connected graph $G$ not belonging to an internal path, then the spectral radius strictly increases by subdividing  $xy$.
\\Let $xy$ be an edge of a connected graph $G\neq T_n$  belonging to an internal path, then the spectral radius strictly decreases by subdividing $xy$.
\end{theorem}

%
%
%

\begin{theorem}\label{lower_bound}\cite{graffiti}
Let $G$ be a graph with maximal degree $\delta_M$, then $\lambda_1(G)\geq\sqrt{\delta_M}$
\end{theorem}

\begin{figure}[htbp]
\begin{center}
\includegraphics[scale=0.5]{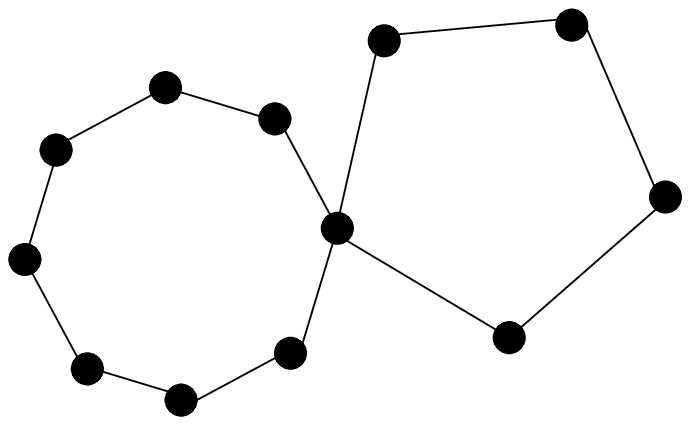}
\caption{B(8,5)}
\label{bouquet}
\end{center}
\end{figure}

Let $B(p,q)$ be the coalescence of two cycles $C_p$ and $C_q$ (see figure \ref{bouquet} for an example).

\begin{theorem}\label{vpB} \cite{On_the_largest_eigenvalue_of_some_homeomorphic_graphs}
For $p\geq3,\ q\geq3$, $\lambda_1(B(p,q))>\frac{4}{\sqrt{3}}>\sqrt{5}$
\end{theorem}

\subsection{Bounds on eigenvalues}

Theorem \ref{homeo} gives the following corollaries: 

\begin{cor}\label{aug_circ}
$\lambda_1(L(p,k))>\lambda_1(L(p+1,k))$
\end{cor}

\begin{cor}\label{aug_paths2}
$\lambda_1(L(p,k))<\lambda_1(L(p,k+1))$
\end{cor}
%

\begin{figure}[htbp]
\begin{center}
\includegraphics[scale=0.5]{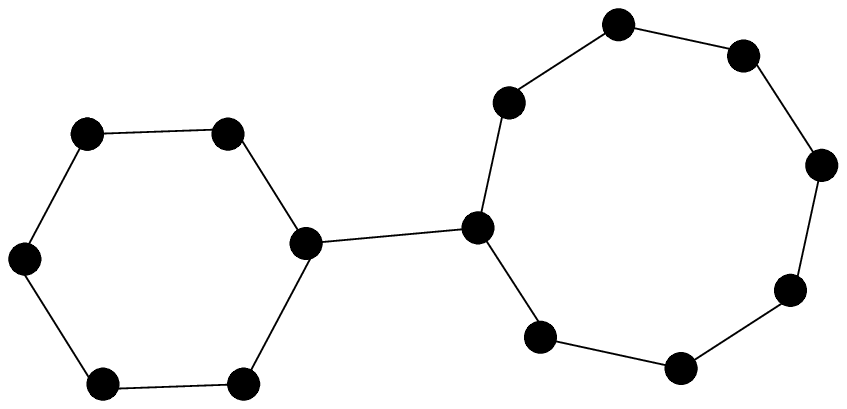}
\caption{H(6,8)}
\label{haltere}
\end{center}
\end{figure}

Given $p\geq3$, $q\geq3$, let $H(p,q)$ be the coalescence of $C_p$ and $L(q,1)$ with the pendant vertex as distinguished vertex (see figure \ref{haltere} for an example).

\begin{theorem}\label{vpH}
$\lambda_1(H(p,q))>\sqrt{5}$.
\end{theorem}

\begin{demo}
Without loss of generality we suppose that $p\geq q$. According to theorem \ref{homeo} we have $\lambda_1(H(p,q))\geq\lambda_1(H(p,p))$ so it is sufficient to prove the theorem for $H(p,p)$. As $\lim_{x\rightarrow +\infty}Q_{H(p,q)}(x)=+\infty$ it is sufficient to prove that $Q_{H(p,p)}(\sqrt{5})<0$
\\
Theorem \ref{isthme} gives:
$$Q_{H(p,p)}(X)=Q_{C_p}(X)Q_{C_p}(X)-Q_{P_{p-1}}(X)Q_{P_{p-1}}(X)$$
$$Q_{H(p,p)}(X)=[Q_{C_p}(X)]^2-[Q_{P_{p-1}}(X)]^2$$
and by property \ref{rec_} we have:
\begin{eqnarray*}
Q_{H(p,p)}(X)&=&[XQ_{P_{p-1}}(X)-2Q_{P_{p-2}}(X)-2]^2-[Q_{P_{p-1}}(X)]^2
\end{eqnarray*}

Let $(u_n)_{n\in\mathbb{N}}$ be the sequence defined by $u_n=Q_{P_n}(\sqrt{5})$. We have (property \ref{rec_}):
$u_n=\sqrt{5}u_{n-1}-u_{n-2}$. Since $u_1=\sqrt{2}$ and $u_2=4$ then $u_n=\beta_1^{n+1}-\beta_2^{n+1}$ where $\beta_1=\frac{\sqrt{5}+1}{2}$ and $beta_2=\frac{\sqrt{5}-1}{2}$.
%

\begin{eqnarray*}
Q_{H(p,p)}(\sqrt{5})&=&[\sqrt{5}u_{p-1}-2u_{p-2}-2]^2-(u_{p-1})^2\\
&=&[(\sqrt{5}+1)\beta_1^{p-1}-2][(\sqrt{5}-1)\beta_2^{p-1}-2]
\end{eqnarray*}

We have $[(\sqrt{5}+1)\beta_1^{p-1}-2]>0$ and $[(\sqrt{5}-1)\beta_2^{p-1}-2]<0$ so $Q_{H(p,p)}(\sqrt{5})<0$.

\end{demo}

\begin{theorem}\label{vp2}
For $k\neq0$ we have $\lambda_1(L(p,k))>2$ and $\lambda_{2}(L(p,k))<2$.
\end{theorem}

\begin{demo}

$\lambda_1(L(p,k))>2$:
the spectral radius of a cycle is $2$ and a cycle is an induced subgraph of $L(p,k)$ so by the interlacing theorem we have $\lambda_1(L(p,k))\geq2$. It remains to show that $\lambda_1(L(p,k))\neq2$. By theorem \ref{isthme} we have $Q_{L(p,k)}(2)=Q_{C_p}(2)Q_{P_k}(2)-Q_{P_{p-1}}(2)Q_{P_{k-1}}(2)=-Q_{P_{p-1}}(2)Q_{P_{k-1}}(2)\neq0$ (because the spectral radius of a path is strictly less than $2$).

$\lambda_{2}(L(p,k))<2$:
the path $P_{p+k-1}$ is an induced subgraph of $L(p,k)$ so by the interlacing theorem we have $\lambda_2(L(p,k))\leq \lambda_1(P_{p+k-1})<2$.
\end{demo}

\begin{theorem}\label{vp_max_L}
We have    $\lambda_1(L(p,k))<\sqrt{5}$.
\end{theorem}

\begin{demo}
By corollary \ref{aug_circ} we have $\lambda_1(L(p,k))\leq \lambda_1(L(3,k))$ so it is sufficient to prove the theorem for $p=3$.
For $k=0$, $\lambda_1(L(3,0))=2<\sqrt{5}$. We now assume that $k>0$.
Using theorem \ref{isthme} and $Q_{C_3}(X)=(X+1)^2(X-2)$ we have:
$$Q_{L(3,k)}(X)=(X+1)^2(X-2)Q_{P_k}(X)-(X-1)(X+1)Q_{P_{k-1}}(X)$$
and
$$Q_{L(3,k)}(\sqrt{5})=(2\sqrt{5}-2)Q_{P_k}(\sqrt{5})-4Q_{P_{k-1}}(\sqrt{5})$$
 Let us suppose that $Q_{L(3,k)}(\sqrt{5})>0$. 

We have $$Q_{L(3,k+1)}(\sqrt{5})=(2\sqrt{5}-2)Q_{P_{k+1}}(\sqrt{5})-4Q_{P_{k}}(\sqrt{5})$$ but $$Q_{P_{k+1}}(\sqrt{5})=\sqrt{5}Q_{P_{k}}(\sqrt{5})-Q_{P_{k-1}}(\sqrt{5})$$ so 
$$Q_{L(3,k+1)}(\sqrt{5})=\frac{2\sqrt{5}-2}{4}\big((2\sqrt{5}-2)Q_{P_k}(\sqrt{5})-4Q_{P_{k-1}}(\sqrt{5})\big)$$
and by induction on $k\geq 1$ we have $Q_{L(3,k+1)}(\sqrt{5})>0$.
\\
\\
Since the polynomial $Q_{L(3,k)}$ has one and only one root in $]2,+\infty[$ (theorem \ref{vp2}) then $Q_{L(3,k)}(2)<0$ and $Q_{L(3,k)}(\sqrt{5})>0$ implies that $\lambda_1(L(3,k))<\sqrt{5}$.
\end{demo}

\begin{theorem}\label{degre_max}
Let $G$ be a graph cospectral with $L(p,k)$, then $max\{\delta(v),\ v\in V(G) \}\leq4$.
\end{theorem}

\begin{demo}
A direct consequence of theorems \ref{vp_max_L} and \ref{lower_bound}.
\end{demo}

\begin{theorem}\label{forbid}
Let $G$ be a graph cospectral with a lollipop. Then, for $p\geq3$ and $q\geq3$, $C_p\cup C_q$ or $H(p,q)$ or $B(p,q)$ cannot be induced subgraphs of $G$
\end{theorem}

\begin{demo}
If $C_p\cup C_q$ is an induced subgraph of $G$ then as $\lambda_2(C_p\cup C_q)=2$ by interlacing theorem we get $\lambda_2(G)\geq2$, impossible by theorem \ref{vp2}.
\\$H(p,q)$ or $B(p,q)$ cannot be induced subgraphs of $G$ because $\lambda_1(G)<\sqrt{5}$ (theorem \ref{vp_max_L}) and $\lambda_1(H(p,q))>\sqrt{5}$ (theorem \ref{vpH}), $\lambda_1(B(p,q))>\sqrt{5}$ (theorem \ref{vpB}).
\end{demo}

\subsection{There are no cospectral non-isomorphic lollipops: revisited proof}\label{sec_pas_2_cosp}

In \cite{lollipop} it is proved that two cospectral lollipops are isomorphic. We revisit here this result in a shortest proof using closed walks.
\begin{theorem}\label{pas_2_cosp}
There are no cospectral non-isomorphic lollipops.
\end{theorem}

\begin{demo}
Let $L(p,k)$ and $L(p',k')$ with $n=p+k=p'+k'$ and $p<p'$ be two non isomorphic lollipops. To show that they have different spectra we  show that there are less closed walks of length $p$ in $L(p',k')$ than in $L(p,k)$.
\\ Let $e$ (resp. $e'$) be an edge of the cycle of $L(p,k)$ (resp. $L(p',k')$) incident to the vertex of degree $3$,  $\mathcal{W}$ (resp $\mathcal{W'}$)  the set of closed walks of length $p$ of $L(p,k)$ (resp. $L(p',k')$),  $\hat{\mathcal{W}}$ (resp $\hat{\mathcal{W'}}$)  the set of closed walks of length $p$ of $L(p,k)$ (resp. $L(p',k')$) not containing $e$ (resp. $e'$) and $\tilde{\mathcal{W}}$ (resp $\tilde{\mathcal{W'}}$)  the set of closed walks of length $p$ of $L(p,k)$ (resp. $L(p',k')$) containing $e$ (resp. $e'$).
\\We have: $|\mathcal{W}|=|\hat{\mathcal{W}}|+|\tilde{\mathcal{W}}|$ (resp. $|\mathcal{W'}|=|\hat{\mathcal{W'}}|+|\tilde{\mathcal{W'}}|$).
It's obvious that $|\hat{\mathcal{W}}|=|\hat{\mathcal{W'}}|$ because $L(p,k)\backslash\{e\}=L(p',k')\backslash\{e'\}=P_{n}$.
We are going to show that $|\tilde{\mathcal{W}}|<|\tilde{\mathcal{W'}}|$ by the following  equation:
$$|\tilde{\mathcal{W}}|=\sum_{M\in\mathcal{M}_p,\ e\in E(M)} w_p(M)|M(G)|$$
where $E(M)$ is the set of the edges of $M$.

We denote by $M^e$ a motif $M$ containing $e$. The motifs containing $e$ (resp $e'$) with at least one $p$-covering closed walk are exactly :
\\\textbullet\  the $P_i$'s for $2\leq i\leq\frac{p}{2}+1$ and we have  $|P_i^{e'}(L(p',k'))|\leq |P_i^e(L(p,k))|$.
\\\textbullet\  the $S_{a,b,c}$'s with $a+b+c\leq\frac{p}{2}$ and we have $|S_{a,b,c}^{e'}(L(p',k'))|\leq |S_{a,b,c}^e(L(p,k))|$.
\\\textbullet\  the $C_p$'s and $0=|C_p^{e'}(L(p',k'))|<|C_p^e(L(p,k))|=1$.

So, $|\tilde{\mathcal{W}}|<|\tilde{\mathcal{W'}}|$ and $|\mathcal{W}|<|\mathcal{W'}|$ which concludes the proof.
\end{demo}

\subsection{The odd lollipop is determined by its spectrum: revisited proof}\label{sec_impair}

We revisit here the proof that the odd lollipop is determined by its spectrum. 
 The aim of the proof is to determine the degree distribution. We already know that there are no vertices of degree greater or equal than 5 (theorem \ref{degre_max}). 
 
 \begin{lemme}
Let $G$ be a graph cospectral with $L(p,k)$, $p$ odd. Then $G$ has no isolated vertices.
 \end{lemme}
 
 \begin{demo}
We have to show that $0$ is not an eigenvalue of $L(p,k)$  that is the constant coefficient, $a_n$, of the characteristic polynomial of $L(p,k)$ is non-zero.
According to theorem \ref{coef_poly} we have:
 $$a_n=\sum_{U\in\mathcal{U}_n}(-1)^{p(U)}2^{c(U)}$$
But $|\mathcal{U}_n|=1$ :
\begin{itemize}
 \item if $k$ is odd, then $\mathcal{U}_n$ is the disjoint union of $\frac{p+k}{2}$ paths $P_2$, and $a_n=(-1)^{\frac{p+k}{2}}\neq0$.
 \item If $k$ is even then  $\mathcal{U}_n$ is the disjoint union of  $\frac{k}{2}$ paths $P_2$ and a cycle $C_p$, and $a_n=(-1)^{\frac{k}{2}+1}2\neq0$.
 \end{itemize}
 \end{demo}

\begin{lemme}\label{odd_4_cycle}
Let $G$ be a graph cospectral with $L(p,k)$, $p$ odd. Then there are no $4$-cycles in $G$.
\end{lemme}

\begin{demo}
Let us remark that an odd closed walk necessary runs through an odd cycle.
As $G$ and $L(p,k)$ have the same characteristic polynomial, according to theorem \ref{coef_poly_impair}, the length of the shortest odd cycle of $G$ is $p$ and there is only one such cycle,  so $\mathcal{M}_{p+2}(G)\subset\{C_p,L(p,1),C_{p+2}\}$.
Using equation (\ref{eq_closed_walk}) we have:
\begin{eqnarray}\label{eq1000}
\sum_{\lambda_i\in\Sp(G)}\lambda_i^{p+2}&=&w_{p+2}(C_p)|C_p(G)|+w_{p+2}(L(p,1))|L(p,1)(G)|+w_{p+2}(C_{p+2})|C_{p+2}(G)|\nonumber\\
&=&w_{p+2}(C_p)+(2p+4)|L(p,1)(G)|+(2p+4)|C_{p+2}(G)|
\end{eqnarray}
and
\begin{eqnarray}\label{eq2000}
\sum_{\lambda_i\in\Sp(L(p,k))}\lambda_i^{p+2}&=&w_{p+2}(C_p)+(2p+4)
\end{eqnarray}
If $|L(p,1)(G)|=0$ then $C_p$ or $C_p$ with (at least) a chord is a connected component of $G$. But the first case is impossible because $2$ is not an eigenvalue of $G$ and the second case is impossible because there are no odd cycles of length less than $p$ in $G$. So the equality of (\ref{eq1000}) and (\ref{eq2000}) implies that $|L(p,1)(G)|=1$ and $|C_{p+2}(G)|=0$.
If we suppose that there is a $4$-cycle in $G$, since $|L(p,1)(G)|=1$ the subgraph induced by $C_p$ and $C_4$ is $C_p\cup C_4$ or $H(p,4)$ but this is impossible by theorem \ref{forbid}.

\end{demo}

Now, we can prove the main theorem of this section: 

\begin{theorem}\label{th_odd}
Let $G$ be a graph cospectral with $L(p,k)$, $p$ odd. Then $G$ is isomorphic to $L(p,k)$.
\end{theorem}

\begin{demo}
Let $n_i$ be the number of vertices of degree $i$ for $i\in\{1,2,3,4\}$. We have $n=n_1+n_2+n_3+n_4$ and $2n=n_1+2n_2+3n_3+4n_4$ (the sum of the degrees is twice the number of edges), so $n_1=n_3+2n_4$.\\
Moreover by proposition \ref{sum4}, $\sum_{\lambda_i\in\Sp(G)} \lambda_i^4=8|C_4(G)|+2m+4(n_2+3n_3+6n_4)$ and by theorem \ref{odd_4_cycle}, $|C_4(G)|=0$. As $\sum_{\lambda_i\in\Sp(G)} \lambda_i^4=\sum_{\lambda_i\in\sp(L(p,k))} \lambda_i^4$ we get $n_2+3n_3+6n_4=n+1$ and then $1=-n_1+2n_3+5n_4$.\\
So we have $1=n_3+3n_4$ and then $n_4=0$, $n_3=1$,  $n_1=1$, $n_2=n-2$.

As the sum of the degrees of a graph is even, the vertex of degree $1$ and the vertex of degree $3$ belongs to the same connected component. If $G$ is not connected there is a 2-regular connected component (\textit{ie} a cycle) which is  impossible ($2$ is not an eigenvalue of $G$). 
As a result, $G$ is a connected graph with degree distribution equal to $(1,2,2,2,...,2,2,3)$, so $G$ is a lollipop and, by theorem \ref{pas_2_cosp}, $G$ is isomorphic to $L(p,k)$.
\end{demo}


\section{The even lollipop is determined by its spectrum.}\label{sec_pair}

Following the same method as the one used for the odd case, to prove that the even lollipop is determined by its spectrum we  show that a graph cospectral with an even lollipop:
\begin{itemize}
\item   is connected (and then it contains no isolated vertices). 
\item has a $4$-cycle if and only if it is cospectral with a $L(4,k)$. 
\end{itemize}
For the second point the difficulty is to prove that a graph cospectral with a $L(4,k)$ has a $4$-cycle.

To lighten the section some technical proofs have been detailed in appendix.


\subsection{Connectivity}\label{sec_isole}

Using results of section \ref{sec_connus} we easily obtain the following property:
\begin{propriete}\label{P(2)}
$\forall a,b,c\in\n$, $Q_{C_p}(2)=0$, $Q_{P_k}(2)=k+1$,  $Q_{S_{a,b,c}}(2)=a+b+c+2-abc$,  $Q_{S_{1,1,a}}(2)=4$
\end{propriete}

The following theorem gives a better bound than the theorem \ref{vp_max_L} on spectral radius of a lollipop $L(p,k)$ when $p\geq4$.

\begin{theorem}\label{2.17}
i) Let $G$ be a graph cospectral with $L(p,k)$ with $p\geq6$, then $\lambda_1(G)<2.17$.
\\ii) Let $G$  be a graph cospectral with $L(4,k)$, then $\lambda_1(G)<\sqrt{2+2\sqrt{2}}$.
\end{theorem}

\begin{demo}
Just follow the proof of theorem \ref{vp_max_L} \textit{mutatis mutandis}.
\end{demo}

\begin{figure}[htbp]
\begin{center}
\includegraphics[scale=0.5]{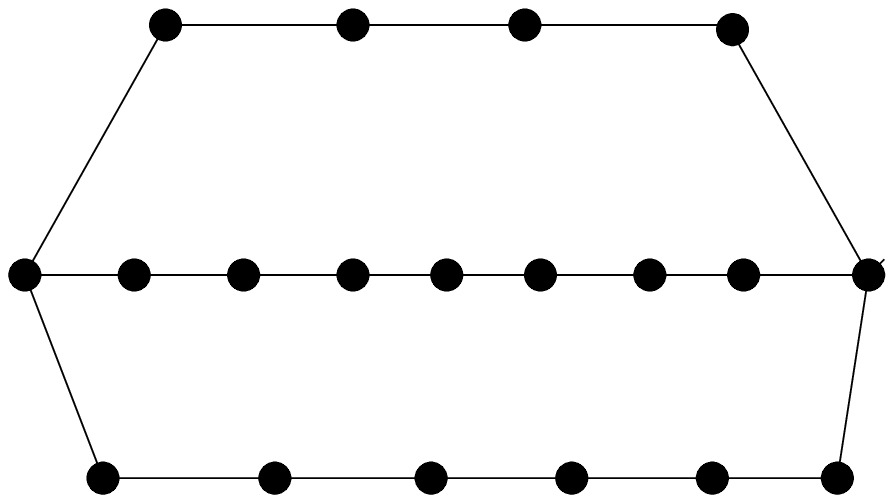}
\caption{P(4,7,6)}
\label{P3}
\end{center}
\end{figure}

Let $P(p_1,p_2,p_3)$ be the graph obtained by identifying the three pendant vertices of $S_{p_1+1,p_2+1,p_3+1}$ (an example is given in  figure \ref{P3}).

\begin{theorem}\label{vpP-2}
The graph $P(p_1,p_2,p_3)$ cannot be an induced subgraph of a graph $G$ cospectral with an even lollipop.
\end{theorem}

\begin{demo}
Sketch of the proof : 

We first show that for some values of $p_1,p_2$ and $p_3$ we have $\lambda_1(P(p_1,p_2,p_3))>\sqrt{2+2\sqrt{2}}$ and in these cases $P(p_1,p_2,p_3)$ cannot be an induced subgraph of G.

For the others cases we compute $Q_{P(p_1,p_2,p_3)}(2)$. 
\begin{itemize}
\item if  $Q_{P(p_1,p_2,p_3)}(2)\geq 0$ then $P(p_1,p_2,p_3)$ and \textit{a fortiori} $G$ (interlacing theorem) possesses two eigenvalues greater than $2$ which contradicts that $G$ is cospectral with a lollipop (theorem \ref{vp2}) . 

\item if $Q_{P(p_1,p_2,p_3)}(2)<0$ then we show that $P(p_1,p_2,p_3)$ cannot be a connected component of $G$ so there is a vertex $x$ not in $P(p_1,p_2,p_3)$ adjacent to a vertex $y$ of   $P(p_1,p_2,p_3)$ and we prove that this graph  so constructed cannot be an induced subgraph of $G$. 
\end{itemize}

A detailed proof is given in appendix \ref{annexe}.

\end{demo}

\begin{theorem}
Let $G$ be a graph cospectral with an even lollipop. Then $G$ is connected.
\end{theorem}

\begin{demo}
The graph $G$ has as many edges as vertices, so if $G$ is not connected, it possesses at least two cycles. The subgraph induced by the two cycles of minimal length is $C_a\cup C_b$, $B(a,b)$, $H(a,b)$ or $P(p_1,p_2,p_3)$ but this is impossible (theorems \ref{forbid} and \ref{vpP-2}).
\end{demo}

\begin{cor}\label{unicyclique}
A graph cospectral with an even lollipop is unicyclic.
\end{cor}
%

\subsection{The even lollipop $L(p,k)$, $p\geq 6$, is determined by its spectrum}

Let $G$ be a graph cospectral with an even lollipop $L(p,k)$, $p\geq6$.
In order to copy the proof of theorem \ref{th_odd} concerning the odd lollipop we have to show that $|C_4(G)|=0$ ($G$ does not have a $4$-cycle), this is the aim of the following proposition.
\begin{proposition}\label{L6knot4cycle}
A graph cospectral with an even lollipop $L(p,k)$, $p\geq 6$ does not have a $4$-cycle.
\end{proposition}

\begin{demo}
Let $G$ be a graph cospectral with an even lollipop $L(p,k)$, $p\geq6$ and suppose that $G$ has a $4$-cycle. As $G$ is connected, unicyclic and has at least $6$ vertices then one of the graph drawn in figure \ref{4} is an induced subgraph of $G$ and we check that the spectral radius of theses graphs is greater than $2.17$.

\begin{figure}[htbp]
\begin{center}
\includegraphics[scale=0.5]{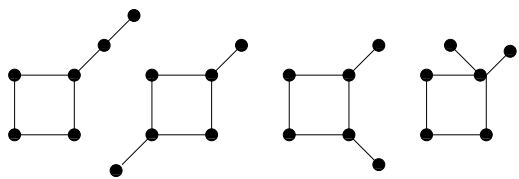}
\caption{Unicyclic graphs with six vertices and having a 4-cycle }
\label{4}
\end{center}
\end{figure}

This contradicts theorem \ref{2.17} .

\end{demo}

We can now state:
\begin{theorem}
The even lollipop $L(p,k)$, $p\geq 6$, is determined by its spectrum.
\end{theorem}

\subsection{The even lollipop $L(4,k)$  is determined by its spectrum} \label{sec_4cycle}

Let $G$ be a graph cospectral with $L(4,k)$, the main point is to show that the converse implication of previous proposition \ref{L6knot4cycle} holds, that is $G$ has a $4$-cycle. The key theorem of this part requires to study the cospectrality of some classes of unicyclic graphs with a lollipop $L(4,k)$, this is done in the sections \ref{sec_G1}, \ref{sec_G2}, \ref{sec_G3}, \ref{sec_G4}.

\subsubsection{Our toolbox: some results on $L(4,k)$}

In the following we are going to prove that $L(4,k)$ is not cospectral with some unicyclic graphs. For that purpose we use several tools detailed in this section: counting closed walks of length $6$, $8$ or $10$, evaluating the characteristic polynomial in $1$ or $2$, using the fact that a lollipop has only one eigenvalue greater than $2$.

\begin{proposition}\label{cor0}
i) 
For $L(4,k)$, $k>1$ we have: $$\sum_i\lambda_i^6=20n+96$$
\\
ii)
For $L(4,k)$, $k>2$ we have: $$\sum_i\lambda_i^8=70n+596$$
\\
iii)
For $L(4,k)$, $k>3$ we have $$\sum_i \lambda_i^{10}=252n+3360$$
%
%
%
%
%
%
\end{proposition}
 \begin{demo}
 Counting closed walks, we check that
\\$i)$ For $k>1$,  $|P_2(L(4,k))|=n$, $|P_3(L(4,k))|=n+1$, $|P_4(L(4,k))|=n+2$, $|C_4(L(4,k))|=1$, $|L(4,1)(L(4,k))|=1$.
\\ $ii)$ Moreover, for $k>2$, $|P_5(L(4,k))|=n-1$, $|S_{1,1,1}(L(4,k))|=3$, $|S_{1,1,2}(L(4,k))|=3$, $|L(4,2)(L(4,k))|=1$.
\\ $iii)$ Moreover, for $k>3$, $|P_6(L(4,k))|=n-2$, $|S_{1,2,2}(L(4,k))|=2$, $|S_{1,1,3}(L(4,k))|=1$, $|L(4,2)(L(4,k))|=1$, $|L(4,3)(L(4,k))|=1$
%
\\and apply proposition \ref{sum6}.
\end{demo}

\begin{propriete}\label{Pp(1)}
We have $Q_{P_p}(1)=Q_{P_{\overline{p}}}(1)$ and $Q_{C_p}(1)=Q_{C_{\overline{p}}}(1)$ where $\overline{p}$ is $p$ modulo $6$ and: 
$$\begin{array}{cc}

Q_{P_{\overline{0}}}(1)=1 & Q_{C_{\overline{0}}}(1)=0\\
Q_{P_{\overline{1}}}(1)=1 & Q_{C_{\overline{1}}}(1)=-1\\
Q_{P_{\overline{2}}}(1)=0 & Q_{C_{\overline{2}}}(1)=-3\\
Q_{P_{\overline{3}}}(1)=-1 & Q_{C_{\overline{3}}}(1)=-4\\
Q_{P_{\overline{4}}}(1)=-1 & Q_{C_{\overline{4}}}(1)=-3\\
Q_{P_{\overline{5}}}(1)=0 & Q_{C_{\overline{5}}}(1)=-1\\
\end{array}$$
\end{propriete}
%
%

\begin{demo}
According to property  \ref{rec_}, $Q_{P_p}(1)=Q_{P_{p-1}}(1)-Q_{P_{p-2}}(1)=-Q_{P_{p-3}}(1)=Q_{P_{p-6}}(1)$ and $Q_{C_p}(1)=Q_{P_{p-1}}(1)-2Q_{P_{p-2}}(1)-2$. Then we can easily compute $Q_{P_i}$ and $Q_{C_i}$ for $0\leq i\leq 5$.
\end{demo}

\begin{propriete}\label{P(0)}
We have:
$$Q_{P_k}(0)=\left\{\begin{array}{c} (-1)^{\frac{k}{2}}\ \textrm{ if $k$ is even}\\ 0 \textrm{ if $k$ is odd}\end{array}\right.$$
and if $k$ is odd we have $R(0)=(-1)^{\frac{k-1}{2}}\frac{k+1}{2}$ where $R(X)=\frac{Q_{P_k}(X)}{X}$.
\end{propriete}

\begin{demo}
Proofs by induction with the relation $Q_{P_k}(X)=XQ_{P_{k-1}}(X)-Q_{P_{k-2}}(X)$.
\end{demo}

\begin{proposition}\label{Pl4(1)}
We have:
$$Q_{L(4,k)}(1)=\left\{\begin{array}{l} 
1\ \textrm{ if } n\equiv0[6] \\
3\ \textrm{ if } n\equiv1[6]\\
2\ \textrm{ if } n\equiv2[6] \\
-1\ \textrm{ if } n\equiv3[6] \\
-3\ \textrm{ if } n\equiv4[6] \\
-2\ \textrm{ if } n\equiv5[6] \\
\end{array}\right.$$
\end{proposition}

\begin{demo}
Theorem \ref{isthme} gives $Q_{L(4,k)}(X)=Q_{C_4}(X)Q_{P_k}(X)-Q_{P_3}(X)Q_{P_{k-1}}(X)$ so $Q_{L(4,k)}(1)=-3Q_{P_k}(1)+Q_{P_{k-1}}(1)$ and we conclude with property \ref{Pp(1)}.
\end{demo}

\begin{proposition}\label{Pl4(2)}
$Q_{L(4,k)}(2)=-4n+16$.
\end{proposition}

\begin{demo}
$Q_{L(4,k)}(X)=Q_{C_4}(X)Q_{P_k}(X)-Q_{P_3}(X)Q_{P_{k-1}}(X)$ and with property \ref{P(2)} we have $Q_{L(4,k)}(2)=-4k=-4n+16$.
\end{demo}

\textbf{Remark} : This proposition can be generalized for all lollipops : $Q_{L(p,k)}(2)=-pk$.

\begin{proposition}\label{Pl4(0)}
If $n=4+k$ is even then $0$ is an eigenvalue of $L(4,k)$ with multiplicity $2$ and $R(0)=(-1)^{\frac{k}{2}+1}n$ where $R(X)=\frac{Q_{L(4,k)}(X)}{X^2}$.
\end{proposition}

\begin{demo}
Since $Q_{L(4,k)}(X)=Q_{C_4}(X)Q_{P_k}(X)-P_3(X)Q_{P_{k-1}}(X)$ we have $R(X)=(X^2-4)Q_{P_k}(X)-(X^2-2)\frac{Q_{P_{k-1}}(X)}{X}$ and property \ref{P(0)} gives the result.
\end{demo}
\subsubsection{Unicyclic graphs with exactly three vertices of maximal degree $3$ whose only one belongs to the cycle}\label{sec_G1}

\begin{figure}[htbp]
\begin{center}
\includegraphics[scale=0.5]{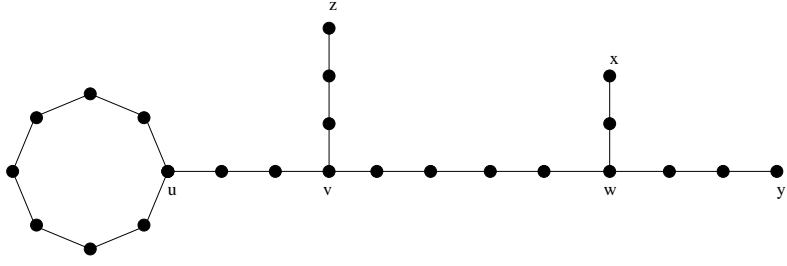}
\caption{A graph $G\in\mathcal{G}_1$}
\label{G1}
\end{center}
\end{figure}

Let $T$ be a tree with exactly two vertices of maximal degree $3$. Let $\mathcal{G}_1$ be  the set of  the coalescences of T with a pendant vertex as distinguished vertex and a cycle $C_p$, $p\geq 6$.  
In the following we assume that the vertex of degree $3$ belonging to the cycle is denoted by $u$ and $v$, $w$ are the other two vertices of degree $3$ such that $v$ is between $u$ and $w$; $x,y,z$ are the pendant vertices of $G$ such that $d(z,v)<d(z,w)$ and $d(x,w)\leq d(y,w)$. An example is given in figure \ref{G1}.

The aim of this section is to show the following theorem whose proof is summed up in table \ref{Proof1}:

\begin{theorem}\label{L4pascospG1}
The lollipop $L(4,k)$ cannot be cospectral with a graph $G\in\mathcal{G}_1$.
\end{theorem}

As $L(4,k)$ cannot be cospectral with a non-bipartite graph we suppose in the following that a graph $G\in\mathcal{G}_1$ is bipartite (the length of the cycle is even).


\begin{table}[h]
\small
\begin{tabular}{|| c|c|c|c|c|c    ||c|c||}
\hline
\multicolumn{6}{|c||}{Graph} & Tool & Prop. \\

\hline
\hline
\multicolumn{6}{||c||}{$p=6$} & $\sum \lambda_i^6$ & \ref{sum6part2} $v)$ \\

\hline

\multirow{12}{*}{ $p\geq 8$ } & \multirow{3}{*}{$\begin{array}{c}d(u,v)\\=1\end{array}$} & \multicolumn{4}{|c||}{$d(v,w)=1$} &  $\sum \lambda_i^6$ &  \ref{sum6part2} $i)$  \\

\cline{3-8}

&&\multirow{2}{*}{$\begin{array}{c} d(v,w)\\>1\end{array}$} & \multicolumn{3}{|c||}{\begin{tabular}{c}$d(v,z)>1$ or $d(x,w)>1$\\ or $d(y,w)>1$\end{tabular}} &  $\sum \lambda_i^6$ & \ref{sum6part2} $ii)$  \\

\cline{4-8}

&&& \multicolumn{3}{|c||}{\begin{tabular}{c}$d(v,z)=1$ and $d(x,w)=1$\\ and $d(y,w)=1$\end{tabular}} &  $Q_G(1)$ &  \ref{duv=1part2}   \\
\cline{2-8}

 & \multirow{8}{*}{$\begin{array}{c} d(u,v)\\=2\end{array}$} 
&  \multirow{6}{*}{$\begin{array}{c} d(x,w)\\=1\end{array}$} & \multicolumn{3}{|c||}{ $d(y,w)=1$ } & $Q_G(2)$ &  \ref{en2} \\
\cline{4-8}
&&                             & \multirow{3}{*}{$\begin{array}{c}d(y,w)\\=2\end{array}$} & \multirow{2}{*}{ $d(v,z)>1$ } & $d(v,w)=1$ &$\sum\lambda_i^6$ &  \ref{sum6part2} $iv)$ \\
\cline{6-8}
&&                               &                             &                             &   $d(v,w)>1$ & $\sum\lambda_i^6$ &  \ref{sum6part2} $iii)$ \\
\cline{5-8}

&&                              &                              & \multirow{3}{*}{$d(v,z)=1$} & $d(v,w)=1$ & $\sum\lambda_i^6$ &  \ref{sum6part2} $iv)$ \\
\cline{6-8}
&&                               &                             &                             & $2\leq d(v,w)\leq3$ & $\sum\lambda_i^8$ &  \ref{VT<=9} \\
\cline{6-8}
&&                               &                             &                             &   $d(v,w)\geq4$ & $\lambda_2\geq2$ &  \ref{T_nsub} $iii)$ \\
\cline{4-8}
&&                               & \multirow{2}{*}{$\begin{array}{c}d(y,w)\\>2\end{array}$} & \multicolumn{2}{|c||}{ $d(v,w)=1$ } & $\sum\lambda_i^6$ &  \ref{sum6part2}, $iv)$ \\
\cline{5-8}
&&                               &                              & \multicolumn{2}{|c||}{ $d(v,w)>1$ } &  $\lambda_2\geq2$ &  \ref{T_nsub} $ii)$ \\
\cline{3-8}

&&  \multirow{2}{*}{$\begin{array}{c} d(x,w)\\>1\end{array}$} & \multicolumn{3}{|c||}{ $d(v,w)=1$ } & $\sum\lambda_i^6$ &  \ref{sum6part2} $iv)$ \\
\cline{4-8}
&&                               & \multicolumn{3}{|c||}{ $d(v,w)>1$ } & $\sum\lambda_i^6$ & \ref{sum6part2} $iii)$ \\
\cline{2-8}

&\multicolumn{5}{|c||}{ $d(u,v)>2$ } & $\lambda_2\geq2$ &  \ref{T_nsub} $i)$ \\
\hline
\hline
\end{tabular}
\normalsize
\caption{Proof of theorem \ref{L4pascospG1} using a case disjunction over the possibilities for the values of $d$.}
\label{Proof1}
\end{table}


\begin{proposition}\label{T_nsub}
Let $G\in\mathcal{G}_1$. If one of the following properties is true:
\\i) $d(u,v)>2$ 
\\ii)  $d(u,v)=2,\ d(v,w)>1$ and $d(y,w)>2$
\\iii)  $d(u,v)=2,\ d(v,w)\geq4$, $d(y,w)\geq 2$
\\ then $G$ is not cospectral with a lollipop.
\end{proposition}

\begin{demo}

Let $p$ be the length of the cycle of $G$. If one of these properties is true then $G$ possesses an induced subgraph with twice the eigenvalue $2$. By the interlacing theorem it cannot be cospectral with a lollipop (theorem \ref{vp2}). 

This subgraph is $C_p\cup T_r$ (for an $r\in\n$) in the case $i)$, $C_p\cup S_{1,3,3}$ in $ii)$ and $C_p\cup S_{1,2,5}$ in $iii)$.
\end{demo}

\begin{proposition}\label{sum6part2}
Let $G\in\mathcal{G}_1$. If one of the following properties is true: 
\\i) $d(u,v)=1$, $d(v,w)=1$,
\\ii) $d(u,v)=1$ and $d(v,w)>1$  and ( $d(v,z)>1$ or $d(x,w)>1$ or $d(y,w)>1$),
\\iii) $d(u,v)>1$  and $d(v,w)>1$ and ( ($d(v,z)>1$ and $d(y,w)>1$) or $d(x,w)>1$),
\\iv) $d(v,w)=1$ and ( $d(v,z)>1$ or $d(y,w)>1$ or $d(x,w)>1$),
\\v) $p=6$.
\\then
$$\sum_{\lambda_i\in\Sp(G)}\lambda_i^6>20n+96$$
and $G$ cannot be cospectral with $L(4,k)$.
\end{proposition}

\begin{demo}
For the cases from $i)$ to $iv)$ we have $|P_2(G)|=n,\ |P_3(G)|=n+3,\ |S_{1,1,1}(G)|=3,\ |P_4(G)|>n+4$ and apply proposition \ref{sum6}.
\\For the case $v)$ we have  $|P_2(G)|=n,\ |P_3(G)|=n+3,\ |S_{1,1,1}(G)|=3,\ |P_4(G)|>n+2,\ |C_6(G)|=1$ and apply proposition \ref{sum6}.
\end{demo}

%

\begin{proposition}\label{duv=1part2}
Let $G\in\mathcal{G}_1$ such that $d(u,v)=1$ and $d(w,x)=d(w,y)=d(v,z)=1$. Then $G$ cannot be cospectral with $L(4,k)$.
\end{proposition}

\begin{demo}
Let $G\in\mathcal{G}_1$, with $n=p+q$ vertices where $p$ is the length of the cycle. 
We have:
 \begin{eqnarray*}
Q_G(X)&=&Q_{C_p}(X)Q_{S_{1,1,q-3}}(X)-XQ_{P_{p-1}}(X)Q_{S_{1,1,q-5}}(X)\\
&=&XQ_{C_p}(X)(Q_{P_{q-1}}(X)-Q_{P_{q-3}}(X))-X^2Q_{P_{p-1}}(X)(Q_{P_{q-3}}(X)-Q_{P_{q-5}}(X))
\end{eqnarray*}

Using property \ref{Pp(1)} we compute $Q_G(1)$, the result depends on $\overline{p}$ and $\overline{q}$ which are $p$ and $q$ modulo $6$ and are summed up into the following table:

$$\begin{array}{||c||c|c|c|c|c|c||}
\hline
\textrm{\backslashbox{$\overline{p}$}{$\overline{q}$}}&\overline{0}&\overline{1}&\overline{2}&\overline{3}&\overline{4}&\overline{5}\\
\hline
\hline
\overline{0}&	0	&	0   &	0	  &	0	&	0	&0        \\ 	
\overline{2}& -1	 &	-5   &	-4	  &	1	&	5	&4        \\ 	
\overline{4}& -5	 &	-7   &	-2	  &	5	&	7	&2         \\ 	
\hline
\hline
\end{array}$$

Comparing this results with proposition \ref{Pl4(1)} ( $\overline{n}=\overline{p}+\overline{q}$) we conclude that $G$ cannot be cospectral with $L(4,k)$.
\end{demo}


%
%

\begin{proposition}\label{VT<=9}
Let $G\in\mathcal{G}_1$ such that $p\geq8$,  $d(u,v)=2$, $2\leq d(v,w)\leq 3$, $d(y,w)=2$, $d(v,z)=1$, $d(x,w)=1$. Then $G$ cannot be cospectral with $L(4,k)$.
\end{proposition}

\begin{demo}
We have $|P_2(G)|=n$, $|P_3(G)|=n+3$, $|P_4(G)|=n+4$, $|S_{1,1,2}(G)|=7$, $|P_5(G)|=n+6 \textrm{ if } d(v,w)=2$ and $|P_5(G)|=n+5 \textrm{ if } d(v,w)\geq 3$ and by proposition \ref{sum6}: 
 
 $$\sum \lambda_i^8=\left\{\begin{array}{l} 
70n+588+16|C_8(G)|\ \textrm{ if }  d(v,w)=2\\
70n+580+16|C_8(G)|\ \textrm{ if }  d(v,w)= 3\\
\end{array}\right.$$
 
 \begin{itemize}
 \item If $d(v,w)=2$ then, by proposition \ref{cor0}, $G$ cannot be cospectral with $L(4,k)$.
 
 \item If $d(v,w)=3$ then, by proposition  \ref{cor0}, $G$ is cospectral with $L(4,k)$ only if $p=8$. We then check that such a graph $G$ (drawn on figure \ref{except}) is not cospectral with $L(4,13)$ by comparing spectral radii (see tables \ref{spect_L4k} and \ref{spect_special} in appendix).
 
 \begin{figure}[h]
 \begin{center}
\includegraphics[scale=0.5]{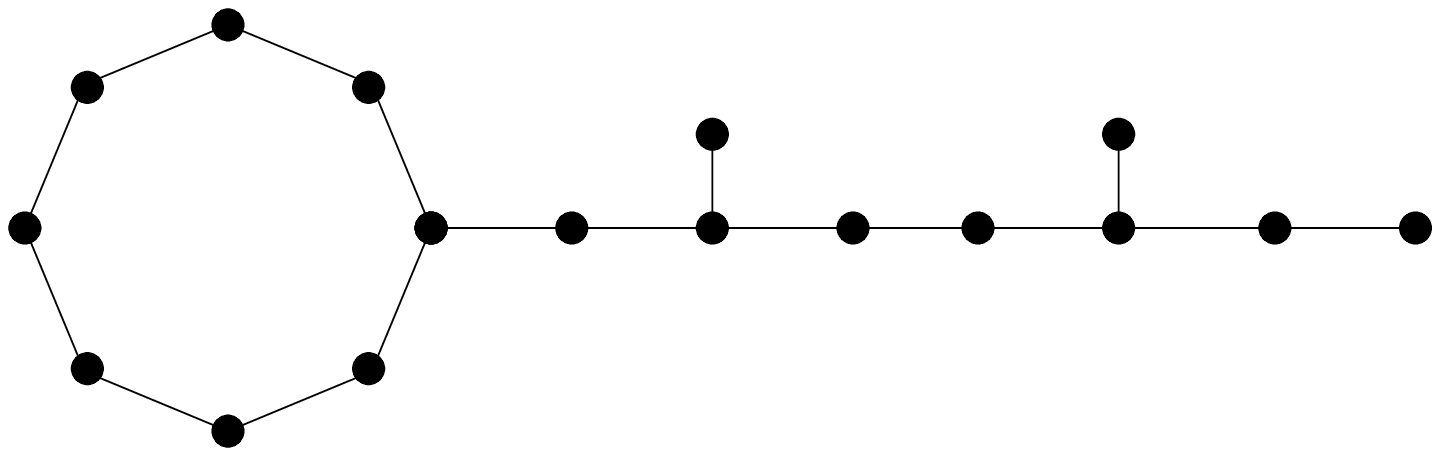}
\caption{}
\label{except}
\end{center}
\end{figure}

 \end{itemize}

\end{demo}

\begin{proposition}\label{en2}
Let $G\in\mathcal{G}_1$  such that $d(u,v)=2$,  $d(x,w)=d(y,w)=1$. Then $Q_G(2)=-4p$ and $G$ cannot be cospectral with a lollipop $L(4,k)$
\end{proposition}

\begin{demo}
Set $b=d(v,w)$ and $a=d(z,v)$, using theorems \ref{isthme} and \ref{del_pendant} we have
$$Q_G(X)=Q_{C_p}(X)Q_T(X)-Q_{P_{p-1}}(X)Q_{S_{1,1,a+b-1}}(X)$$
(where $T$ is a tree) and using property \ref{P(2)} we get $Q_G(2)=0-p\times4=-4n+4(n-p)$. As $n-p>4$, proposition \ref{Pl4(2)} implies that $G$ cannot be cospectral with a lollipop $L(4,k)$.
\end{demo}
%
%

%
%
\subsubsection{Unicyclic graphs with exactly three vertices of maximum degree $3$ whose exactly two belongs to the cycle.}\label{sec_G2}

\begin{figure}[htbp]
\begin{center}
\includegraphics[scale=0.5]{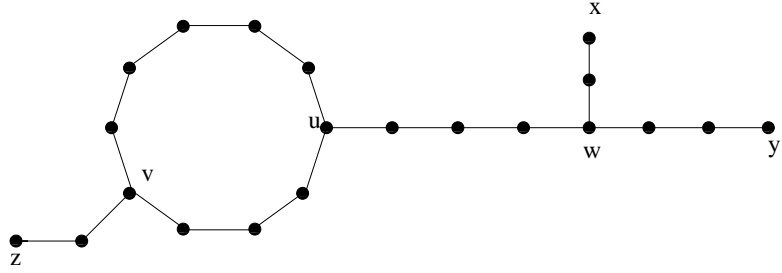}
\caption{A graph $G\in\mathcal{G}_2$}
\label{G2}
\end{center}
\end{figure}

Let $T$ be a tree with exactly one vertex $w$ of maximum degree $3$  and $L(p,k)$, $p\geq6$, a lollipop  (the vertex of degree $3$ is denoted by $v$ and the pendant vertex by $z$). Let $\mathcal{G}_2$ be the set of coalescences of a lollipop with a vertex $u$ of degree $2$  of the cycle as distinguished vertex and $T$ with a pendant vertex as distinguished vertex. The pendant vertices different from $z$ are denoted by $x$ and $y$ such that  $d(x,w)\leq d(y,w)$. Such a graph is drawn in figure \ref{G2}.

The aim of this section is to show the following theorem whose proof is summed up in table \ref{Proof2}.

\begin{theorem}\label{L4pascospG2}
A $L(4,k)$ cannot be cospectral with a graph $G\in\mathcal{G}_2$.
\end{theorem}

As in the previous section we can assume the length of the cycle of $G$ is even.


\begin{table}
\small
\begin{tabular}{|| c|c|c|c|c   ||c|c||}
\hline
\multicolumn{5}{|c||}{Graph} & Tool & \begin{tabular}{c} Propo-\\sition \end{tabular}\\

\hline
\hline
\multicolumn{5}{||c||}{$p=6$} & $\sum \lambda_i^6$ & \ref{sum6part3} $v)$ \\

\hline

\multirow{23}{*}{ $p\geq 8$ } & \multirow{8}{*}{$d(u,v)=1$} & \multicolumn{3}{|c||}{ $d(v,z)>1$ or $d(x,w)>1$ or  $d(y,w)>1$}&  $\sum \lambda_i^6$ &  \ref{sum6part3}  $ii)$ \\

\cline{3-7}

&       & \multicolumn{3}{|c||}{$d(w,u)=1$}&  $\sum \lambda_i^6$ &  \ref{sum6part3}  $iii)$ \\
\cline{3-7}
&        &\multicolumn{3}{|c||}{\begin{tabular}{c} $d(v,z)=1$ and \\ $d(x,w)=1$ and \\ $d(y,w)=1$ and $d(w,u)>1$ \end{tabular}}&  \begin{tabular}{c} $Q_G(2)$\\ and \\ $Q_G(1)$\end{tabular} &  \ref{duv=1}   \\

\cline{2-7}

&\multirow{13}{*}{$d(u,v)>1$} & \multicolumn{3}{|c||}{$d(x,w)>1$ or ($d(y,w)>1$ and $d(z,v)>1$)} &  $\sum \lambda_i^6$ & \ref{sum6part3}  $i)$  \\

\cline{3-7}

&&\multirow{10}{*}{\begin{tabular}{c}$d(x,w)=1$\\and\\$d(y,w)=1$\\ \\ or \\ \\ $d(x,w)=1$\\and\\$d(z,v)=1$ \end{tabular}} & d(u,w)=1 & \begin{tabular}{c}$d(v,z)>1$ or\\$d(y,w)>1$\end{tabular} & $\sum \lambda_i^6$  &  \ref{sum6part3}  $iv)$ \\

\cline{5-7}

&&&  & \begin{tabular}{c}$d(v,z)=1$ and\\$d(y,w)=1$\end{tabular} & $Q_G(2)$  &  \ref{ref2}   \\

\cline{4-7}

&&& d(u,w)>1 & \begin{tabular}{c}$d(v,z)=1$ and\\$d(y,w)=1$\end{tabular} & $\sum \lambda_i^6$  &  \ref{ref1}   \\
\cline{5-7}
&&&  & \begin{tabular}{c}$d(v,z)=1$ and\\$d(y,w)>1$  \\ \end{tabular} &   &  \ref{ref3}   \\
\cline{5-7}
&&&  & \begin{tabular}{c}$d(v,z)>1$ and\\$d(y,w)=1$  \\ \end{tabular} &   &  \ref{ref4}   \\

\hline
\hline
\end{tabular}
\normalsize
\caption{Proof of theorem  \ref{L4pascospG2} using a case disjunction over the possibilities for the values of $d$. An empty cell in the column \emph{tool} means that the proof uses more than three tools.}
\label{Proof2}
\end{table}


\begin{proposition}\label{sum6part3}
Let $G\in\mathcal{G}_2$. If one of the following properties is true
\\i) $d(x,w)>1$ or ($d(y,w)>1$ and $d(z,v)>1$),
\\ii) $d(u,v)=1$ and ($d(z,v)>1$ or $d(y,w)>1$),
\\iii) $d(u,v)=1$ and $d(u,w)=1$,
\\iv) $d(u,w)=1$ and ($d(z,v)>1$ or $d(y,w)>1$),
\\v)p=6,
\\then
$$\sum_i\lambda_i^6>20n+96$$
and $G$ cannot be cospectral with a lollipop $L(4,k)$.
\end{proposition}

\begin{demo}
For all cases we have $|P_2(G)|=n$, $|P_3(G)|=n+3$, $|S_{1,1,1}(G)|=3$.
Moreover, for the cases i) to iv) $|P_4(G)|>n+4$ and for the case v) $|P_4(G)|>n+2$ and $|C_6(G)|=1$ and we apply proposition \ref{sum6}.
\end{demo}

\begin{proposition}\label{ref1}
Let $G\in\mathcal{G}_2$ such that $p\geq8$, $d(u,v)>1,\ d(u,w)>1,\ d(z,v)=1,\ d(w,x)=d(w,y)=1,$ then $\sum_i\lambda_i^6<20n+96$ and $G$ cannot be cospectral with a lollipop $L(4,k)$.
\end{proposition}

\begin{demo}
The subgraphs $M$ of $G$ with $w_6(G)>0$ are $P_2,P_3,P_4,S_{1,1,1}$ and  $|P_2(G)|=n$, $|P_3(G)|=n+3$, $|S_{1,1,1}(G)|=3$, $|P_4(G)|=n+3$ and we apply proposition \ref{sum6}.
\end{demo}

\begin{proposition}\label{duv=1}
Let $G\in\mathcal{G}_2$ such that $d(u,v)=1,\ d(z,v)=1,\ d(w,y)=1,\ d(w,u)>1$, then $G$ is not cospectral with $L(4,k)$.
\end{proposition}

\begin{demo}
Since $d(w,x)\leq d(w,y)=1$ we have $d(w,x)=1$. Let $\alpha=d(u,w)$ (so $n=p+\alpha+3$), by theorem \ref{isthme} we get:
\begin{eqnarray*}
Q_G(X)&=&Q_{L(p,1)}(X)Q_{S_{1,1,\alpha-1}}(X)-Q_{P_{p}}(X)Q_{S_{1,1,\alpha-2}}(X)\\
&=&\big(XQ_{C_p}(X)-Q_{P_{p-1}}(X)\big)Q_{S_{1,1,\alpha-1}}(X)-Q_{P_{p}}(X)Q_{S_{1,1,\alpha-2}}(X)
\end{eqnarray*}

and (with property \ref{P(2)}) $Q_G(2)=-8p-4$. So $Q_G(2)=Q_{L(4,k)}(2)$ if and only if $-8p-4=-4n+16$ that is $\alpha=p+2$.


As a consequence
\begin{eqnarray*}
Q_{G(X)}&=&\big(XQ_{C_p}(X)-Q_{P_{p-1}}(X)\big)Q_{S_{1,1,p+1}}(X)-Q_{P_{p}}(X)Q_{S_{1,1,p}}(X)\\
&=&\big(XQ_{C_p}(X)-Q_{P_{p-1}}(X)\big)X\big(Q_{P_{p+3}}(X)-Q_{P_{p+1}} (X)\big)\\
&&-Q_{P_{p}}(X)X\big(Q_{P_{p+2}}(X)-Q_{P_{p}} (X)\big)
\end{eqnarray*}

By property \ref{Pp(1)} we have (let's note that $n=2p+5$):
\begin{itemize}
\item If $\overline{p}=\overline{0}$ (so $\overline{n}=\overline{5}$) then $Q_G(1)=1$.
\item If $\overline{p}=\overline{2}$ (so $\overline{n}=\overline{3}$) then $Q_G(1)=-4$.
\item If $\overline{p}=\overline{4}$ (so $\overline{n}=\overline{1}$) then $Q_G(1)=0$.
\end{itemize}
where $\overline{p}$ and $\overline{n}$ are $p$ and $n$ modulo $6$.
And by proposition \ref{Pl4(1)}, $G$ is not cospectral with $L(4,k)$.
\end{demo}

\begin{proposition}\label{ref2}
Let $G\in\mathcal{G}_2$ such that $d(u,v)>1,\  d(w,x)=d(w,y)=d(v,z)=d(u,w)=1$, then $G$ cannot be cospectral with $L(4,k)$.
\end{proposition}

\begin{demo}
Set $a=d(u,v)$ and $b=p-a$. We have:
\begin{eqnarray*}
Q_G(X)&=&Q_{L(p,1)}(X)Q_{P_{3}}(X)-X^2Q_{S_{1,a-1,b-1}}(X)\\
&=&\big(XQ_{C_p}(X)-Q_{P_{p-1}}(X)\big)Q_{P_{3}}(X)-X^2Q_{S_{1,a-1,b-1}}(X)
\end{eqnarray*}
and $Q_G(2)=-4p-4(2a+2b-ab)$. As $n=p+4$ we have $Q_{G}(2)+4n-16=-4(2p-ab)$ so $Q_{G}(2)+4n-16=0$ if and only if $ab=2p$.
\begin{itemize}
\item If $a=2$ then $2b=4+2b$, impossible.
\item If $a=3$ then $b=6$ and $p=9$, $p$ odd is impossible.
\item If $a=4$ then $b=4$ and $p=8$ we check that this graph is not cospectral with $L(4,8)$.
\item If $a>4$ then as $p\leq 2b$ we have $2p-ab<0$.
\end{itemize}
As a result $G$ is not cospectral with $L(4,k)$.
\end{demo}

\begin{proposition}\label{ref3}
Let $G\in\mathcal{G}_2$ such that $p\geq 8$, $d(u,v)>1$, $d(w,u)>1$, $d(w,y)>1$, $d(w,x)=d(v,z)=1$. Then $G$ is not cospectral with a lollipop $L(4,k)$.
\end{proposition}

\begin{demo}
Let $a=d(u,v)$, $b=p-a$, $\alpha=d(u,w)$, $\beta=d(w,y)\geq 2$. We have  $a\leq b$ and $p\leq 2b$ and $n=p+\alpha+\beta+2$.
\begin{eqnarray*}
Q_G(X)&=&Q_{L(p,1)}(X)Q_{S_{1,\alpha-1,\beta}}(X)-Q_{S_{1,a-1,b-1}}(X)Q_{S_{1,\alpha-2,\beta}}(X)\\
&=&\big(XQ_{C_p}(X)-Q_{P_{p-1}}(X)\big)Q_{S_{1,\alpha-1,\beta}}(X)-Q_{S_{1,a-1,b-1}}(X)Q_{S_{1,\alpha-2,\beta}}(X)
\end{eqnarray*}
Using property \ref{P(2)} we obtain 
$$Q_G(2)=-p(\alpha+2\beta-\alpha\beta+2)-(2a+2b-ab)(\alpha+3\beta-\alpha\beta+1)$$
The following inequality will be useful: $ab=(a-1)(b-1)+p-1\geq b-1+p-1 \geq \frac{3}{2}p-2$.
The main argument of this proof is that $Q_{G}(2)\neq -4n+16$ so $G$ cannot be cospectral with a lollipop $L(4,k)$ (proposition \ref{Pl4(2)}).
\begin{itemize}
\item Case $\beta=2$.
$Q_{G}(2)-(-4n+16)=(3\alpha-16)p+(7-\alpha)ab+4\alpha$.
\begin{itemize}
\item If $\alpha=2$ then $Q_{G}(2)-(-4n+16)=-10p+5ab+8\neq 0$ (otherwise $5$ divides $8$)
\item If $\alpha=3$ then  $Q_{G}(2)-(-4n+16)=-7p+4ab+12$. If  $a\geq 4$ then $-7p+4ab+12>0$ (because $p \leq 2b$). If $a=3$ then $-7p+4ab+12=5b-9\neq0$ (because $b\in\n$). If $a=2$ then $-7p+4ab+12=b-2\neq 0$ (because $a+b=p\geq8$).
\item If $4\leq \alpha\leq 7$ then 
\begin{eqnarray*}
Q_{G}(2)-(-4n+16)&\geq& (3\alpha-16)p+(7-\alpha)(\frac{3}{2}p-2)+4\alpha\\
&\geq& (\frac{3}{2}\alpha-\frac{11}{2})p-14+6\alpha\\
&\geq& \frac{p}{2}+10>0
\end{eqnarray*}
\item If $\alpha>7$ then the disjoint union $C_p\cup S_{1,2,5}$ is an induced subgraph of $G$ with twice the eigenvalue $2$ and by the interlacing theorem and theorem \ref{vp2}, $G$ is not cospectral with a lollipop.
\end{itemize}

\item Case $\beta \geq 3$ :
\begin{itemize}
\item $\alpha=2$. We have $|P_2(G)|=n$, $|P_3(G)|=n+3$, $|P_4(G)|=n+4$, $|S_{1,1,2}(G)|=7$, $|P_5(G)|=n+6 \textrm{ if } a>2$ and $|P_5(G)|=n+7 \textrm{ if } a=2$. By proposition \ref{sum6} we have $\sum \lambda_i^8=70n+588+16c_8$ if $a>2$ and in that case  $G$ in not cospectral with $L(4,k)$ (proposition \ref{cor0}). If $a=2$ then $Q_G(2)=-4p-4(\beta+3)=-4n+4\neq -4n+16$.

\item $\alpha=3$. 
$Q_G(2)+4n-16=-p(-\beta+5)-(2p-ab)\times4+4(p+\beta+5)-16=p(\beta-9)+4ab+4\beta+4$. But $\beta\geq3$ and $ab\geq\frac{3}{2}p-2$, so $Q_G(2)+4n-16\geq4\beta-4>0$.

\item $\alpha=4$.
\begin{itemize}
\item If $\beta\geq 5$ the disjoint union $C_p\cup S_{1,2,5}$ is an induced subgraph of $G$ with twice the eigenvalue $2$ and by the interlacing theorem and theorem \ref{vp2}, $G$ is not cospectral with a lollipop.
\item If $\beta=4$ then $Q_G(2)=ab>0$ and $Q_{L(4,k)}(2)< 0$ 
\item If $\beta=3$ then $Q_G(2)=2(ab-2p)$, $n=p+9$ and $Q_G(2)-(-4n+16)=2ab+20>0$

%
%

\end{itemize}

\item $\alpha>4$. The disjoint union $C_p\cup S_{1,3,3}$ is an induced subgraph of $G$ with twice the eigenvalue $2$ and by the interlacing theorem and theorem \ref{vp2}, $G$ is not cospectral with a lollipop.
\end{itemize}

\end{itemize}

\end{demo}

\begin{propriete}
Let $r\in\mathbb{R}$, $r>2$, we have $Q_{P_n}(r)=\alpha_1\beta_1^n+\alpha_2\beta_2^n$ with $\beta_1=\frac{r+\sqrt{r^2-4}}{2}>1$, $\beta_2=\frac{r-\sqrt{r^2-4}}{2}<1$, $\alpha_1=\frac{r-\beta_2}{\beta_1-\beta_2}>1$, $\alpha_2=1-\alpha_1<0$.
\end{propriete}
\begin{demo}
Let $(u_n)_{n\in\n}$ be the sequence $u_n=Q_{P_n}(r)$. We have $u_n=ru_{n-1}-u_{n-2}$, so $u_n=\alpha_1\beta_1^n+\alpha_2\beta_2^n$ where $\beta_1,\beta_2$ are roots of $X^2-rX+1$ and we note that $1=u_0=\alpha_1+\alpha_2$, $r=u_1=\alpha_1\beta_1+\alpha_2\beta_2$.
\end{demo}

\begin{lemme}\label{infty}
Let $G\in\mathcal{G}_2$ with $d(u,v)=2,\  d(w,x)=d(w,y)=1,\ d(v,z)>1$, $d(u,w)>d(v,z)$, then $G$ is not cospectral with a $L(4,k)$.
\end{lemme}

\begin{demo}
Let $\alpha=d(u,w)$, $l=d(v,z)$, we have $n=p+\alpha+l+2$. Applying theorem \ref{del_vertex} to the vertex at distance $1$ of $u$ and $v$, we have:
\begin{eqnarray*}
Q_G(X)=XQ_{S_{1,1,n-4}}(X)-Q_{P_l}(X)Q_{S_{1,1,\alpha+p-3}}(X)-Q_{P_{l+p-2}}(X)Q_{S_{1,1,\alpha-1}}(X)-2Q_{P_l}Q_{S_{1,1,\alpha-1}}
\end{eqnarray*}
and applying theorem \ref{del_vertex} to the vertex of degree $2$ of the cycle of $L(4,k)$ we have:
\begin{eqnarray*}
Q_{L(4,k)}(X)=XQ_{S_{1,1,n-4}}(X)-2Q_{P_{n-2}}(X)-2Q_{P_{n-4}}(X)
\end{eqnarray*}
Noting that $Q_{S_{1,1,c}}=X(Q_{P_{c+2}}(X)-Q_{P_c}(X))$ and $Q_{P_{n-2}}(X)+Q_{P_{n-4}}(X)=XQ_{P_{n-3}}(X)$ we have:

\begin{eqnarray*}
Q_G(X)-Q_{L(4,k)}(X)&=&-XQ_{P_l}(X)Q_{P_{\alpha+p-1}}(X)+XQ_{P_l}(X)Q_{P_{\alpha+p-3}}(X)\\
&&-XQ_{P_{l+p-2}}(X)Q_{P_{\alpha+1}}(X)+XQ_{P_{l+p-2}}(X)Q_{P_{\alpha-1}}(X)\\
&&-2Q_{P_l}(X)Q_{P_{\alpha+1}}(X)+2Q_{P_l}(X)Q_{P_{\alpha-1}}(X)+2XQ_{P_{n-3}}(X)
\end{eqnarray*}

According to the previous property, we have for $r>2$:

\begin{eqnarray*}
Q_G(r)-Q_{L(4,k)}(r)&=&-r\alpha_1^2\beta_1^{n-3}-r\alpha_2^2\beta_2^{n-3}-r\alpha_1\alpha_2\beta_1^{\alpha+p-1-l}-r\alpha_1\alpha_2\beta_2^{\alpha+p-1-l}\\
&&+r\alpha_1^2\beta_1^{n-5}+r\alpha_2^2\beta_2^{n-5}+r\alpha_1\alpha_2\beta_1^{\alpha+p-3-l}+r\alpha_1\alpha_2\beta_2^{\alpha+p-3-l}\\
&&-r\alpha_1^2\beta_1^{n-3}-r\alpha_2^2\beta_2^{n-3}-r\alpha_1\alpha_2\beta_1^{l+p-2}\beta_2^{\alpha+1}-r\alpha_1\alpha_2\beta_2^{l+p-2}\beta_1^{\alpha+1}\\
&&+r\alpha_1^2\beta^{n-5}+r\alpha_2^2\beta^{n-5}+r\alpha_1\alpha_2\beta_1^{l+p-2}\beta_2^{\alpha-1}+r\alpha_1\alpha_2\beta_2^{l+p-2}\beta_1^{\alpha-1}\\
&&+2r\alpha_1\beta_1^{n-3}+2r\alpha_2\beta_2^{n-3}
\end{eqnarray*}
Let $x=\alpha+p-l-1$ and $y=|l+p-\alpha-1|$, we have $x>y$.
\begin{eqnarray*}
Q_G(r)-Q_{L(4,k)}(r)&=&2r\left( (\alpha_1-\alpha_1^2)\beta_1^2+\alpha_1^2\right)\beta_1^{n-5}+2r\left( (\alpha_2-\alpha_2^2)\beta_2^2+\alpha_2^2\right)\beta_2^{n-5}\\
&&-r\alpha_1\alpha_2(\beta_1^x-\beta_1^{x-2}-\beta_1^y+\beta_1^{y-2})
-r\alpha_1\alpha_2(\beta_2^x-\beta_2^{x-2}-\beta_2^y+\beta_2^{y-2})
\end{eqnarray*}
but we have the four following equalities:
\begin{eqnarray*}
&\alpha_1\alpha_2=\alpha_1-\alpha_1^2=\frac{-1}{r^2-4}\\
&(\alpha_1-\alpha_1^2)\beta_1^2+\alpha_1^2=0\\
&(\alpha_2-\alpha_2^2)\beta_2^2+\alpha_2^2=0\\
&\beta_2=\beta_1^{-1}
\end{eqnarray*}
so
\begin{eqnarray*}
&Q_G(r)-Q_{L(4,k)}(r)=\frac{r}{r^2-4}\left( \beta_1^x-\beta_1^{x-2}-\beta_1^y+\beta_1^{y-2}+\beta_1^{-x}-\beta_1^{-x+2}-\beta_1^{-y}+\beta_1^{-y+2}\right)
\end{eqnarray*}
As
\begin{eqnarray*}
&\lim_{r\rightarrow + \infty}\frac{\beta_1^x-\beta_1^{x-2}-\beta_1^y+\beta_1^{y-2}+\beta_1^{-x}-\beta_1^{-x+2}-\beta_1^{-y}+\beta_1^{-y+2}}{r}=+\infty  \textrm{ (note that $x>2$)}
\end{eqnarray*}
we have
$$\lim_{r\rightarrow + \infty} Q_G(r)-Q_{L(4,k)}(r) =+\infty$$
and $G$ is not cospectral with $L(4,k)$.

\end{demo}

\begin{proposition}\label{ref4}
Let $G\in\mathcal{G}_2$ with $p\geq 8$, $d(u,v)>1,\  d(w,x)=d(w,y)=1,\ d(v,z)>1$, $d(u,w)>1$, then $G$ is not cospectral with a $L(4,k)$.
\end{proposition}

\begin{demo}

We distinguish the following cases :
\begin{itemize}
\item case 1 : $d(u,v)>2$ and $d(u,w)>2$ and $d(z,v)>2$
\item case 2 : $d(u,v)>2$ and $d(u,w)>2$ and $d(z,v)=2$
\item case 3 : $d(u,v)>2$ and $d(u,w)=2$ and $d(z,v)>2$
\item case 4 : $d(u,v)>2$ and $d(u,w)=2$ and $d(z,v)=2$
\item case 5 : $d(u,v)=2$ 
\end{itemize}
\textbullet\ For  cases  1 and 4 we have $|P_2(G)|=n,\ |P_3(G)|=n+3,\ |P_4(G)|=n+4$,  $|S_{1,1,1}(G)|=3,\ |S_{1,1,2}(G)|=7,\  |P_5(G)|=n+6$ , $|L(4,1)(G)|=0,\ |L(4,2)(G)|=0$ so (proposition \ref{sum6}) $\sum\lambda_i^8=70n+588+16c_8$ and $G$ is not cospectral with $L(4,k)$ (proposition \ref{cor0}).
\\\textbullet\ For  cases 2, 3 and 5, let us compute $Q_G(2)$. Let $a=d(u,v)$, $b=p-a$, $\alpha=d(u,w)$, $l=d(v,z)$.

\begin{eqnarray*}
Q_G(X)&=&Q_{L(p,l)}(X)Q_{S_{1,1,\alpha-1}}(X)-Q_{S_{a-1,b-1,l}}(X)Q_{S_{1,1,\alpha-2}}(X)\\
&=&\big(P_{{C_p}}(X)Q_{P_l}(X)-Q_{P_{p-1}}(X)Q_{P_{l-1}}(X)\big)Q_{S_{1,1,\alpha-1}}(X)\\
&&-Q_{S_{a-1,b-1,l}}(X)Q_{S_{1,1,\alpha-2}}(X)
\end{eqnarray*}

Using property \ref{P(2)} we have $Q_G(2)+4n-16=-8lp+4abl+4\alpha+4l-8$ and $G$ is cospectral with $L(4,k)$ only if  $Q_G(2)+4n-16=0$ that is  $\alpha=l(2p-ab-1)+2$. 
\begin{itemize}
\item For  case 3 we have $\alpha=2$ so $2p-ab+1=0$ and $a$ is odd. If $a=3$ then $b=5$ and $p=8$. We have $|P_2(G)|=n$, $|P_3(G)|=n+3$, $|P_4(G)|=n+4$, $|P_5(G)|=n+7$, $|S_{1,1,1}(G)|=3$,  $|S_{1,1,2}(G)|=7$, $|C_8(G)|=1$.
So $\sum \lambda_i^{8}=70n+612$ and in this case $G$ is not cospectral with $L(4,k)$ (proposition \ref{cor0}).
If $a\geq 5$ then  $2p-ab-1\leq 4b-5b-1<0$ and this finishes the case 3.
\item For  case 2, $|P_2(G)|=n,\ |P_3(G)|=n+3,\ |P_4(G)|=n+4$,  $|S_{1,1,1}(G)|=3,\ |S_{1,1,2}(G)|=7,\  |P_5(G)|=n+5$ , $|L(4,1)(G)|=0,\ |L(4,2)(G)|=0$  so (proposition \ref{sum6}) $\sum\lambda_i^8=70n+580+16c_8$ and $G$ is cospectral with $L(4,k)$ only if $p=8$. We have $l=2$ and $\alpha=l(2p-ab-1)+2$ so the graphs that can be cospectral with $L(4,k)$ are the ones with  $a=3, b=5$ so $\alpha=2$, impossible, or  $a=4,b=4$ so $\alpha=0$, impossible.
\item For  case 5, $G$ is cospectral with $L(4,k)$ only if  $\alpha=3l+2$, but this is impossible according to lemma \ref{infty}.
\end{itemize}

\end{demo}

%
%

\subsubsection{Unicyclic graphs with exactly three vertices of maximum degree $3$, all of them belonging to the cycle.}\label{sec_G3}\label{deg3}

Let $\mathcal{G}_3$ be the set of the graphs $G$ obtained in the following way:
\begin{itemize}
\item Do the coalescence of a lollipop $L(p,k)$, $p\geq 6$, $k\geq1$ with a vertex of degree $2$ of the cycle  as distinguished vertex and  a path with a pendant vertex as distinguished vertex.
\item Do the coalescence of the previous graph with a vertex of the cycle of degree $2$ as distinguished vertex and  a path with a pendant vertex as distinguished vertex.
\end{itemize}

\begin{figure}
\begin{center}
\includegraphics[scale=0.5]{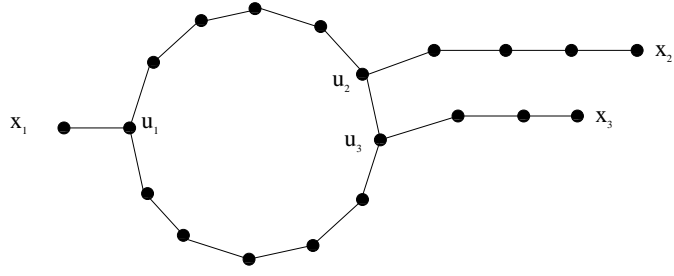}
\caption{A graph $G\in\mathcal{G}_3$}
\label{G3}
\end{center}
\end{figure}

We denote by $u_1,u_2,u_3$ the three vertices of degree $3$ and by $x_1,x_2,x_3$ the pendant vertices such that $d(x_i,u_i)=\min_j {d(x_i,u_j)}$. Un example is given in figure \ref{G3}


The aim of this section is to show the following theorem whose proof is summed up in table \ref{proof3} :

\begin{theorem}\label{L4pascosp}
A lollipop  $L(4,k)$ cannot be cospectral with a graph $G\in\mathcal{G}_3$.
\end{theorem}

As in the previous sections we assume that the cycle of $G$ is even.


\begin{landscape}
\begin{table}
\begin{tabular}{|| c|c|c|c|c|c    ||c|c||}
\hline
\multicolumn{6}{|c||}{Graph} & Tool & \begin{tabular}{c} Propo-\\sition \end{tabular}\\

\hline
\hline
\multicolumn{6}{||c||}{$p=6$} & $\sum \lambda_i^6$ & \ref{cor} $iv)$ \\

\hline

\multicolumn{6}{||c||}{$\exists i,j,\ i\neq j :\ d(x_i,u_i)>1$ and $d(x_j,u_j)>1 $} & $\sum \lambda_i^6$ & \ref{cor} $i)$ \\

\hline

\hline

\multicolumn{6}{||c||}{$\exists r,s,t,\ r\neq s, s\neq t, r\neq t, :\ d(u_r,u_s)=1$ and $d(u_s,u_t)=1 $} & $\sum \lambda_i^6$ & \ref{cor} $ii)$\\

\hline

\multirow{17}{*}{\begin{tabular}{c} $p\geq 8$\\ and \\$\exists i,j,k$ \\two by two\\ distinct\\$\exists r,s,t$ \\two by two\\ distinct:\\ $d(x_i,u_i)=1$\\$d(x_j,u_j)=1$\\$d(x_k,u_k)\geq1$\\$d(u_r,u_s)>1$\\$d(u_s,u_t)>1$\end{tabular}  } & \multirow{3}{*}{$\begin{array}{c}d(x_k,u_k)=1\end{array}$} & \multicolumn{4}{|c||}{$d(u_r,u_t)>1$} &  $\sum \lambda_i^6$ &  \ref{corbis}   \\

\cline{3-8}

&&\multirow{2}{*}{$d(u_r,u_t)=1$} & \multicolumn{3}{|c||}{$d(u_r,u_s)=2$ or $d(u_s,u_t)=2$} &  $\sum \lambda_i^8$ & \ref{cor2}   \\

\cline{4-8}

&&& \multicolumn{3}{|c||}{$d(u_r,u_s)>2$ and $d(u_s,u_t)>2$ }   & $\sum \lambda_i^8$ &  \ref{cor2}   \\
\cline{2-8}

 & \multirow{15}{*}{$\begin{array}{c} d(x_k,u_k)>1\end{array}$} 
&   \multicolumn{4}{|c||}{ $d(u_r,u_t)=1$ } & $\sum \lambda_i^6$ &  \ref{cor} $iii)$ \\
\cline{3-8}
&&  \multirow{14}{*}{$d(u_r,u_t)>1$} & \multicolumn{3}{|c||}{ $d(x_k,u_k)=2$ } & $\sum\lambda_i^8$ &  \ref{cor3} \\
\cline{4-8}
&&                              &    \multirow{13}{*}{ $d(x_k,u_k)>2$}      & \multicolumn{2}{|c||}{$\forall l_1,l_2,\ d(u_{l_1},u_{l_2})>2$} & $\sum\lambda_i^8$ &  \ref{cor4} $iii)$ \\
\cline{5-8}
&&                               &                             & \multicolumn{2}{|c||}{\begin{tabular}{c}$\exists r,s,t$  \\$d(u_{r},u_{s})=2$ and \\$d(u_{r},u_{t})>2$ and \\$d(u_{s},u_{t})>2$\end{tabular}} & $\sum\lambda_i^8$ &  \ref{cor4} $ii)$ \\  
\cline{5-8}

&&                               & & \multirow{3}{*}{\begin{tabular}{c} $d(u_{i},u_{j})=2$ \\and\\ $d(u_{j},u_{k})=2$ \\ \end{tabular}} 
&$p=8$ & $\sum\lambda_i^8$ &  \ref{cor4} $i)$ \\
\cline{6-8}
&&                               & & &  &  & \\
&&                               & &  & $p\geq 10$ & $\sum\lambda_i^{10}$ &  \ref{sum10_3} \\

\cline{5-8}

&&                               & & \multicolumn{2}{|c||}{\begin{tabular}{c} $d(u_{i},u_{k})=2$ \\and\\ $d(u_{j},u_{k})=2$ \\ \end{tabular}} &  \begin{tabular}{c} $Q_G(2)$ \\and\\ $R(0)$ \\ \end{tabular}  &  \ref{new_prop} \\

\hline
\hline
\end{tabular}
\caption{Proof of theorem \ref{L4pascosp} using a case disjunction over the possibilities for the values of $d$. $R$ denotes the polynomial $R(X)=\frac{Q_G(X)}{X^2}$}
\label{proof3}
\end{table}
\end{landscape}


\begin{proposition}\label{cor}
Let $G \in\mathcal{G}_3$. If one of the following properties is true:
\\i) $\exists i,j$, $i\neq j$ : $d(x_i,u_i)>1$, $d(x_j,u_j)>1$,
\\ii) $\exists r,s,t$, $r\neq s, r\neq t, s\neq t$ : $d(u_r,u_s)=d(u_s,u_t)=1$,
\\iii) $\exists i,r,t$ : $d(x_i,u_i)>1$, $d(u_r,u_t)=1$,
\\iv) $p=6$,
\\then
$$\sum_i\lambda_i^6>20n+96$$
and $G$ cannot be cospectral with a lollipop $L(4,k)$.
\end{proposition}

\begin{demo}
For  cases i) to iii) we have $|P_2(G)|=n,\ |P_3(G)|=n+3,\ |S_{1,1,1}(G)|=3,\ |P_4(G)|>n+4$ and we apply proposition \ref{sum6}.
\\For  case iv) we have  $|P_2(G)|=n,\ |P_3(G)|=n+3,\ |S_{1,1,1}(G)|=3,\ |P_4(G)|>n+2,\ |C_6(G)|=1$ and we apply proposition \ref{sum6}.
\end{demo}

\begin{proposition}\label{corbis}
Let $G \in\mathcal{G}_3$ such that $p>6$, $\forall i,r,s ,\ d(u_i,x_i)=1,\ d(u_r,u_s)>1$. Then
$$\sum_i\lambda_i^6=20n+90$$
and $G$ cannot be cospectral with a lollipop $L(4,k)$.
\end{proposition}

\begin{demo}
We have $|P_2(G)|=n$, $|P_3(G)|=n+3$, $|P_4(G)|=n+3$, $|S_{1,1,1}(G)|=3$ and no $p$-cycle for $p\leq6$. We conclude with proposition \ref{sum6}.
\end{demo}

The following three propositions compute $\sum\lambda_i^8$ for some $G\in\mathcal{G}_3$, their proofs are based on counting motifs in $\mathcal{M}_8(G)$  which is done in a summary table \ref{proof8}.

\begin{proposition}\label{cor2} 
Let $G\in\mathcal{G}_3$ such that $p\geq 8$, $\forall i,\ d(u_i,x_i)=1$, and $\exists r,s,t$ two by two distinct $:\ d(u_r,u_t)=1,\ d(u_r,u_s)>1,\ d(u_s,u_t)>1$. 
Then: $$\sum_i\lambda_i^8=\left\{\begin{array}{l} 
70n+588+16c_8\ \textrm{ if } d(u_r,u_s)=2 \textrm{ or } d(u_s,u_t)=2\\
70n+580+16c_8\ \textrm{ otherwise } \\

\end{array}\right.$$
and $G$ cannot be cospectral with a $L(4,k)$.
\end{proposition}

\begin{demo}
Using table \ref{proof8}, we apply proposition \ref{sum6} to compute $\sum_i\lambda_i^8$. The only case for which $\sum \lambda_i^8=70n+596$ is when $\forall i,\ d(u_i,x_i)=1$, $\exists r,s,t$ two by two distinct $:\ d(u_r,u_t)=1,\ d(u_r,u_s)>2,\ d(u_s,u_t)>2$ and $c_8=1$. This case is drawn in figure \ref{C8_3} and we check that it is not cospectral with $L(4,7)$ by comparing spectral radii (see tables \ref{spect_L4k} and \ref{spect_special} in appendix).

\begin{figure}[htbp]
\begin{center}
\includegraphics[scale=0.4]{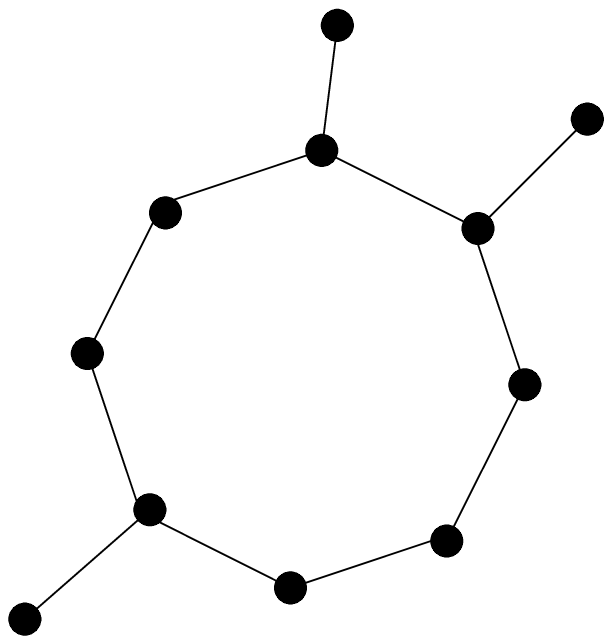}
\caption{}
\label{C8_3}
\end{center}
\end{figure}

\end{demo}

\begin{proposition}\label{cor3}
Let $G\in\mathcal{G}_3$ such that $p\geq 8$,  $\exists i,j,k:\ d(u_i,x_i)=d(u_j,x_j)=1,\ d(u_k,x_k)=2$. We distinguish the three following cases 
\begin{itemize}
\item case 1 : $\exists r,s,t$, $r\neq s, r\neq t, s\neq t$: $d(u_r,u_s)=d(u_s,u_t)=2$.
\item case 2 : $\exists r,s,t$, $r\neq s, r\neq t, s\neq t$: $d(u_r,u_s)=2$ and $d(u_r,u_t)>2$ and $d(u_s,u_t)>2$.
\item case 3:  $\forall s,t,\  d(u_{s},u_{t})>2$.
\end{itemize} 
Then: 
$$\sum_i\lambda_i^8=\left\{\begin{array}{l} 
70n+588+16c_8\ \textrm{ for the case 1 }  \\
70n+580+16c_8\ \textrm{ for the case 2 }  \\
70n+572\ \textrm{ for the case 3 }  \\
\end{array}\right.$$

and $G$ cannot be cospectral with a lollipop $L(4,k)$.

%

\end{proposition}

\begin{demo}
Using table \ref{proof8}, we apply proposition \ref{sum6} to compute $\sum_i\lambda_i^8$. Under the hypotheses of the proposition, the only cases for which $\sum \lambda_i^8=70n+596$ is when $c_8=1$ in case 2. These cases are drawn in figure \ref{C8_3b} and we check that they are not cospectral with $L(4,8)$ by comparing spectral radii (see tables \ref{spect_L4k} and \ref{spect_special} in appendix).

\begin{figure}[htbp]
\begin{center}
\includegraphics[scale=0.4]{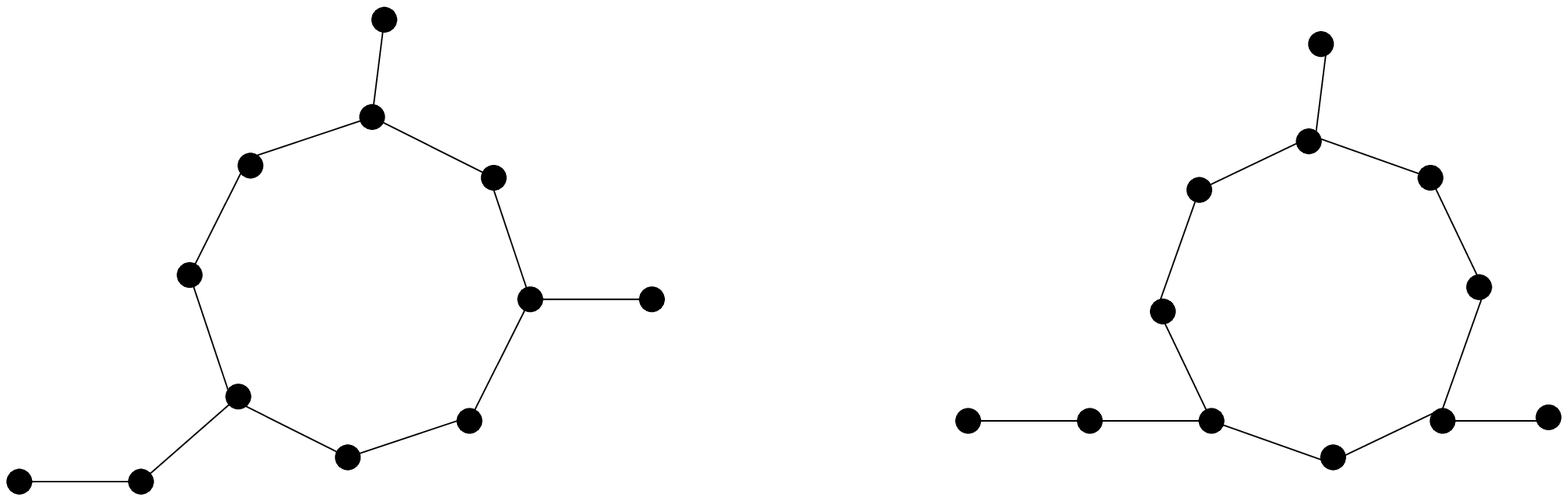}
\caption{}
\label{C8_3b}
\end{center}
\end{figure}

\end{demo}

\begin{proposition}\label{cor4}
Let $G\in\mathcal{G}_3$ such that $p\geq 8$, $\exists i,j,k:\ d(u_i,x_i)=d(u_j,x_j)=1,\ d(u_k,x_k)>2$.  We distinguish the three following cases 

\begin{itemize}
\item case 1 : $\exists r,s,t$, $r\neq s, r\neq t, s\neq t$: $d(u_r,u_s)=d(u_s,u_t)=2$.
\item case 2 : $\exists r,s,t$, $r\neq s, r\neq t, s\neq t$: $d(u_r,u_s)=2$ and $d(u_r,u_t)>2$ and $d(u_s,u_t)>2$.
\item case 3:  $\forall s,t,\  d(u_{s},u_{t})>2$.
\end{itemize}

Then: 
$$\sum_i\lambda_i^8=\left\{\begin{array}{l} 
70n+596+16c_8\ \textrm{ for the case 1 }  \\
70n+588+16c_8\ \textrm{ for the case 2 }  \\
70n+580\ \textrm{ for the case 3 }  \\
\end{array}\right.$$
and $G$ cannot be cospectral with a lollipop in the cases 2 and 3 and in the case 1 if $c_8=1$.

\end{proposition}


\begin{table}
$$
\begin{array}{||c|c||c|c|c||}
\hline
M & w_8(M) &  |M(G_a )| & |M(G_b )| & |M(G_c )| \\
\hline
\hline
P_2 & 2  &   n &n&n\\
\hline
P_3 & 28 &  n+3 &n+3&n+3\\
\hline
P_4 & 32 &  n+4 & n+4 &n+4\\
\hline
P_5 &8 & 

\begin{array}{c} n+4  \textrm{ case 1 }\\ n+3 \textrm{ case 2 }   \end{array} &
\begin{array}{c} n+6  \textrm{ case 1 }\\n+5  \textrm{ case 2 }  \\ n+4 \textrm{ case 3 } \end{array} &
\begin{array}{c} n+7  \textrm{ case 1 }\\n+6  \textrm{ case 2 }  \\ n+5 \textrm{ case 3 } \end{array}\\
\hline 
S_{1,1,1} & 72 & 3 &3&3\\
\hline
S_{1,1,2}   & 16 & 8&7&7\\
\hline
C_4 & 264 &  0 & 0& 0\\
\hline
L(4,1) & 112  & 0 & 0& 0 \\
\hline
L(4,2) & 16 &0 &0&0\\
\hline
C_8 & 16 & c_8 & 
\begin{array}{c} c_8  \textrm{ case 1 }\\c_8  \textrm{ case 2 }  \\0\textrm{ case 3 }  \end{array}&
\begin{array}{c} c_8  \textrm{ case 1 } \\c_8  \textrm{ case 2 }  \\0\textrm{ case 3 } \end{array}\\
\hline\hline
\multicolumn{2}{||c||}{\sum_i \lambda_i^8=}&
\begin{array}{c} 70n+588+16c_8 \\\textrm{ for the case 1 }   \\\\ 70n+580+16c_8 \\\textrm{ for the case 2 } \end{array} &
\begin{array}{c}  70n+588+16c_8 \\\textrm{ for the case 1 }  \\\\70n+580+16c_8 \\\textrm{ for the case 2 }  \\\\70n+572\\\textrm{ for the case 3 } \end{array}& 
\begin{array}{c}  70n+596+16c_8 \\\textrm{ for the case 1 }\\\\70n+588+16c_8 \\\textrm{ for the case 2 }  \\\\ 70n+580\\\textrm{ for the case 3 }  \end{array} 
\\\hline
\end{array}
$$

\caption{Count of the motifs of some graphs $G\in\mathcal{G}_3$. We denote by $G_a$ (resp. $G_b$, $G_c$) a graph described in proposition \ref{cor2} (resp. \ref{cor3}, \ref{cor4}). }
\label{proof8}
\end{table}

%

The two following propositions solve the case 1 of proposition \ref{cor4} when $c_8=0$.
\begin{proposition}\label{sum10_3}
Let $G\in\mathcal{G}_3$ such that $p\geq 10$, $\exists i,j,k:\ d(u_i,x_i)=d(u_j,x_j)=1$, $d(u_k,x_k)>2$, $d(u_i,u_j)=d(u_j,u_k)=2$.
\\Then:
$$\sum_i\lambda_i^{10}=\left\{\begin{array}{l} 
252n+3340+20c_{10}\ \textrm{ if }  d(u_k,x_k)=3\\
252n+3350+20c_{10}\ \textrm{ if }  d(u_k,x_k)>3\\
\end{array}\right.$$
where $c_{10}=|C_{10}(G)|$. And $G$ cannot be cospectral with $L(4,k)$.
\end{proposition}

\begin{demo}
We have $|P_2(G)|=n$, $|P_3(G)|=n+3$, $|P_4(G)|=n+4$, $|P_5(G)|=n+7$, $|P_6(G)|=n+6$ if $d(u_k,x_k)=3$, $|P_6(G)|=n+7$ if $d(u_k,x_k)>$3, $|S_{1,1,1}(G)|=3$,  $|S_{1,1,2}(G)|=7$,  $|S_{1,2,2}(G)|=5$,  $|S_{1,1,3}(G)|=11$,  and no others subgraphs $M$ such that $w_k(M)>0$. We then apply proposition \ref{sum6}. The only case for which $\sum\lambda_i^{10}=252n+3360$ is for the graph of figure \ref{C10_113}, and we check that it is not cospectral with $L(4,11)$ by comparing spectral radii (see tables \ref{spect_L4k} and \ref{spect_special} in appendix).

\begin{figure}[htbp]
\begin{center}
\includegraphics[scale=0.3]{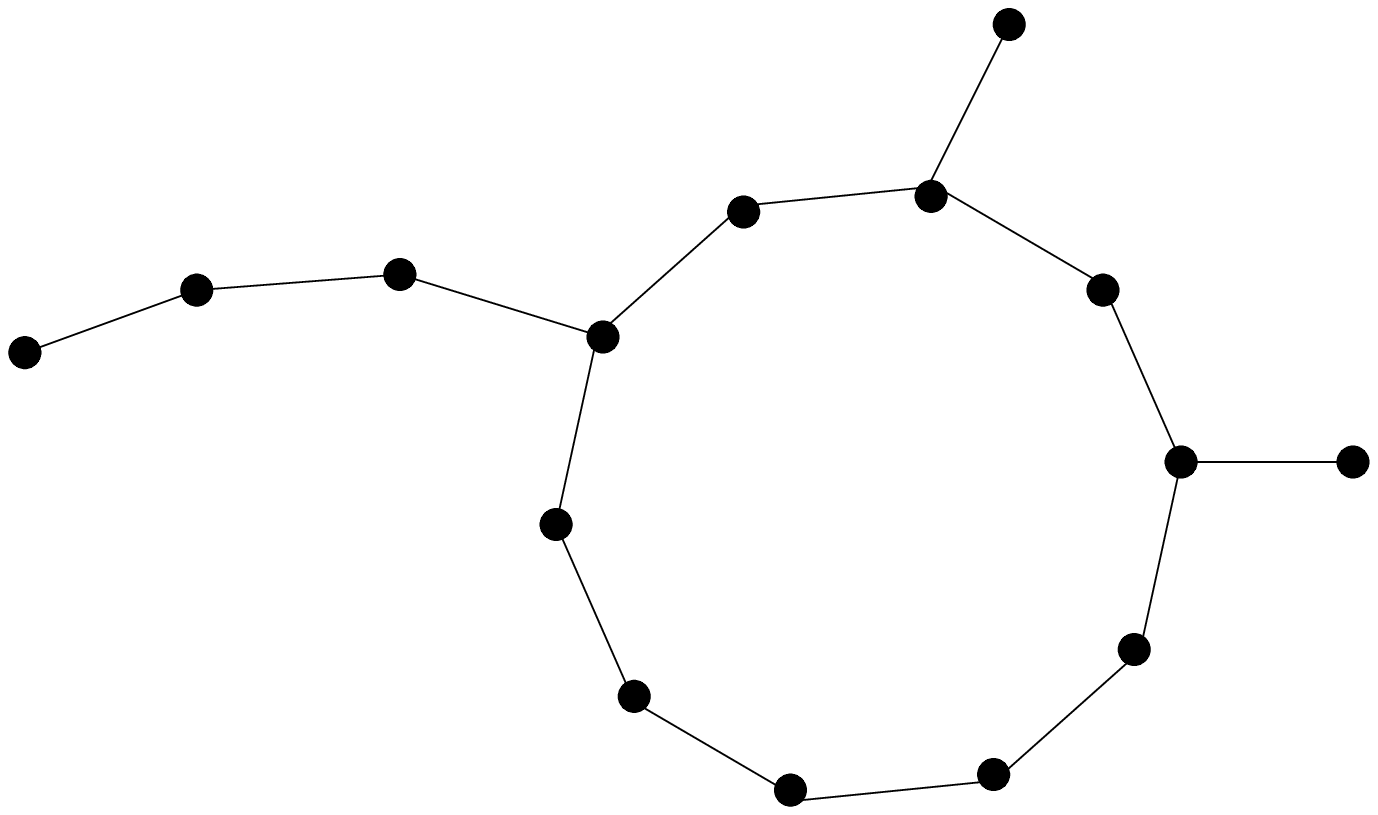}
\caption{}
\label{C10_113}
\end{center}
\end{figure}

\end{demo}

\begin{proposition}\label{new_prop}
Let $G\in\mathcal{G}_3$ such that $p\geq 10$, $\exists i,j,k:\ d(u_i,x_i)=d(u_j,x_j)=1$, $d(u_k,x_k)>2$, $d(u_i,u_k)=d(u_j,u_k)=2$.
Then $G$ cannot be cospectral with $L(4,k)$.
\end{proposition}
 
\begin{demo}
Let $G$ be a graph cospectral with $L(4,k)$ and let $q=d(u_k,x_k)$ (we have $n=p+q+2$). Applying theorem \ref{del_vertex} to the vertex $u_k$, we have:
$$Q_G(X)=XQ_{T_{p+1}}(X)Q_{P_q}(X)-2Q_{S_{1,1,p-3}}(X)Q_{P_q}(X)-Q_{T_{p+1}}(X)Q_{P_{q-1}}(X)-2X^2Q_{P_q}(X)$$
Property \ref{P(2)} gives $Q_G(2)=-16(q+1)$ and according to proposition \ref{Pl4(2)} $G$ is cospectral with a lollipop $L(4,k)$ only if $-16(q+1)=-4n+16$ \textit{ie}  $p=3q+6$ and $q$ is necessarily even.
\\

Using $Q_{S_{1,1,c}}(X)=X(Q_{P_{c+2}}(X)-Q_{P_{c}}(X))$ we have that if $c$ is odd then $0$ is an eigenvalue of $S_{1,1,c}$ with multiplicity $2$ and if $R(X)=\frac{Q_{S_{1,1,c}}(X)}{X}$ then $R(0)=(-1)^{\frac{c+1}{2}}(c+2)$. 
The relation $Q_{T_n}(X)=XQ_{S_{1,1,n-4}}(X)-XQ_{S_{1,1,n-2}}(X)$ implies that $0$ is an eigenvalue of $T_n$ with multiplicity $2$.
\\Let $R(X)=\frac{Q_G(X)}{X^2}$. Property \ref{P(0)} gives $$R(0)=\left\{\begin{array}{c} -2p \textrm{ if $q\equiv0[4]$}\\-2p+4 \textrm{ if $q\equiv0[4]$}\end{array}\right.$$
If $q\equiv0[4]$ then   according to proposition \ref{Pl4(0)}, $G$ is cospectral with a lollipop $L(4,k)$ only if $-2p=-n$ \textit{ie} $p=q+2$ which contradicts $p=3q+6$.
\\If $q\equiv2[4]$ then  according to proposition \ref{Pl4(0)}, $G$ is cospectral with a lollipop $L(4,k)$ only if $-2p+4=-n$ \textit{ie} $p=q+6$ which contradicts $p=3q+6$.
\end{demo}

\subsubsection{Unicyclic graphs without vertices of degree $3$ and only one vertex of maximum degree $4$}\label{deg4}\label{sec_G4}

\begin{figure}[htbp]
\begin{center}
\includegraphics[scale=0.5]{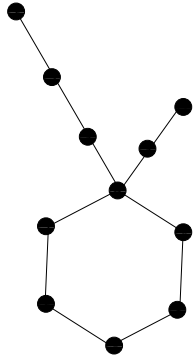}
\caption{ $\gamma_{6,2,3}$ }
\label{Mpk1k2un}
\end{center}
\end{figure}

The graph $\gamma_{p,k_1,k_2}$ is the coalescence of a lollipop $L(p,k_1)$ with the vertex of degree $3$ as distinguished vertex and a path $P_{k_2+1}$ with a pendant vertex as distinguished vertex (cf figure \ref{Mpk1k2un} for an example).


\begin{proposition}\label{cor00}
For a graph $\gamma_{p,k_1,k_2}$ with $p>4$  we have: 
$$\sum_i\lambda_i^6=\left\{\begin{array}{l} 
20n+96+12c_6\ \textrm{ if } k_1=k_2=1\\
20n+108+12c_6\ \textrm{ if } k_1>1,k_2=1\\
20n+120+12c_6\ \textrm{ if } k_1>1, k_2>1\\
\end{array}\right.$$
where $c_6=|C_6(G)|$.

\end{proposition}

\begin{demo}
We have $|P_2(G)|=n$, $|P_3(G)|=n+3$, $|S_{1,1,1}(G)|=4$ and
\begin{itemize}
\item  $|P_4(G)|=n+2$ if $k_1=k_2=1$ 
\item  $|P_4(G)|=n+4$ if $k_1>k_2=1$
\item  $|P_4(G)|=n+6$ if $k_1\geq k_2>1$
\end{itemize}
and we apply proposition \ref{sum6}.
\end{demo}

\begin{proposition}\label{lemL4kMp}
A lollipop $L(4,k)$ cannot be cospectral with a graph $\gamma_{p,1,1}$.
\end{proposition}

\begin{demo}
 The graphs $L(4,k)$ and $\gamma_{p,1,1}$ have $n=k+4=p+2$ vertices. Let us show that $Q_{L(4,k)}(2)\neq P_{\gamma_{p,1,1}}(2)$.
Using twice the theorem \ref{del_pendant}:
 \begin{eqnarray*}
 P_{\gamma_{p,1,1}}(X)&=&XQ_{L(p,1)}(X)-XQ_{P_{p-1}}(X)\\
 &=&X(XQ_{C_p}(X)-Q_{P_{p-1}}(X))-XQ_{P_{p-1}}(X)\\
\end{eqnarray*}
%
%

And by proposition \ref{P(2)}, $P_{\gamma_{p,1,1}}(2)=-4p=-4n+8$ which contradicts  $Q_{L(4,k)}(2)=-4n+16$ (proposition \ref{Pl4(2)}).
%
\end{demo}

\begin{theorem}\label{gamma}
A lollipop $L(4,k)$ cannot be cospectral with $\gamma_{p,k_1,k_2}$, $p>4$.
\end{theorem}

\begin{demo}
It is a straightforward consequence of propositions \ref{cor00} and \ref{lemL4kMp}.
\end{demo}

\subsubsection{Key theorem}

\begin{theorem}
Let $G$ be a graph cospectral with a lollipop $L(4,k)$ then $G$ possesses a $4$-cycle.
\end{theorem}

\begin{demo}

Let $G$ be a graph cospectral with $L(4,k)$ then $G$ is connected, unicyclic and bipartite (so the length of the cycle is even).
Let $n_j$ be the number of vertices of degree $j$, $j\in\{1,2,3,4\}$, of $G$ (remind by theorem \ref{degre_max} that the maximum degree of $G$ is less than or equal to $4$).
We have (by proposition \ref{sum4}), 

$$\sum_i\lambda_i^4=8c+2n+4(n_2+3n_3+6n_4)$$ were $c=1$ if $G$ has a $4$-cycle and $c=0$ otherwise. Moreover for $L(4,k)$ we have 
$$\sum_i\lambda_i^4=8+2n+4(n+1)$$ so 
$$4n+12=4(n_2+3n_3+6n_4)+8c$$
We know that $n=n_1+n_2+n_3+n_4$ and $2m=2n=n_1+2n_2+3n_3+4n_4$ (the sum of the degrees is twice the number of edges) so $n=n_2+2n_3+3n_4$ and $n_1=n_3+2n_4$. We get:
$$4n+12=4(n+n_3+3n_4)+8c$$
and then $2c=3-n_3-3n_4$.
\\ If $c=0$ then there are two cases:
\begin{itemize}
\item  $n_4=1, n_3=0$, so $n_1=2$ and $G=\gamma_{p,k_1,k_2}$ with $n=p+k_1+k_2$ 
%

By theorem \ref{gamma}, $G$ cannot be cospectral with $\gamma_{p,k_1,k_2}$; this case is impossible.
\item $n_4=0,n_3=3$, so $n_1=3$ and $G\in\mathcal{G}_1\cup\mathcal{G}_2\cup\mathcal{G}_3$. But by theorems \ref{L4pascospG1}, \ref{L4pascospG2} and \ref{L4pascosp} this is impossible.
\end{itemize}

As a result $c\neq0$ and $G$ has a $4$-cycle.
\end{demo}

Following the proof of theorem \ref{th_odd} for odd lollipop, we can now state:

\begin{theorem}
The lolipop $L(4,k)$ is determined by its spectrum.
\end{theorem}

\section{Conclusion}
In this paper we give a way to count closed walks, which is relevant to show that two graphs cannot be cospectral. 

That provides a new approach to show that the odd lollipops are determined by their spectrum  and following this same idea we have proved that even lollipops are also determined by their spectrum. However this is far to be as simple as the odd case and we had to develop several tools to show the non-cospectrality of two given graphs. The most difficult case, as it was noted in \cite{lollipop,developments}, is for the lollipops $L(4,k)$ where connectivity and presence of a $4$-cycle are quite long to establish.

\newpage
\appendix
\section{Appendix}\label{annexe}

\subsection{Counting covering closed walks}

\begin{center}
\begin{table}[h!]
\begin{center}
$\begin{array}{|c|c|c|c|}
M & w_6(M) & w_8(M) & w_{10}(M) \\
\hline
P_2 &   2  &  2     &    2      \\
P_3 &   12  &  28   &    60      \\
P_4 &   6  &  32   &    120      \\
P_5 &   0  &  8    &    60     \\
P_6 &   0  &  0     &    10      \\
C_4 &   48  &  264     &  1320       \\
C_6 &   12  &       &          \\
C_8 &   0  &  16     &         \\
C_{10} &   0  &  0     &    20      \\
S_{1,1,1}& 12 & 72    &    300      \\
S_{1,1,2} & 0 & 16    &    140      \\
S_{1,1,3} & 0 & 0    &    20     \\
S_{1,2,2} & 0 & 0    &    20      \\
L(4,1) &   12  &  112     & 840      \\
L(4,2) &   0  &  16   &    180     \\
L(4,3) &   0  &  0    &    20      \\
\hline
\end{array}$
\end{center}
\caption{Number of covering closed walks on a given graph.}
\label{tab_marches}
\end{table}
\end{center}

\subsection{Proof of theorem \ref{vpP-2} }
First, we notice the following relations which will be useful to prove lemmas  \ref{33p} and \ref{dist2} and whose proof is straightfoward  by induction on $p$.
\begin{equation}\label{beta}
\forall p>0,\ Q_{P_p}(\alpha)>\beta Q_{P_{p-1}}(\alpha)
\end{equation}
where $\alpha=\sqrt{2+2\sqrt{2}}$ and
$\beta=\frac{\sqrt{2}}{2}\alpha$. Obviously equation \ref{beta} is true if we replace $\beta$ by $\beta'\leq\beta$.

\begin{lemme}\label{vpP2}
$\lambda_1(P(p_1,p_2,p_3))>2$.
\end{lemme}

\begin{demo}
On one hand $\lambda_1(P(0,1,1))>2$ and $\lambda_1(P(1,1,1))>2$.
On the other hand, if there exists $p_i\geq 2$ (we assume $p_3\geq 2$)  then  the lollipop $L(p_1+p_2+2,1)$ is an induced subgraph of $P(p_1,p_2,p_3)$. Since $\lambda_1(L(p_1+p_2+2,1))>2$ (theorem \ref{vp2}) the interlacing theorem gives the result.
\end{demo}

Applying theorem \ref{del_vertex} to a vertex of degree $3$ of $P(p_1,p_2,p_3)$ we can get the following expression of the characteristic polynomial of $P(p_1,p_2,p_3),\ p_i>0$ which will be useful for the next results.

\begin{equation}
\begin{split}
Q_{P(p_1,p_2,p_3)}(X) = XQ_{S_{p_1,p_2,p_3}}(X)-Q_{S_{p_1-1,p_2,p_3}}(X)-Q_{S_{p_1,p_2-1,p_3}}(X)\\
-Q_{S_{p_1,p_2,p_3-1}}(X)-2Q_{P_{p_1}}(X)-2Q_{P_{p_2}}(X)-2Q_{P_{p_3}}(X)
\end{split}
\label{eq_pol_car_P}
\end{equation}

where

\begin{equation}
\begin{split}
Q_{S_{a,b,c}}(X)=XQ_{P_a}(X)Q_{P_b}(X)Q_{P_c}(X)-Q_{P_{a-1}}(X)Q_{P_b}(X)Q_{P_c}(X)\\
-Q_{P_a}(X)Q_{P_{b-1}}(X)Q_{P_c}(X) -Q_{P_a}(X)Q_{P_b}(X)Q_{P_{c-1}}(X)
\end{split}
\label{eq_pol_car_S}
\end{equation}

\begin{lemme}\label{33p}
If $p_1\leq 3$, $p_2\leq 3$ then $\forall p\in\n :\ \lambda_1(P(p_1,p_2,p))>\sqrt{2+2\sqrt{2}}$.
\end{lemme}

\begin{demo}
According to theorem \ref{homeo} it is sufficient to prove the result for $p_1=3$, $p_2=3$. Let $\alpha=\sqrt{2+2\sqrt{2}}$.
We shall show that $Q_{P(3,3,p)}(\alpha)<0$.
Using equations (\ref{eq_pol_car_P}) and (\ref{eq_pol_car_S}) and $Q_{P_{p-2}}(X)=XQ_{P_{p-1}}(X)-Q_{P_{p}}(X)$, $Q_{P_{2}}(X)=X^2-1$,  $Q_{P_{3}}(X)=X^3-2X$ and $Q_{P_{4}}(X)=X^4-3X^2+1$, we get:

\begin{eqnarray*}
Q_{P(3,3,p)}(X)&=&Q_{P_{p}}(X)\big( X^8-9 X^6+ 24 X^4-20 X^2 \big)\\
&&+Q_{P_{p-1}}(X)\big( -X^7+8 X^5-16 X^3+ 8 X \big)-4(X^3-2X)\\
\end{eqnarray*}
so 
\begin{eqnarray*}
Q_{P(3,3,p)}(\alpha)&=&(16-16\sqrt{2})Q_{P_{p}}(\alpha)+\alpha(16-8\sqrt{2})Q_{P_{p-1}}(\alpha)-8\sqrt{2}\alpha\\
&=&\left(-16+16\sqrt{2}\right)\left(-Q_{P_{p}}(\alpha)+\frac{\alpha}{\sqrt{2}}Q_{P_{p-1}}(\alpha)\right)-8\sqrt{2}\alpha<0 \textrm{ (by eq.(\ref{beta}) ) }
\end{eqnarray*}

As a result  $\lambda_1(P(3,3,p))>\alpha$.
\end{demo}


\begin{lemme}\label{24p}
If $p_1\leq 2$, $p_2\leq 4$ then $\forall p\in\n :\ \lambda_1(P(p_1,p_2,p))>2.2>\sqrt{2+2\sqrt{2}}$
\end{lemme}

%
%

\begin{demo}
\textit{Mutatis mutandis} the proof is the same as the one of lemma \ref{33p}.
\end{demo}

\begin{lemme}\label{P1pq}
For $p_2,p_3>0$, $p_1\in\{0,1\}$ we have $\lambda_1(P(p_1,p_2,p_3))>\sqrt{2+2\sqrt{2}}$. 
\end{lemme}

\begin{demo}
Let $\alpha=\sqrt{2+2\sqrt{2}}$.
According to theorem \ref{homeo} it is sufficient to prove the result for $P(1,p,p)$ where $p=\max(p_2,p_3)$.
Applying theorem \ref{del_vertex} to a vertex at distance one of the two vertices of degree $3$ we have:

\begin{eqnarray*}
Q_{P(1,p,p)}(X)&=& XQ_{C_{2p+2}}(X)-2Q_{P_{2p+1}}(X)-4Q_{P_p}(X)\\
&=& X\big( XQ_{P_{2p+1}}(X) -2Q_{P_{2p}}(X)-2  \big)-2Q_{P_{2p+1}}(X)-4Q_{P_p}(X)\\
&=& (X^2-2)Q_{P_{2p+1}}(X) -2XQ_{P_{2p}}(X)-4Q_{P_p}(X)-2X
\end{eqnarray*}

But (theorem \ref{del_vertex} applied to a vertex at distance $p$ of a pendant vertex in the graphs $P_{2p+1}$ and $P_{2p}$ ):
$$Q_{P_{2p+1}}(X)=X\big(Q_{P_p}^2(X)\big)-2Q_{P_p}(X)Q_{P_{p-1}}(X)$$ and $$Q_{P_{2p}}(X)=XQ_{P_{p}}(X)Q_{P_{p-1}}(X)-Q_{P_{p-1}}^2(X)-Q_{P_{p}}(X)Q_{P_{p-2}}(X)$$ 

So
\begin{eqnarray*}
Q_{P(1,p,p)}(X) &=& \big(X^2-2\big)\big(XQ_{P_{p}}^2(X)-2XQ_{P_{p}}(X)Q_{P_{p-1}}(X)\big)\\
&&-2X\big(X Q_{P_{p}}(X)Q_{P_{p-1}}(X)-Q_{P_{p-1}}^2(X)-Q_{P_{p}}(X)Q_{P_{p-2}}(X)\big)\\
&&-4Q_{P_{p}}(X)-2X\\
&=&Q_{P_{p}}(X)\big( (X^3-4X)Q_{P_{p}}(X)+(4-2X^2) Q_{P_{p-1}}(X) -4\big)\\
&&+2XQ_{P_{p-1}}^2(X)-2X
\end{eqnarray*}

Using $Q_{P_{p}}(\alpha)>\beta Q_{P_{p-1}}(\alpha)$ (equation (\ref{beta})), we get

\begin{eqnarray*}
Q_{P(1,p,p)}(\alpha) &<&Q_{P_{p}}(\alpha)\big( (\alpha^3-4\alpha)Q_{P_{p}}(\alpha)+(4-2\alpha^2+\frac{2\alpha}{\beta}) Q_{P_{p-1}}(\alpha) -4\big)-2\alpha
\end{eqnarray*}

%
%

we then notice that $\frac{4-2\alpha^2+\frac{2\alpha}{\beta}}{-\alpha^3+4\alpha}=\beta$ and by equation (\ref{beta}) we have $Q_{P(1,p,p)}(\alpha)<0$.
\end{demo}

\begin{lemme}\label{dist2}
Given $P(2,p_2,p_3)$ with $p_3\geq 3$,  denote by $u$ and $v$ the two vertices of degree $3$.  Let $y$ be a vertex at distance $2$ from $u$ and at distance greater than or equal to 2 from $v$, we define  $\tilde{P}(2,p_2,p_3)$ as the graph obtained by adding to $P(2,p_2,p_3)$ a pendant vertex $x$ to $y$. We have $\lambda_1(\tilde{P}(2,p_2,p_3))>\sqrt{2+2\sqrt{2}}$.
\end{lemme}

\begin{demo}
Let $\alpha=\sqrt{2+2\sqrt{2}}$. By theorem \ref{homeo} it is sufficient to prove the result for $p_2=p_3=p=\max\{p_2,p_3\}$. The aim of the proof is to show that $Q_{\tilde{P}(2,p,p)}(\alpha)<0$.
The following equations will be useful:
$$Q_{S_{2,a,b}}(\alpha)=(\alpha^2-1)Q_{P_{a+b+1}}(\alpha)-\alpha Q_{P_a}(\alpha)Q_{P_b}(\alpha)$$
$$Q_{P_{2p+1}}(\alpha)=\alpha Q_{P_p}^2(\alpha)-2Q_{P_p}(\alpha)Q_{P_{p-1}}(\alpha)$$
$$Q_{P_{2p}}(\alpha)=Q_{P_p}^2(\alpha)-Q_{P_{p-1}}^2(\alpha)$$
$$Q_{P_{2p-1}}(\alpha)=\alpha Q_{P_{2p}}(\alpha)-Q_{P_{2p+1}}(\alpha)=-\alpha Q_{P_{p-1}}^2(\alpha)+2Q_{P_p}(\alpha)Q_{P_{p-1}}(\alpha)$$
and we deduce
$$Q_{S_{2,p,p}}(\alpha)=(\alpha^3-2\alpha)Q_{P_p}^2(\alpha)-2(\alpha^2-1)Q_{P_p}(\alpha)Q_{P_{p-1}}(\alpha)$$
$$Q_{S_{2,p,p-1}}(\alpha)=(\alpha^2-1)Q_{P_p}^2(\alpha)-(\alpha^2-1)Q_{P_{p-1}}^2(\alpha)-\alpha Q_{P_p}(\alpha)Q_{P_{p-1}}(\alpha)$$
$$Q_{S_{2,p+1,p-1}}(\alpha)=(\alpha^3-\alpha)Q_{P_p}^2(\alpha)+(-3\alpha^2+2)Q_{P_p}(\alpha)Q_{P_{p-1}}(\alpha)+\alpha Q_{P_{p-1}}^2(\alpha)$$

Theorem \ref{del_vertex} gives
$$Q_{\tilde{P}(2,p,p)}(\alpha)=\alpha Q_{P(2,p,p)}(\alpha)-Q_H(\alpha)$$
where $H=\tilde{P}(2,p,p)\backslash\{x,y\} $
Equation \ref{eq_pol_car_P} gives 
$$Q_{P(2,p,p)}(\alpha)=\alpha Q_{S_{2,p,p}}(\alpha)-Q_{S_{1,p,p}}(\alpha)-2Q_{S_{2,p,p-1}}(\alpha)-4Q_{P_p}(\alpha)-2Q_{P_2}(\alpha)$$
but $Q_{S_{1,p,p}}(\alpha)=\frac{1}{\alpha}(Q_{S_{2,p,p}}(\alpha)+Q_{S_{0,p,p}}(\alpha))$ so

\begin{eqnarray*}
\alpha Q_{P(2,p,p)}(\alpha)&=&(\alpha^2-1) Q_{S_{2,p,p}}(\alpha)-2\alpha Q_{S_{2,p,p-1}}(\alpha)- Q_{P_{2p+1}}(\alpha)\\
&&-4\alpha Q_{P_p}(\alpha)-2\alpha Q_{P_2}(\alpha)\\
&=&(\alpha^5-5\alpha^3+3\alpha)Q_{P_p}^2(\alpha)+(2\alpha^3-2\alpha)Q_{P_{p-1}}^2(\alpha)\\
&&+(-2\alpha^4+6\alpha^2)Q_{P_p}(\alpha)Q_{P_{p-1}}(\alpha)-4\alpha Q_{P_p}(\alpha)-2\alpha Q_{P_2}(\alpha)\\
\end{eqnarray*}

Theorem \ref{del_vertex} gives
$$Q_H(\alpha)=\alpha^2Q_{S_{2,p,p-2}}(\alpha)-\alpha Q_{S_{1,p,p-2}}(\alpha)-\alpha Q_{S_{2,p-1,p-2}}(\alpha)-Q_{S_{2,p,p-2}}(\alpha)-2\alpha Q_{P_{p-2}}(\alpha)$$
but $Q_{S_{2,p,p-2}}(\alpha)=\alpha Q_{S_{2,p,p-1}}(\alpha)-Q_{S_{2,p,p}}(\alpha)$, $\alpha Q_{S_{1,p,p-2}}(\alpha)=Q_{S_{2,p,p-2}}(\alpha)+Q_{S_{0,p,p-2}}(\alpha)$ and $Q_{S_{2,p-1,p-2}}(\alpha)=(\alpha^2-1) Q_{S_{2,p,p-1}}(\alpha)-\alpha Q_{S_{2,p-1,p+1}}(\alpha)$ so

\begin{eqnarray*}
Q_H(\alpha)&=&-\alpha Q_{S_{2,p,p-1}}(\alpha)-(\alpha^2-2)Q_{S_{2,p,p}}(\alpha)+\alpha^2Q_{S_{2,p-1,p+1}}(\alpha)\\
&&-Q_{P_{2p-1}}(\alpha)-2\alpha Q_{P_{p-2}}(\alpha)\\
&=&(2\alpha^3-3\alpha)Q_{P_p}^2(\alpha)+(2\alpha^3)Q_{P_{p-1}}^2(\alpha)\\
&&+(-\alpha^4-3\alpha^2+2)Q_{P_p}(\alpha)Q_{P_{p-1}}(\alpha)-2\alpha Q_{P_{p-2}}(\alpha)\\
\end{eqnarray*}

So we have:
\begin{eqnarray*}
Q_{\tilde{P}(2,p,p)}(\alpha)&= &(\alpha^5-7\alpha^3+6\alpha)Q_{P_p}^2(\alpha)-2\alpha Q_{P_{p-1}}^2(\alpha)
\\&&+( -\alpha^4+9\alpha^2-2 )Q_{P_p}(\alpha)Q_{P_{p-1}}(\alpha)
\\&&+2\alpha Q_{P_{p-2}}(\alpha)-4\alpha Q_{P_p}(\alpha)-2\alpha Q_{P_2}(\alpha)
\end{eqnarray*}

Equation (\ref{beta}) gives $2\alpha Q_{P_{p-2}}(\alpha)-4\alpha Q_{P_p}(\alpha)<0$. \\
Lets us show  that $x Q_{P_p}^2(\alpha)+y Q_{P_{p-1}}^2(\alpha)+zQ_{P_p}(\alpha)Q_{P_{p-1}}(\alpha)<0$ 
with $x=\alpha^5-7\alpha^3+6\alpha$, $y=-2\alpha$, $z= -\alpha^4+9\alpha^2-2$. Note that $\frac{y}{z+\beta x}=-\beta$, where $\beta$ is defined in equation (\ref{beta}).


\begin{eqnarray*}
x Q_{P_{p}}^2(\alpha)+y Q_{P_{p-1}}^2(\alpha)+zQ_{P_{p}}(\alpha)Q_{P_{p}}(\alpha)=\\
Q_{P_{p}}(\alpha)\left(xQ_{P_{p}}(\alpha)-\beta xQ_{P_{p-1}}(\alpha)\right)+Q_{P_{p-1}}(\alpha)\left((z+\beta x)Q_{P_{p}}(\alpha)+yQ_{P_{p-1}}(\alpha)\right)=\\
Q_{P_{p}}(\alpha)x\left(Q_{P_{p}}(\alpha)-\beta Q_{P_{p-1}}(\alpha)\right)+Q_{P_{p-1}}(\alpha)(z+\beta x)\left(Q_{P_{p}}(\alpha)-\beta Q_{P_{p-1}}(\alpha)\right)<0\\
\end{eqnarray*}

because $\frac{(z+\beta x)\left(Q_{P_{p}}(\alpha)-\beta Q_{P_{p-1}}(\alpha)\right)}{-x\left(Q_{P_{p}}(\alpha)-\beta Q_{P_{p-1}}(\alpha)\right)}=\frac{z+\beta x}{-x}<\beta$ and we use equation (\ref{beta}).
\end{demo}

\begin{lemme}\label{dist1}
Given $P(2,p_2,p_3)$ with $p_3\geq 3$,  denote by $u$ and $v$ the two vertices of degree $3$.  Let $y$ be a vertex at distance $1$ from $u$ and at distance greater than or equal to 1 from $v$, we denote by $\hat{P}(2,p_2,p_3)$  the graph obtained by adding to $P(2,p_2,p_3)$ a pendant vertex $x$ to $y$. We have $\lambda_1(\hat{P}(2,p_2,p_3))>\sqrt{2+2\sqrt{2}}$.
\end{lemme}

\begin{demo}
A direct consequence of theorem \ref{homeo} and lemma \ref{dist2}.
\end{demo}


\begin{theorem_annexe*}
For $p_1,p_2,p_3>0$, $P(p_1,p_2,p_3)$ cannot be an induced subgraph of a graph cospectral with an even lollipop.
\end{theorem_annexe*}

\begin{demo}
Without loss of generality we  assume $p_1\leq p_2\leq p_3$.
In order to lead a proof by contradiction, let $P(p_1,p_2,p_3)$ be an induced subgraph of $G$ cospectral with an even lollipop.  As $G$ is bipartite, $P(p_1,p_2,p_3)$ doesn't have odd cycles and the $p_i$'s are all odd or all even. 
Using equations (\ref{eq_pol_car_P}) and property \ref{P(2)} we obtain:
$$Q_{P(p_1,p_2,p_3)}(2)=p_1p_2p_3-p_1p_2-p_1p_3-p_2p_3-3p_1-3p_2-3p_3-5$$ 
i) First assume that $p_1,p_2,p_2$ are odd.\\\\
By lemma \ref{P1pq} we have $p_1\geq3$ and by lemma \ref{33p} we have $p_2\geq5$.
\begin{itemize}
\item If $p_1=3$ and $p_2=5$ then $Q_{P(p_1,p_2,p_3)}(2)=4p_3-44 \geq 0$ if $p_3\geq11$
\item If $p_1=3$ and $p_2\geq 7$ then $Q_{P(p_1,p_2,p_3)}(2)\geq 2p_3-14 \geq 0$ (because  $p_3\geq p_2\geq7$)
\item If $5\leq p_1\leq p_2\leq p_3$ then $Q_{P(p_1,p_2,p_3)}(2)\geq p_3-5\geq 0$.
\end{itemize}
$Q_{P(p_1,p_2,p_3)}(2)\geq 0$ implies that $P(p_1,p_2,p_3)$ has two eigenvalues greater than or equal to $2$ (we already know by lemma \ref{vpP2} that $P(p_1,p_2,p_3)$ has at least one eigenvalue strictly greater than $2$) and since a lollipop has only one eigenvalue greater than $2$ (theorem \ref{vp2}), the interlacing theorem provides a contradiction except when $p_1=3$, $p_2=5$ and $p_3\in\{5,7,9\}$.

Assume now that $p_1=3$, $p_2=5$ and $p_3\in\{5,7,9\}$. According to table \ref{tab_eig_P3}, $\lambda_1(P(3,5,p_3))>2.17$ and so $P(3,5,p_3)$   cannot be an induced subgraph of a graph cospectral with  $L(p,k)$ for $p\geq 6$ (theorem \ref{2.17}).
Moreover $P(3,5,p_3)$ cannot be a connected component of a graph cospectral with $L(4,k)$ because  $\lambda_1(P(3,5,p_3))<2.195$ while  $\lambda_1(L(4,k))\geq\lambda_1(L(4,5))> 2.195$ when $k\geq5$.
So there is a new vertex $x$ adjacent to one vertex $y$ of $P(3,5,p_3)$ (and only one because otherwise there exists $r,s\in\n$ such that $P(1,r,s)$ is an induced subgraph of $G$ which is impossible by lemma \ref{P1pq}). Let $H$ be the subgraph induced by $P(3,5,p_3)$ and $x$, denote by $u$ and $v$ the two vertices of degree $3$ in $P(3,5,p_3)$. 
\begin{enumerate}
\item If $y=u$ or $y=v$ then the graph $T$ drawn on figure  \ref{usb} is an induced subgraph of $H$ and $\lambda_1(T)\geq 2.20>\sqrt{2+2\sqrt{2}}>\lambda_1(L(p,k))$ and $H$ cannot be an induced subgraph of $G$.

\begin{figure}[htbp]
\begin{center}
\includegraphics[scale=0.5]{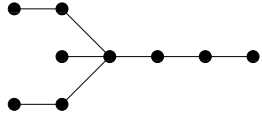}
\caption{Tree $T$ whose spectral radius is greater than 2.20}
\label{usb}
\end{center}
\end{figure}

\item If $\min\{d(y,u), d(y,v)\}\geq5$ the disjoint union of a cycle and $S_{1,3,3}$ is an induced subgraph of $H$ with twice the eigenvalue $2$, so $H$ cannot be an induced subgraph of $G$ (by the interlacing theorem and theorem \ref{vp2}). 

\item The cases where $1\leq \min\{d(y,u),d(y,v)\}\leq4$ are summed up in table \ref{tab_P+v}. For all these cases $H$ cannot be an induced subgraph of  $G$ because either $H$ has two eigenvalues greater than $2$ or $H$ has a spectral radius greater than $\sqrt{2+2\sqrt{2}}$.
\end{enumerate}
 As a result $P(p_1,p_2,p_3)$ with $p_i$'s odd  cannot be an induced subgraph of $G$.
\\
\\
\\ii) We now assume that $p_1,p_2,p_3$ are even.\\\\
By lemma \ref{P1pq} we have $p_1\geq 2$. 
\begin{itemize}
\item If $p_1=2$ and $p_2\leq 4$ then by lemma \ref{24p} $P(p_1,p_2,p_3)$ cannot be an induced subgraph of $G$.
\item If $p_1=2$ and $p_2=6$ then $Q_{P(p_1,p_2,p_3)}(2)=p_3-41 \geq 0$ if $p_3\geq42$
\item If $p_1=2$ and $p_2=8$ then $Q_{P(p_1,p_2,p_3)}(2)=3p_3-51 \geq 0$ if $p_3\geq18$
\item If $p_1=2$ and $p_2=10$ then $Q_{P(p_1,p_2,p_3)}(2)=5p_3-61 \geq 0$ if $p_3\geq14$
\item If $p_1=2$ and $p_2\geq 12$ then $Q_{P(p_1,p_2,p_3)}(2)\geq 2p_3-11 \geq 0$ (because $p_3\geq 12$)
\item If $p_1=4$ and $p_2=4$ then $Q_{P(p_1,p_2,p_3)}(2)=5p_3-45 \geq 0$ if $p_3\geq10$
\item If $p_1=4$ and $p_2\geq 6$ then $Q_{P(p_1,p_2,p_3)}(2)\geq4p_3-17 \geq 0$ (because $p_3\geq6$)
\item If $6\leq p_1\leq p_2\leq p_3$ then $Q_{P(p_1,p_2,p_3)}(2)\geq 9p_3-5\geq 0$.
\end{itemize}

As in the proof of the odd case, if $Q_{P(p_1,p_2,p_3)}(2)\geq 0$ then $P(p_1,p_2,p_3)$ has two eigenvalues greater than or equal to $2$ and cannot be an induced subgraph of $G$. We are now going to study the remaining cases for $p_1=2$ and $p_1=4$. 
\\\textbf{First case $p_1=2$ : }\\ 
The only unsolved cases we are going to consider here are for $p_2\in\{6,8,10\}$ with the corresponding constraints on $p_3$.
According to table \ref{tab_eig_P2}, the spectral radius of these remaining cases is greater than $2.17$ and so the corresponding graphs cannot be an induced subgraph of a graph cospectral with $L(p,k)$, $p\geq 6$.
%
As it was detailed in the proof of the odd case, none of these graphs is a connected component of a graph cospectral with $L(4,k)$ and so there is a new vertex $x$ adjacent to one and only one vertex $y$ of $P(2,p_2,p_3)$. Let $H$ be the subgraph induced by $P(2,p_2,p_3)$ and $x$. With the same notations and arguments as for the odd case, $H$ cannot be   an induced subgraph of $G$ when  $\min\{d(y,u), d(y,v)\}\geq5$ or $y=u$ or $y=v$. Moreover if $\min\{d(y,u), d(y,v)\}\leq 2$ then by lemmas \ref{dist1} and \ref{dist2}, $\lambda_1(H)>\sqrt{2+2\sqrt{2}}$ so $H$ cannot be an induced subgraph of a lollipop. We are now going to examine the two last tricky cases: $\min\{d(y,u),\ d(y,v)\}=3$ and $\min\{d(y,u),\ d(y,v)\}=4$.

\begin{itemize}

\item  If $\min\{d(y,u),\ d(y,v)\}=3$, we can assume that $d(y,v)=3$. Let $\{b,c\}=\{p_2,p_3\}$ such that $y$ is a vertex belonging to a path of length $c+1$ of  $P(2,b,c)$ between $u$ and $v$. Then applying theorem \ref{del_pendant} to $x$ we get $Q_H(X)=XQ_{P(2,b,c)}(X)-Q_{P(2,b,c)\backslash\{y\}}(X)$ and applying theorem \ref{del_vertex} to  $v$ we have:
\begin{eqnarray*}
Q_{P(2,b,c)\backslash\{y\}}(X)&=&XQ_{P_2}(X)Q_{S_{2,b,c-3}}(X)-Q_{P_2}(X)Q_{S_{1,b,c-3}}(X)\\
&&-Q_{P_2}(X)Q_{S_{2,b-1,c-3}}(X)-Q_{P_1}(X)Q_{S_{2,b,c-3}}(X)\\
&&-2Q_{P_2}(X)Q_{P_{c-3}}(X)
\end{eqnarray*} 
Using equation (\ref{eq_pol_car_P}) and property  \ref{P(2)} which gives the value in $2$ of the characteristic polynomials of paths and $T$-shape trees  we obtain:
$$Q_H(2)=bc-5b+4c-56$$
\begin{itemize}
\item If $b\leq c$
\begin{itemize}
\item If $b=6$ (so $c\geq 6$) then $Q_H(2)=10c-86$ so if $c\geq 10$, $H$ has two eigenvalues greater than $2$  and cannot be and induced subgraph of $G$. Otherwise for $c=8$ we check that $\lambda_1(H)\sim 2.2050 >\sqrt{2+2\sqrt{2}}$ and so $H$ cannot be an induced subgraph of $G$ for $c\leq 8$.
\item If $c\geq b\geq 8$ then  $Q_H(2) \geq 7c-56\geq 0$ and $H$ has two eigenvalues greater than $2$ and cannot be an induced subgraph of $G$.
\end{itemize} 
\item If $b\geq c$
\begin{itemize}
\item If $c=6$ then $Q_H(2)=b-32$ so if $b\geq 32$ then $H$ has two eigenvalues greater than $2$ and cannot be an induced subgraph of $G$. Otherwise we check that for $b=30$ we have $\lambda_1(H)\sim 2.2071>\sqrt{2+2\sqrt{2}}$ and so $H$ cannot be an induced subgraph of $G$ for $b\leq 30$.
\item If $8\leq c\leq b$ then $Q_H(2)\geq 4c-32\geq 0$ and $H$  has two eigenvalues greater than $2$ and cannot be an induced subgraph of $G$.
\end{itemize}
\end{itemize}

\item If $\min\{d(y,u),d(y,v)\}=4$, note that $c\geq 8$  (otherwise $y$ is at distance less than $4$ from $u$ or $v$). In the same way as previously we compute $Q_H(2)$:
$Q_H(2)=b+9c-86$.
\begin{itemize}
\item If $b\leq c$
\begin{itemize}
\item If $b=6$ then $Q_H(2)=9c-80$. So if $c\geq 10$ then $H$  has two eigenvalues greater than $2$ and cannot be an induced subgraph of $G$. Otherwise we check that for $c=8$ we have $\lambda_1(H)\sim 2.2014>\sqrt{2+2\sqrt{2}}$.
\item If $b=8$ then $Q_H(2)=9c-78$. So if $c\geq 10$ then $H$  has two eigenvalues greater than $2$ and cannot be an induced subgraph of $G$. The case $c=b=8$ is considered further in the proof.
\item If  $10\leq b\leq c$ then $Q_H(2)>0$ and $H$  has two eigenvalues greater than $2$ and cannot be an induced subgraph of $G$.
\end{itemize} 
\item If $b\geq c$ 
\begin{itemize}
\item If $c=8$ then $Q_H(2)=b-14$. So if $b\geq 14$ then $H$  has two eigenvalues greater than $2$ and cannot be an induced subgraph of $G$. Otherwise we check for $c=8$ and $8\leq b\leq 12$ that $\lambda_1(H)<2.196$ so $H$ cannot be a connected component of  $G$ because for $k\geq 6$ $\lambda_1(L(4,k))\geq\lambda_1(L(4,6))> 2.196$. And so there is a new vertex $x'$ adjacent to a vertex $y'$ of $H$. Let $H'$ be the graph induced by $H$ and $x'$.
\begin{itemize}
\item If $y'=y$ then $x'$ is not adjacent to another vertex of $P(2,a,b)$ otherwise there exists $r,s\in\n$ such that $P(1,r,s)$ is an induced subgraph of $G$ which is impossible by lemma \ref{P1pq} and $x'$ is not adjacent to $x$ otherwise $G$  contains a triangle (impossible because $G$ is bipartite). Hence $x'$ is a pendant in $H'$. The graph $H'$ then contains $C_q\cup S_4$ (for $q\geq3$) as an induced subgraph and so has two eigenvalues greater than $2$ which is impossible.
 \item Assume that $y'=x$. If $x'$ is adjacent to another vertex $z$ of $H$ distinct from $y'$ and $y$, then by the previous cases we necessarily have  $\min\{d(z,u),d(z,v)\}=4$. Either the graph $S_{1,3,3}\cup S_{2,2,2}$ or $C_4\cup C_q$ is an induced subgraph of $H'$ and has two eigenvalues greater than $2$ and cannot be an induced subgraph of $G$.
\item If $y'\neq y$ and $y'\neq x$  then by the previous cases we necessarily have $\min\{d(y',u), d(y',v)\}=4$. 
\\If $x'$ is adjacent to another vertex $z$ in $H$ distinct from $y'$ and $y$ then by the previous cases we necessarily have  $\min\{d(z,u),d(z,v)\}=4$ and either $S_{2,2,2}\cup S_{1,2,5}$ or $C_r\cup C_s$ is an an induced subgraph of $H'$  and has two eigenvalues greater than $2$ and cannot be an induced subgraph of $G$.
\\If $x'$ is not adjacent to another vertex of $H$ then the graph $T_n\cup C_q$ or the graph $S_{1,3,3}\cup S_{1,3,3}$ is an induced subgraph of $H'$ and has two eigenvalues greater than $2$ and cannot be an induced subgraph of $G$.
\end{itemize}
\item  If  $10\leq b\leq c$ then $Q_H(2)>0$ and $H$  has two eigenvalues greater than $2$ and cannot be an induced subgraph of $G$.
\end{itemize}

\end{itemize}

\textbf{Second case: $p_1=4$. }\\ We  have $p_2=4$ and $p_3\in\{4,6,8\}$. 

According to table \ref{tab_eig_P4}, $\lambda_1(P(4,4,p_3))> 2.17$ and so $P(4,4,p_3)$   cannot be an induced subgraph of a graph cospectral with  $L(p,k)$ for $p\geq 6$ (theorem \ref{2.17}).
Moreover $\lambda_1(P(4,4,4))>\sqrt{2+2\sqrt{2}}$ and $P(4,4,4)$ cannot be an induced subgraph  of a graph cospectral with  $L(4,k)$.
When $p_3\in\{6,8\}$, $P(4,4,p_3)$ cannot be a connected component of a graph cospectral with $L(4,k)$ because  $\lambda_1(P(4,4,p_3))<2.1854$ while  $\lambda_1(L(4,k))\geq\lambda_1(L(4,3))> 2.1888$ when $k\geq3$.
So there is a new vertex $x$ adjacent to one vertex $y$ of $P(4,4,p_3)$ (and only one because otherwise there exists $r,s\in\n$ such that $P(1,r,s)$ is an induced subgraph of $G$ which is impossible by lemma \ref{P1pq}). Let $H$ be the subgraph induced by $P(4,4,p_3)$ and $x$, these graphs $H$ are summed up in table \ref{tab_P+v2} which shows that that $H$ cannot be an induced subgraph of  $G$ because either $H$ has two eigenvalues greater than $2$ or $H$ has a spectral radius greater than $\sqrt{2+2\sqrt{2}}$.
%

\end{itemize}
\end{demo}

\subsection{Tables of some graphs eigenvalues}

\begin{center}
\begin{table}[htbp]
$\begin{array}{cc}

\includegraphics[scale=0.7]{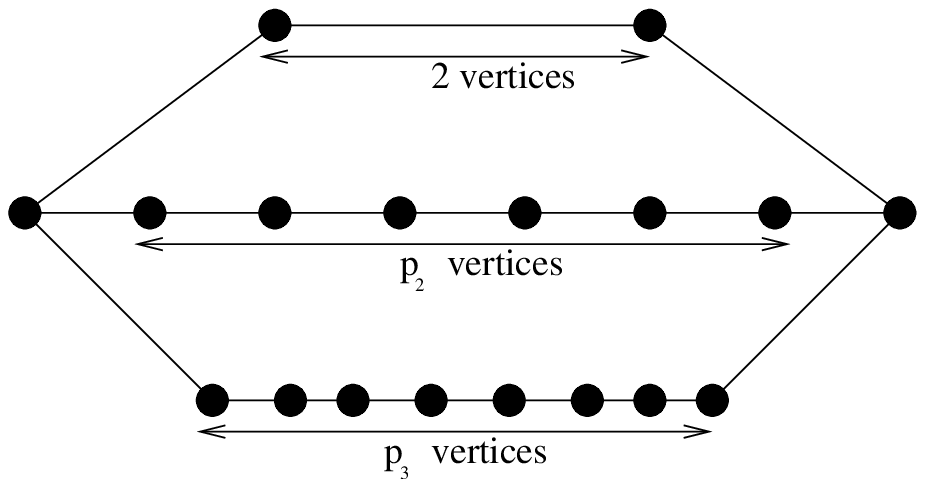}

&
\begin{array}{|c||c|c|c|}
\hline
\textrm{\backslashbox{$p_3$}{$p_2$}} &  6  &  8  &  10  \\
\hline
\hline
6        &  \begin{array}{c}2.1987\\1.9122 \end{array}    &   \begin{array}{c}2.1921\\1.9426 \end{array}   &   \begin{array}{c}2.1891\\1.19604 \end{array}        \\
\hline

8        &  \begin{array}{c}2.1921\\1.9426 \end{array}    &  \begin{array}{c}2.1853\\1.9666 \end{array}    &      \begin{array}{c}2.1822\\1.19805 \end{array}      \\
\hline

10        &  \begin{array}{c}2.1891\\1.9604 \end{array}    &  \begin{array}{c}2.1822\\1.9805 \end{array}     &     \begin{array}{c}2.1790\\1.9922 \end{array}        \\
\hline

12       &  \begin{array}{c}2.1878\\1.9716 \end{array}   &   \begin{array}{c}2.1808\\1.9891 \end{array}    &      \begin{array}{c}2.1776\\1.9994 \end{array}      \\
\hline

14        &  \begin{array}{c}2.1872\\1.9790 \end{array}   &   \begin{array}{c}2.1802\\1.9947 \end{array}    &     \begin{array}{c}2.1770\\2.0041 \end{array}       \\
\hline

16        &  \begin{array}{c}2.1870\\1.9842 \end{array}    &  \begin{array}{c}2.1800\\1.9986 \end{array}     &        \begin{array}{c}2.1767\\2.0072 \end{array}    \\
\hline

40        &   \begin{array}{c}2.1868\\1.9999 \end{array}     &     &          \\
\hline
\end{array}
\end{array}$

\caption{The two largest eigenvalues of $P(2,p_2,p_3)$ with a $4$ decimal place accuracy.}
\label{tab_eig_P2}
\end{table}
\end{center}

\begin{center}
\begin{table}[htbp]
$\begin{array}{cc}

\includegraphics[scale=0.7]{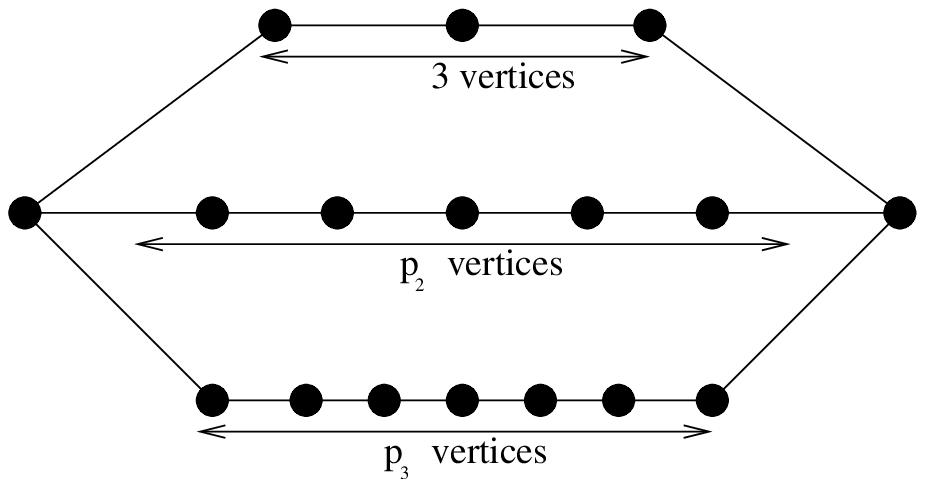}

&
\begin{array}{|c||c|c|}
\hline
\textrm{\backslashbox{$p_3$}{$p_2$}} &  5  &  7    \\
\hline\hline
5        &   \begin{array}{c}2.1940\\1.9319 \end{array}  &    \begin{array}{c}2.1847\\1.9696 \end{array}    \\
\hline
7        &    \begin{array}{c}2.1847\\1.9696 \end{array}   &    \begin{array}{c}2.1753\\2.0000 \end{array}   \\
\hline
9        &   \begin{array}{c}2.1804\\1.9890 \end{array}  &  \begin{array}{c}2.1709\\2.0153 \end{array}     \\
\hline
11       &   \begin{array}{c}2.1785\\2.0000 \end{array}  &   \begin{array}{c}2.1689\\2.0237 \end{array}   \\
\hline
\end{array}
\end{array}$

\caption{The two largest eigenvalues of $P(3,p_2,p_3)$ with a $4$ decimal place accuracy.}
\label{tab_eig_P3}
\end{table}
\end{center}

\begin{center}
\begin{table}[htbp]
$\begin{array}{cc}

\includegraphics[scale=0.7]{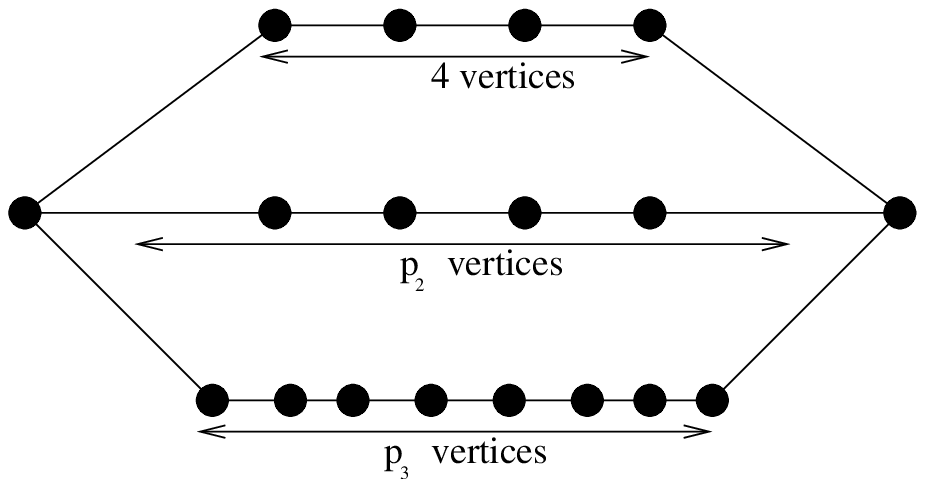}

&
\begin{array}{|c||c|c|}
\hline
\textrm{\backslashbox{$p_3$}{$p_2$}}&  4  &  6    \\
\hline
\hline
4       &  \begin{array}{c}2.1987\\1.9122 \end{array}   &   \begin{array}{c}2.1853\\1.9666 \end{array}    \\
\hline
6       &  \begin{array}{c}2.1853\\1.9666 \end{array}    &   \begin{array}{c}2.1723\\2.0102 \end{array}   \\
\hline
8       &   \begin{array}{c}2.1790\\1.9922 \end{array}   &      \begin{array}{c}2.1660\\2.0300 \end{array}  \\
\hline
10       &   \begin{array}{c}2.1762\\2.0058 \end{array}    &    \begin{array}{c}2.1631\\2.0401 \end{array}     \\
\hline
\end{array}
\end{array}$

\caption{The two largest eigenvalues of $P(4,p_2,p_3)$ with a $4$ decimal place accuracy.}
\label{tab_eig_P4}
\end{table}
\end{center}

\begin{center}
\begin{table}[htbp]
$\begin{array}{cc}

\begin{array}{cc}
\textrm{Graph}  & \textrm{Eigenvalues} \\

\includegraphics[scale=0.3]{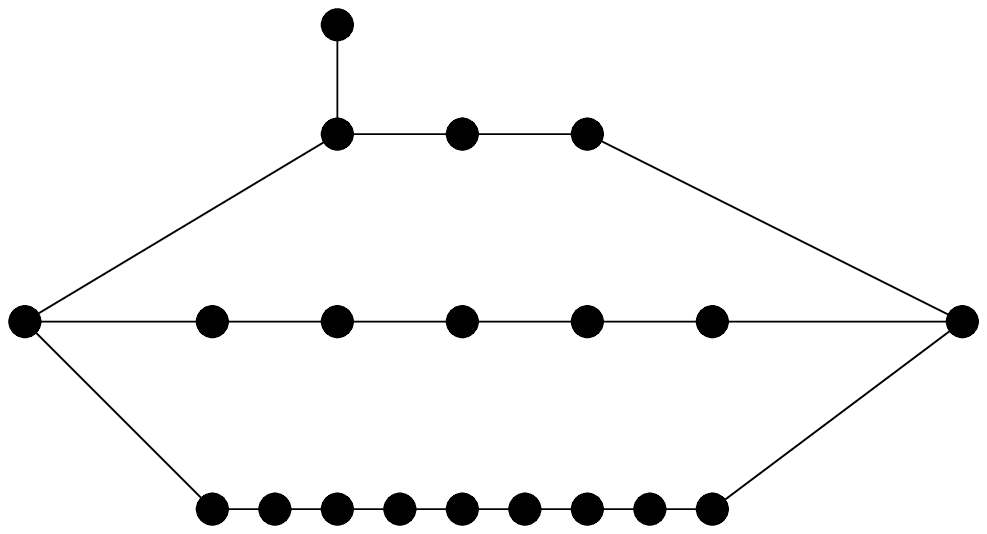} & \begin{array}{c}2.2346\\2.0117 \end{array} \\
\includegraphics[scale=0.3]{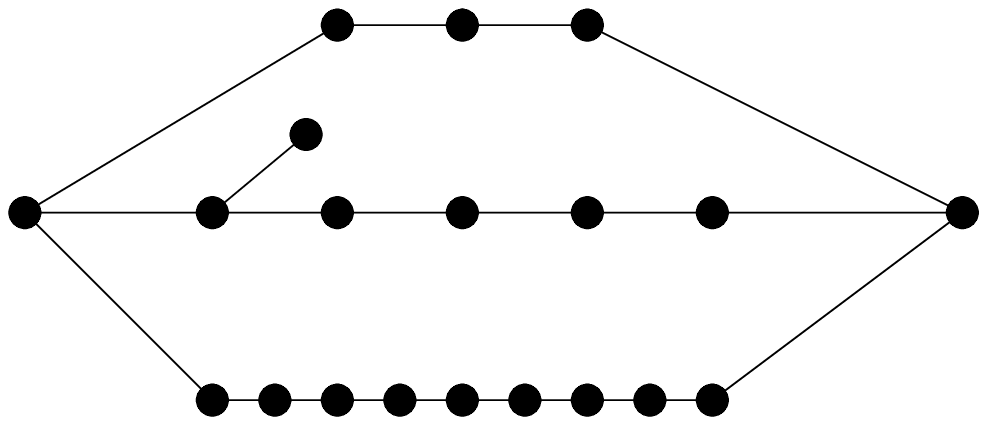} & \begin{array}{c}2.2288\\2.0287 \end{array} \\
\includegraphics[scale=0.3]{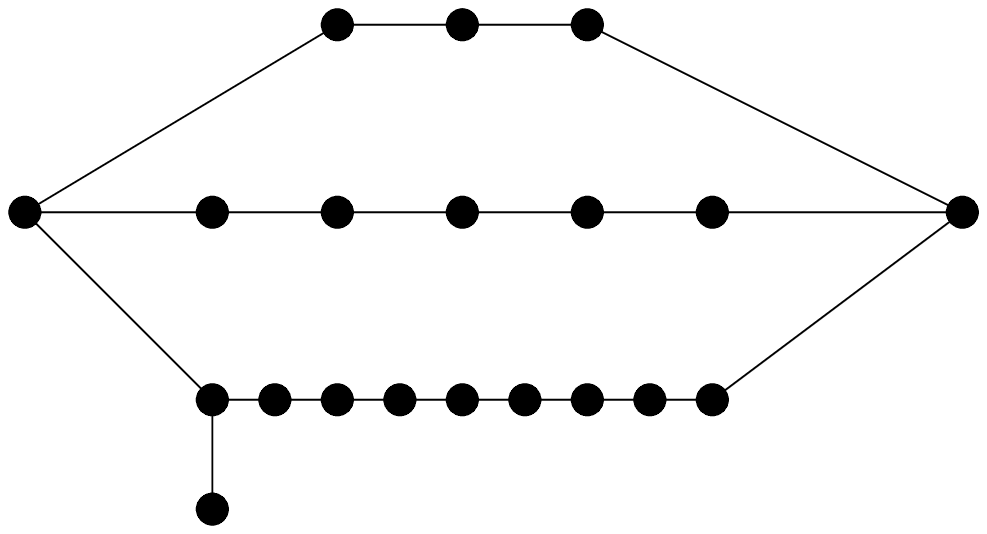} & \begin{array}{c}2.2249\\2.0433 \end{array} \\

\includegraphics[scale=0.3]{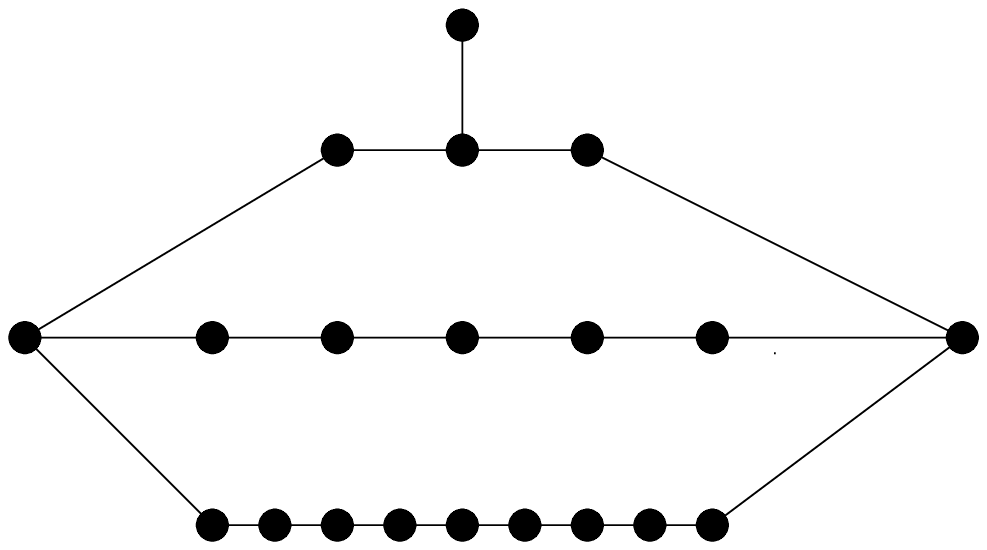} & \begin{array}{c}2.2247\\1.9890 \end{array} \\
\includegraphics[scale=0.3]{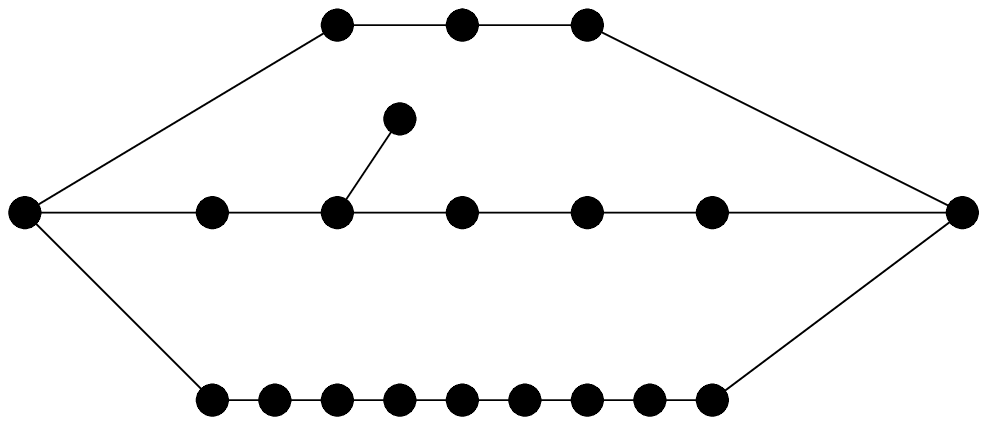} & \begin{array}{c}2.2109\\2.0043 \end{array} \\
\includegraphics[scale=0.3]{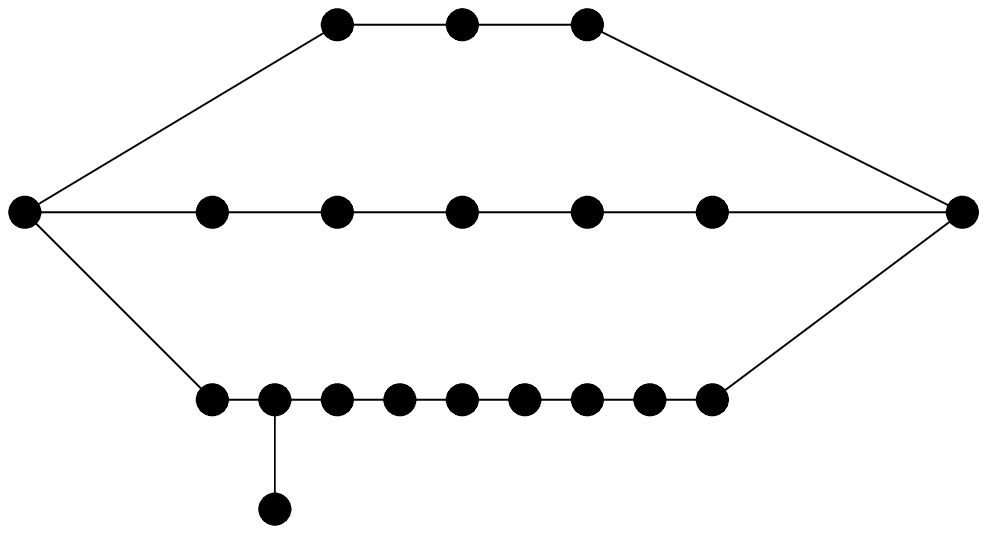} & \begin{array}{c}2.2026\\2.0352 \end{array} \\
\end{array}
&
\begin{array}{cc}
\textrm{Graph}  & \textrm{Eigenvalues} \\

\includegraphics[scale=0.3]{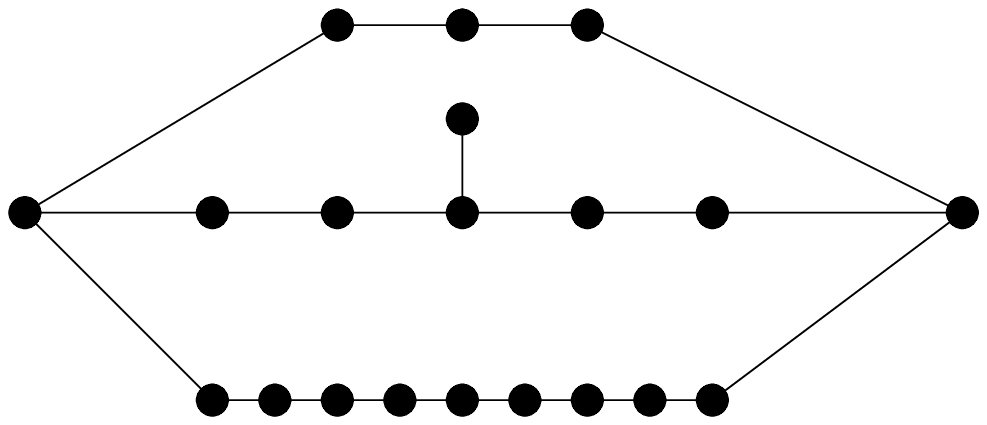} & \begin{array}{c}2.2062\\1.9890\end{array} \\
\includegraphics[scale=0.3]{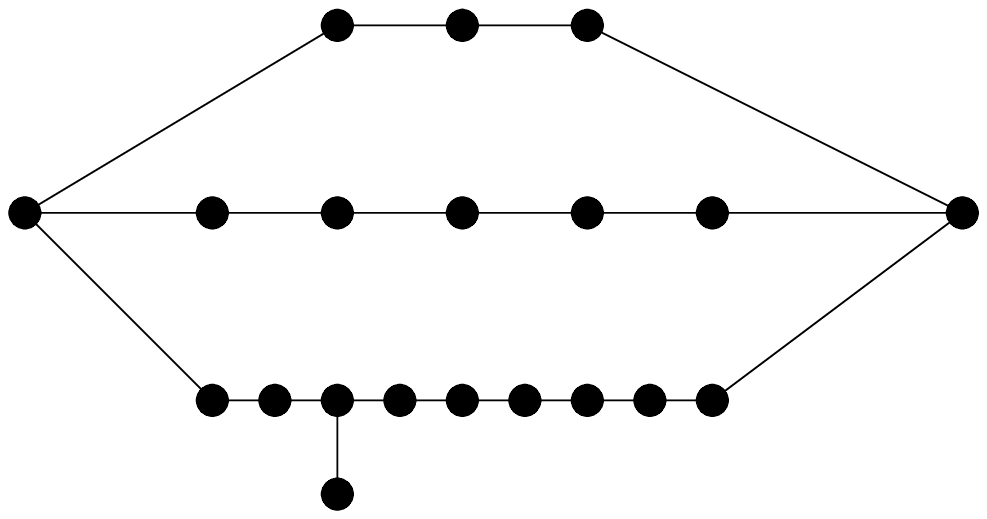} & \begin{array}{c}2.1925\\2.0253 \end{array} \\
\includegraphics[scale=0.3]{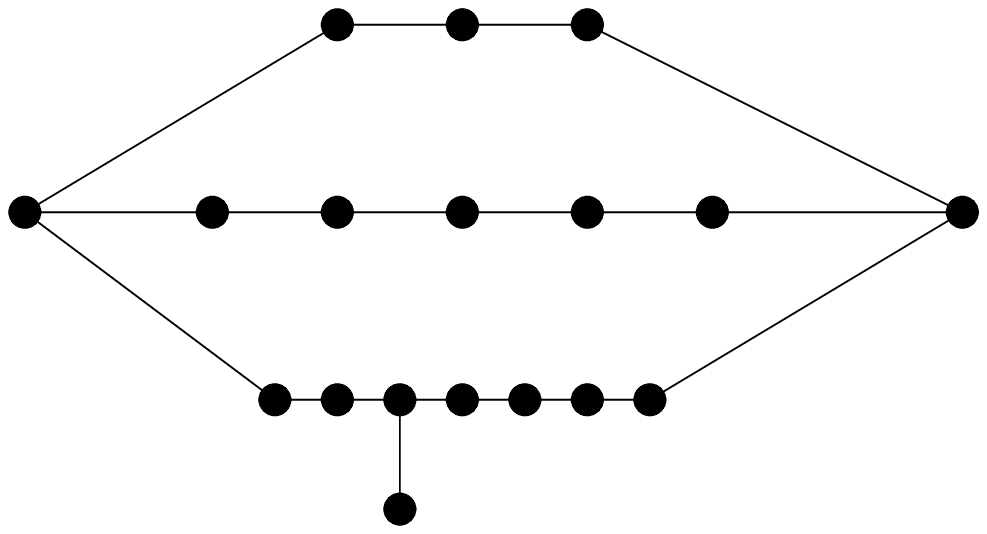} & \begin{array}{c}2.1999\\1.9909 \end{array} \\
\includegraphics[scale=0.3]{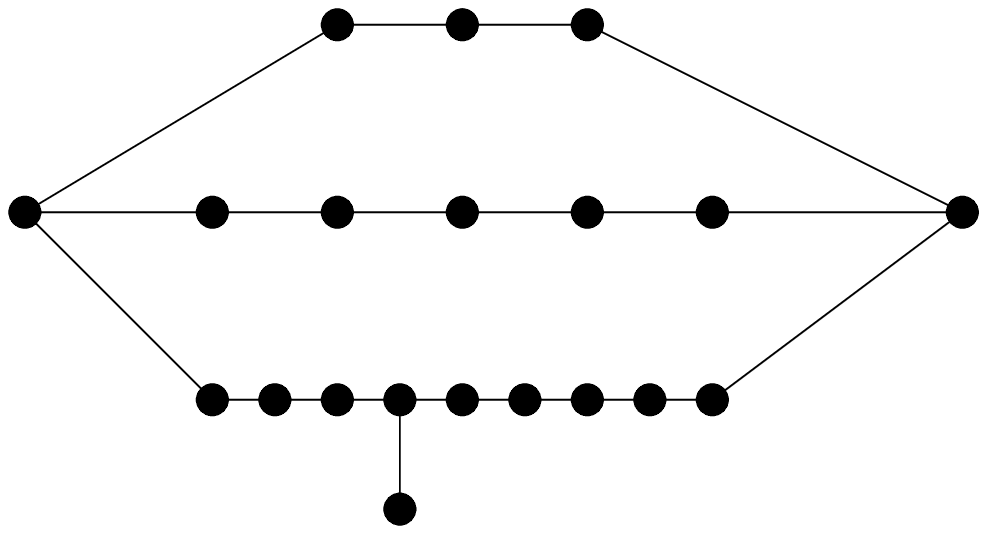} & \begin{array}{c}2.1882\\2.0126 \end{array} \\
\includegraphics[scale=0.3]{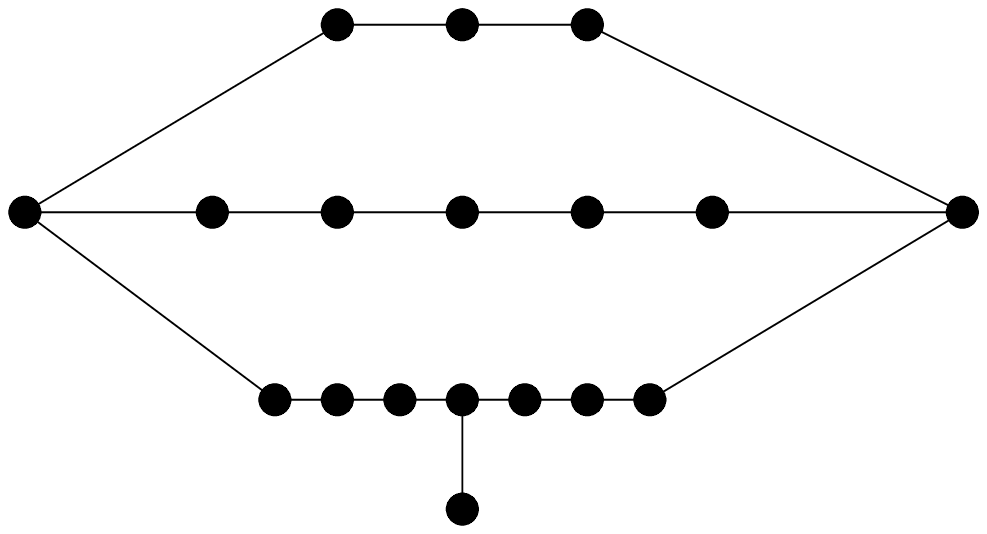} & \begin{array}{c}2.1976\\1.9696 \end{array} \\
\includegraphics[scale=0.3]{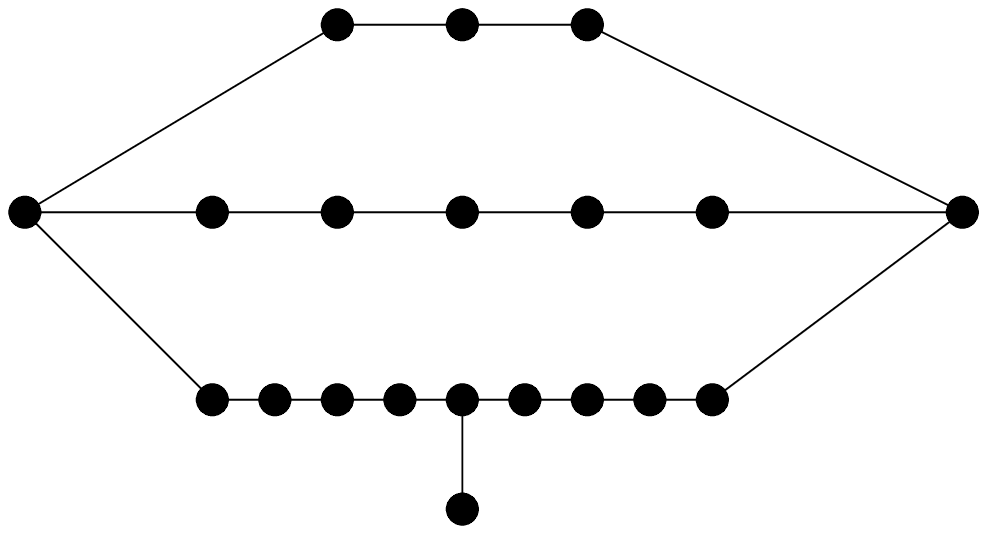} & \begin{array}{c}2.1870\\2.0000 \end{array} \\
\end{array}

\end{array}$

\caption{The two largest eigenvalues of some graphs $H$ with a $4$ decimal place accuracy. Note that the spectral radius increases when  the number of vertices between two vertices of degree $3$ decreases (theorem \ref{homeo}).}
\label{tab_P+v}
\end{table}
\end{center}

\begin{center}
\begin{table}[htbp]
$\begin{array}{cc}

\begin{array}{cc}
\textrm{Graph}  & \textrm{Eigenvalues} \\
\includegraphics[scale=0.3]{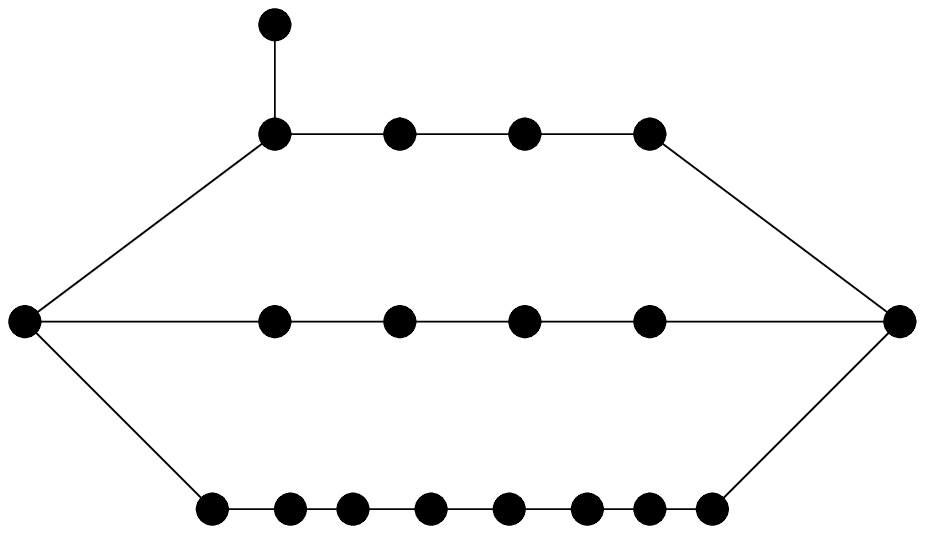} & \begin{array}{c}2.2292 \\2.0262  \end{array} \\
\includegraphics[scale=0.3]{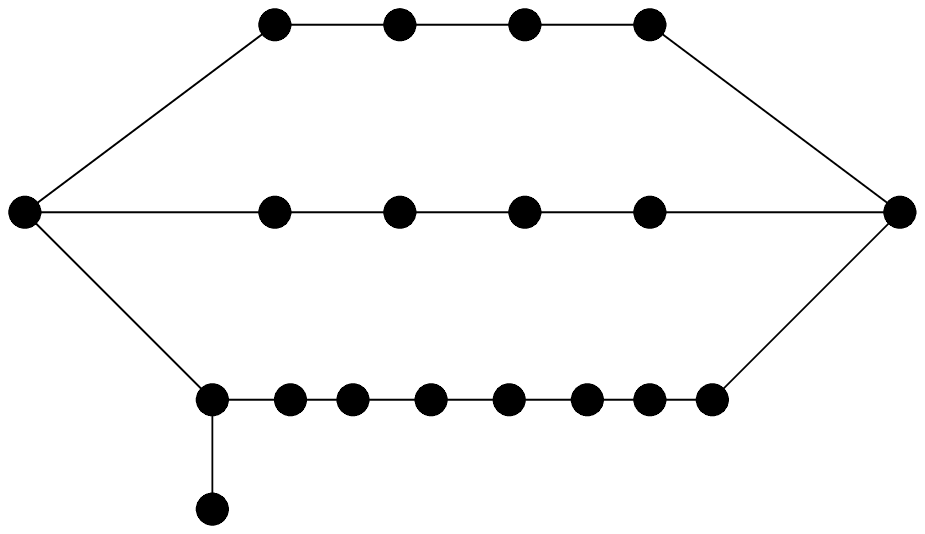} & \begin{array}{c}2.2236 \\2.0460  \end{array} \\
\includegraphics[scale=0.3]{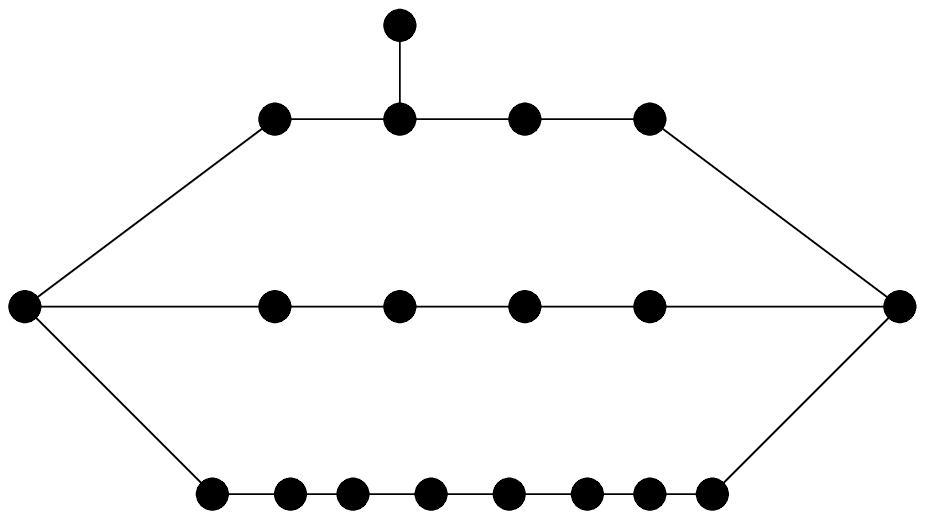} & \begin{array}{c}2.2143\\1.9977\end{array} \\
\includegraphics[scale=0.3]{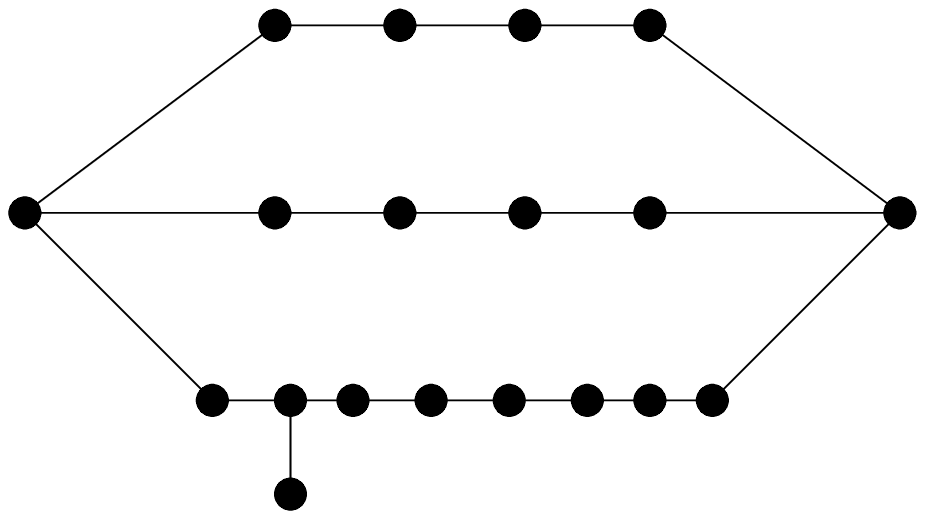} & \begin{array}{c}2.2020\\2.0341 \end{array} \\

\end{array}
&
\begin{array}{cc}
\textrm{Graph}  & \textrm{Eigenvalues} \\
\includegraphics[scale=0.3]{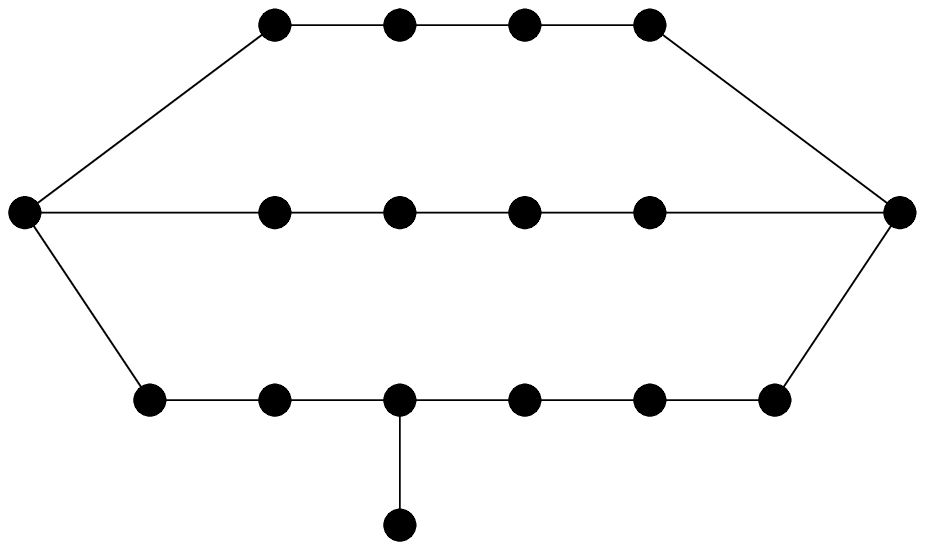} & \begin{array}{c}2.2037\\1.9737 \end{array} \\
\includegraphics[scale=0.3]{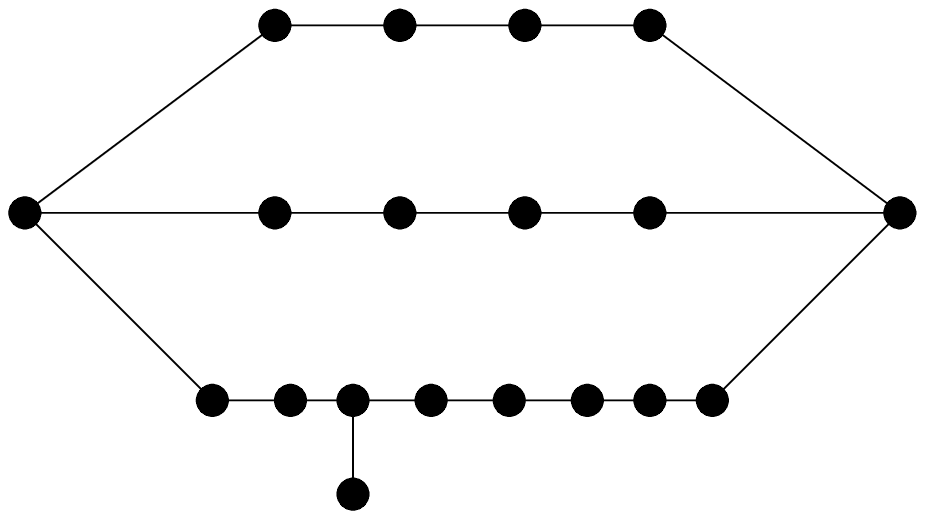} & \begin{array}{c}2.1925\\2.0191 \end{array} \\
\includegraphics[scale=0.3]{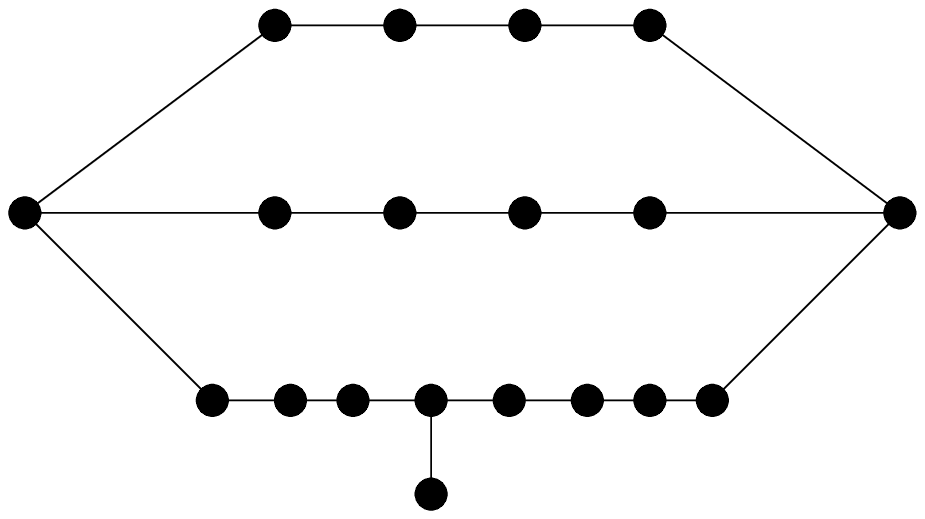} & \begin{array}{c}2.1889\\2.0000 \end{array} \\

\end{array}

\end{array}$

\caption{The two largest eigenvalues of some graphs $H$ with a $4$ decimal place accuracy. Note that the spectral radius increases when  the number of vertices between two vertices of degre $3$ decreases (theorem \ref{homeo}).}
\label{tab_P+v2}
\end{table}
\end{center}

\begin{center}
\begin{table}[htbp]
$\begin{array}{c}

\includegraphics[scale=0.7]{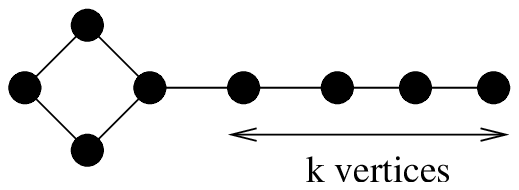}

\\

\begin{array}{c|c|c|c|c|c|c|c|}
k                   &  1  &  2  &  3   & 4  &  5  &  6  &  7   \\
\hline
\lambda_1(L(4,k))     &  2.1358   & 2.1753 & 2.1889 & 2.1940  & 2.1960  & 2.1968  &  2.1971       \\
\hline
\end{array}

\\\\

\begin{array}{c|c|c|c|c|c|c|c|}
k                  &  8  &  9  &  10  &  11   & 12  &  13  &  14     \\
\hline
\lambda_1(L(4,k))    &2.1973 &  2.1973   &  2.1974 & 2.1974 & 2.1974  & 2.1974  & 2.1974   \\
\hline
\end{array}

\end{array}
$

\caption{Spectral radius of $L(4,k)$ with a $4$ decimal place accuracy.}
\label{spect_L4k}
\end{table}
\end{center}

\begin{center}
\begin{table}[htbp]
$\begin{array}{cc}
\begin{array}{cc}
\textrm{Graph}  & \textrm{Spectral radius} \\
\includegraphics[scale=0.25]{except.eps} & 2.1856 \\
\includegraphics[scale=0.25]{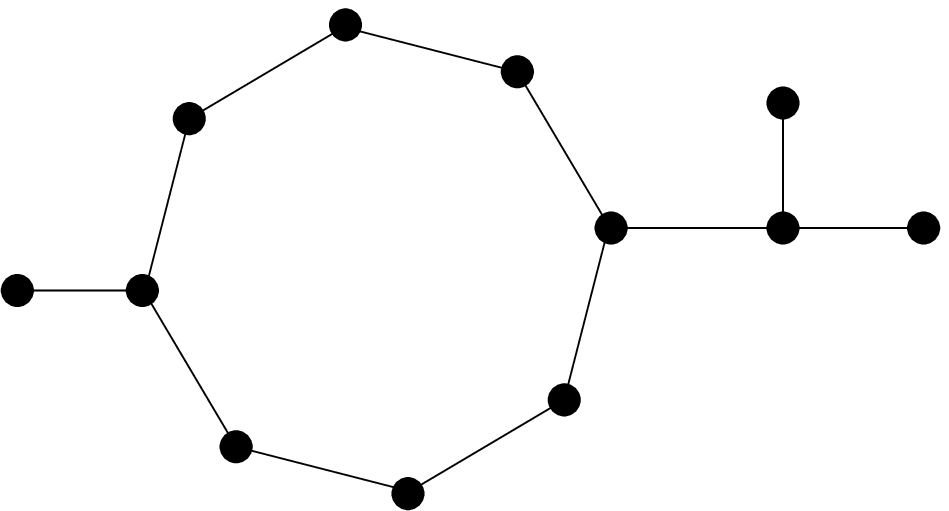} & 2.2005\\
\includegraphics[scale=0.25]{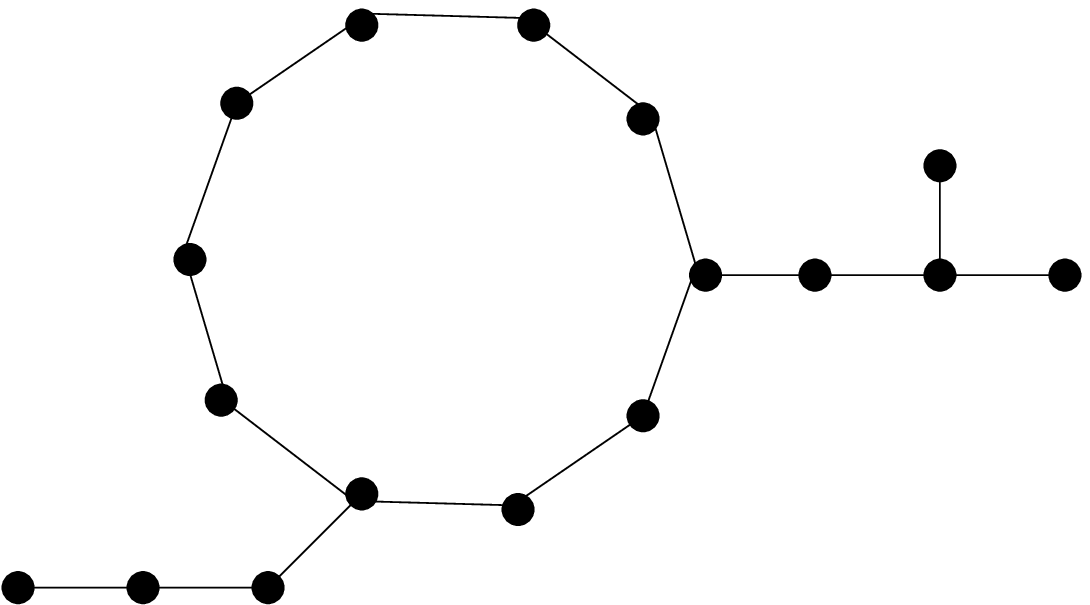} & 2.1927 \\
\includegraphics[scale=0.25]{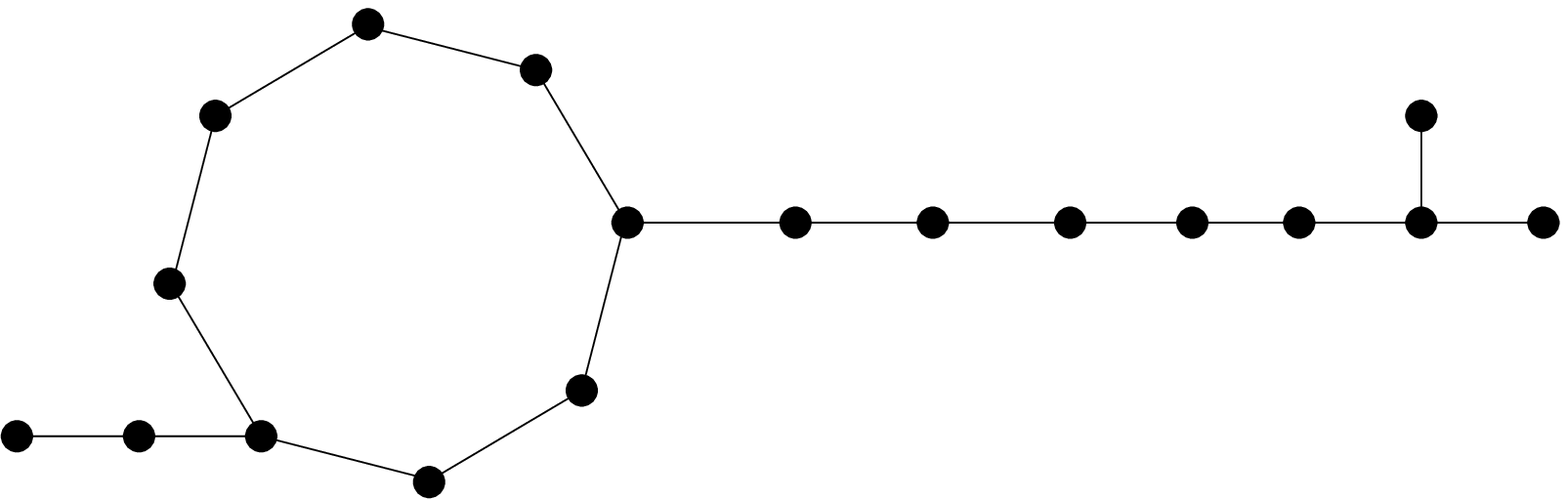} & 2.1922 \\
\includegraphics[scale=0.25]{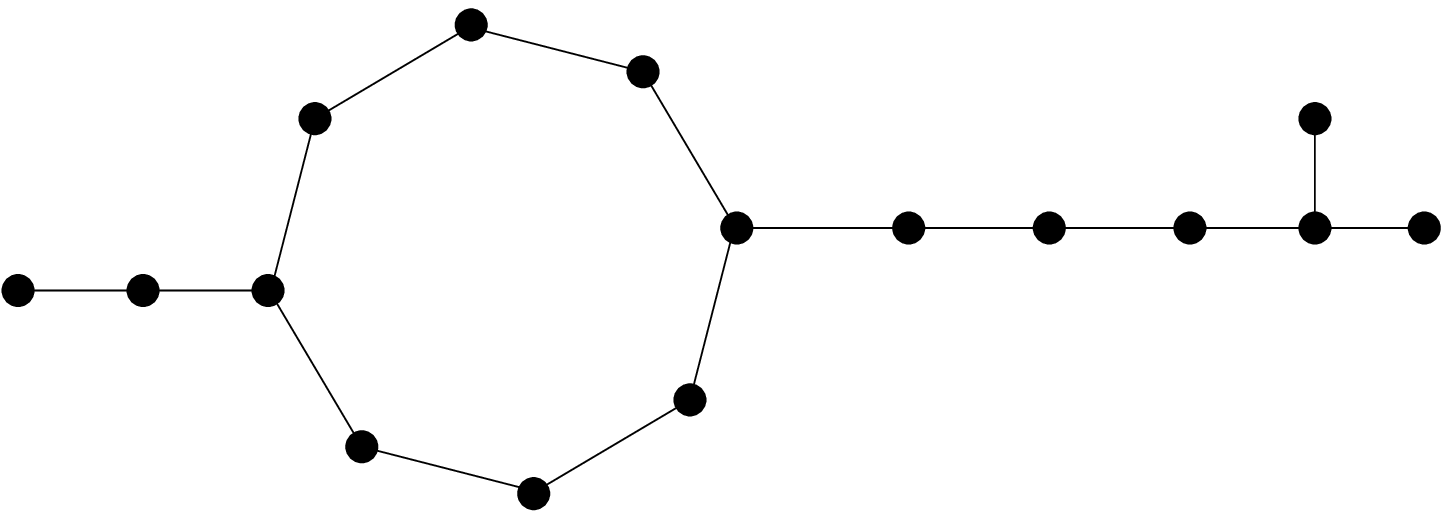} & 2.1894 \\
\includegraphics[scale=0.25]{C8_3.eps} & 2.2025 \\
\end{array}
&
\begin{array}{cc}
\textrm{Graph}  & \textrm{Spectral radius} \\
\includegraphics[scale=0.25]{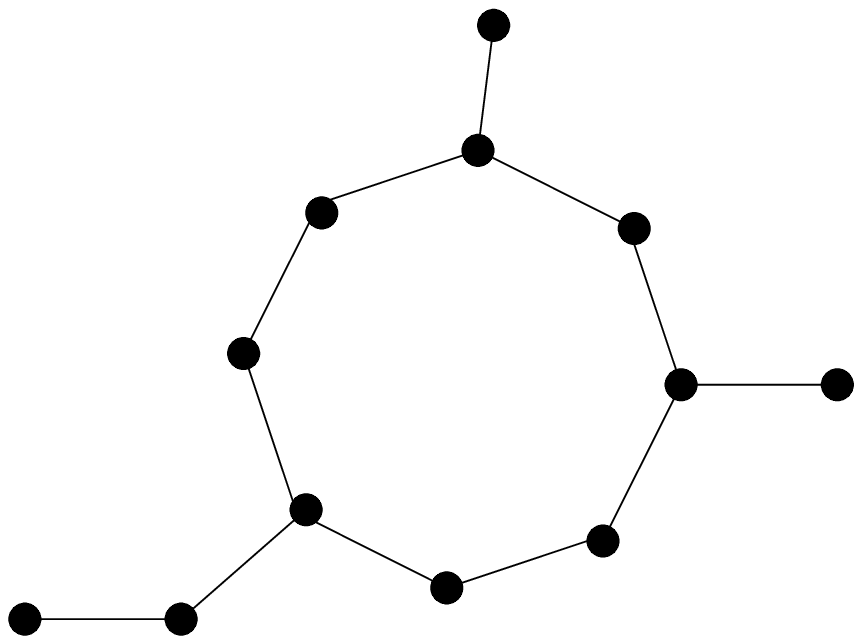} & 2.2047 \\
\includegraphics[scale=0.25]{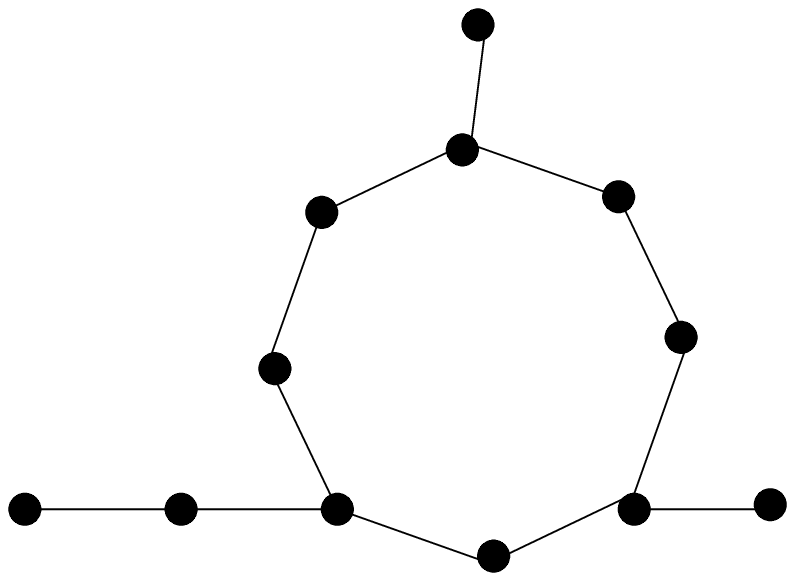} & 2.2075 \\
\includegraphics[scale=0.25]{C10_113.eps} & 2.1987 \\
\includegraphics[scale=0.3]{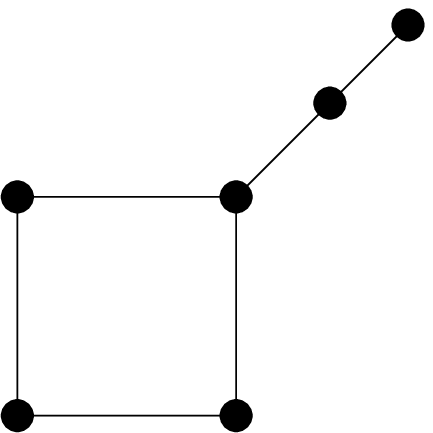} & 2.1753\\
\includegraphics[scale=0.3]{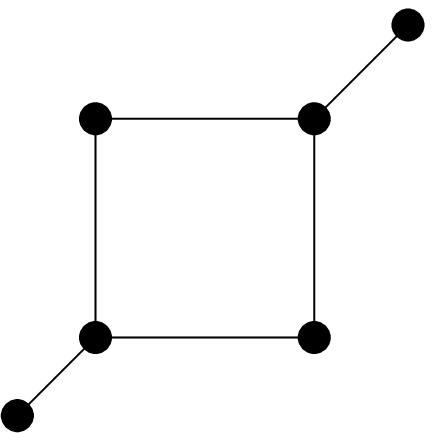} & 2.2361 \\
\includegraphics[scale=0.3]{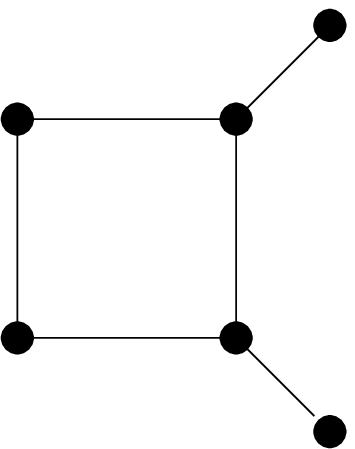} & 2.2470 \\
\includegraphics[scale=0.3]{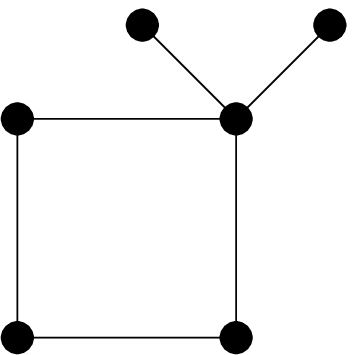} & 2.2883 \\
\end{array}

\end{array}$

\caption{Spectral radius of some unicyclic graphs with a $4$ decimal place accuracy.}
\label{spect_special}
\end{table}
\end{center}

\newpage
\bibliographystyle{plain}
\bibliography{./MaBiblio}

\newpage
\tableofcontents
\end{document}